\pdfoutput=1

\documentclass[preprint,10pt]{elsarticle}

\usepackage{fullpage}
\usepackage[colorlinks=true]{hyperref}
\usepackage{float}

\usepackage{amsmath,amssymb,amsfonts,amsthm}
\theoremstyle{definition}

\theoremstyle{lemma}

\theoremstyle{theorem}
\newtheorem{theorem}{Theorem}
\theoremstyle{assumption}

\usepackage[titletoc,toc,title]{appendix}
\usepackage{array} 
\usepackage{listings}
\usepackage{mathtools}
\usepackage{pdfpages}
\usepackage[textsize=footnotesize,color=green]{todonotes}
\usepackage{bm}
\usepackage{bbm}
\usepackage{mathabx}

\usepackage{tikz}
\usetikzlibrary{shapes,arrows}
\usepackage[normalem]{ulem}
\usepackage{hhline}

\usepackage{graphicx}
\usepackage{subfig}
\usepackage{color}

\usepackage{pgfplots}
\usepackage{pgfplotstable}
\definecolor{markercolor}{RGB}{124.9, 255, 160.65}
\pgfplotsset{
compat=1.3,
width=10cm,
tick label style={font=\small},
label style={font=\small},
legend style={font=\small}
}

\usetikzlibrary{calc}
\usetikzlibrary{intersections} 

\newcommand{\logLogSlopeTriangleFlip}[5]
{

    \pgfplotsextra
    {
        \pgfkeysgetvalue{/pgfplots/xmin}{\xmin}
        \pgfkeysgetvalue{/pgfplots/xmax}{\xmax}
        \pgfkeysgetvalue{/pgfplots/ymin}{\ymin}
        \pgfkeysgetvalue{/pgfplots/ymax}{\ymax}

        \pgfmathsetmacro{\xBrel}{#1-#2}
        \pgfmathsetmacro{\yBrel}{#3}
        \pgfmathsetmacro{\xCrel}{#1}

        \pgfmathsetmacro{\lnxB}{\xmin*(1-(#1-#2))+\xmax*(#1-#2)} 
        \pgfmathsetmacro{\lnxA}{\xmin*(1-#1)+\xmax*#1} 
        \pgfmathsetmacro{\lnyA}{\ymin*(1-#3)+\ymax*#3} 
        \pgfmathsetmacro{\lnyC}{\lnyA+#4*(\lnxA-\lnxB)}
        \pgfmathsetmacro{\yCrel}{\lnyC-\ymin)/(\ymax-\ymin)} 

	\pgfmathsetmacro{\xArel}{\xBrel}
        \pgfmathsetmacro{\yArel}{\yCrel}

        \coordinate (A) at (rel axis cs:\xArel,\yArel);
        \coordinate (B) at (rel axis cs:\xBrel,\yBrel);
        \coordinate (C) at (rel axis cs:\xCrel,\yCrel);

        \draw[#5]   (A)-- node[pos=0.5,anchor=east] {#4}
                    (B)-- 
                    (C)-- node[pos=0.5,anchor=south] {}
                    cycle;
    }
}

\usepackage{stmaryrd}
\newcommand{\myfrac}[2]{\displaystyle \frac{#1}{#2}}

\newcommand{\LRp}[1]{\left( #1 \right)}

\newcommand{\LRc}[1]{\left\{ #1 \right\}}

\newcommand{\jump}[1] {\ensuremath{\llbracket#1\rrbracket}}
\newcommand{\avg}[1] {\ensuremath{\LRc{\!\{#1\}\!}}}

\newcommand{\eval}[2][\right]{\relax
  \ifx#1\right\relax \left.\fi#2#1\rvert}

\def\etal{{\it et al.~}}

\newcolumntype{C}[1]{>{\centering\let\newline\\\arraybackslash\hspace{0pt}}m{#1}}

\makeatletter
\renewcommand\d[1]{\mspace{6mu}\mathrm{d}#1\@ifnextchar\d{\mspace{-3mu}}{}}
\makeatother

\graphicspath{{./figs/}}

\begin{document}


\begin{frontmatter}
\title{A high order discontinuous Galerkin method for the symmetric form of the anisotropic viscoelastic wave equation}

\author{Khemraj Shukla \corref{cor}}
\ead{rajexplo@gmail.com}

\author{Jesse Chan \corref{}}
\ead{jesse.chan@rice.edu}

\author{Maarten V. de Hoop\corref{}}
\ead{mdehoop@rice.edu}

\cortext[cor]{corresponding author}

\address{Department of Computational and Applied Mathematics, Rice University, 6100 Main St, Houston, TX, 77005}
\begin{abstract}
Wave propagation in real media is affected by various non-trivial physical phenomena, e.g., anisotropy, an-elasticity and dissipation. 
Assumptions on the stress-strain relationship are an integral part of seismic modeling and determine the deformation and relaxation of the medium. Stress-strain relationships based on simplified rheologies will incorrectly predict seismic amplitudes, which are used for quantitative reservoir characterization. Constitutive equations for the rheological model include the generalized Hooke's law and Boltzmann's superposition principal with dissipation models based on standard linear solids or a Zener approximation.  

In this work, we introduce a high-order discontinuous Galerkin finite element method for wave equation in inhomogeneous and anisotropic dissipative medium. 
This method is based on a new symmetric treatment of the anisotropic viscoelastic terms, as well as an appropriate memory variable treatment of the stress-strain convolution terms. Together, these result in a symmetric system of first order linear hyperbolic partial differential equations. 
The accuracy of the proposed numerical scheme is proven and verified using convergence studies against analytical plane wave solutions and analytical solutions of viscoelastic wave equation. Computational experiments are shown for various combinations of homogeneous and heterogeneous viscoelastic media in two and three dimensions.
\end{abstract}
\end{frontmatter}


\section{Introduction}
Numerical simulations of seismic wave propagation are essential for various imaging problems arising at different scales. At a global scale, seismic waves, traveling thorough the entire Earth allow geophysicists to infer properties of the Earth interior. At a macro-scale, seismic wave propagation can be used to image and characterize oil and gas reservoirs. On a micro or laboratory scale, seismic waves play a major role in studying the micro-structure of materials. To simulate wave propagation accurately, the input model should be able to accommodate arbitrary variations of petrophysical and lithological properties, as they play an important role, particularly in the targets of exploration geophysics, i.e., reservoir rocks. To study reservoir monitoring and evaluation of rock properties in a laboratory setting, lithological and reservoir properties become more important. Reservoir rocks such as cracked limestones can show effective anisotropy in the seismic band. Furthermore, fluid-filled cracked rocks and porous sandstones show considerable attenuation properties. Experimental work also shows that anisotropy effects of attenuation are more pronounced than anisotropic elastic effects \cite{hosten1987, arts1992}. Thus, a realistic rheology should be able to model  anisotropic attenuation characteristics.

Various dissipation mechanisms  (e.g. Kelvin-Voigt, Maxwell, Zener \cite{carcione2014}) can be modeled by a viscoelastic constitutive relation. Attenuation of energy is caused by a large variety of dissipation mechanisms. it is difficult, if not impossible, to build a general microstructure theory incorporating all these mechanisms. Modeling of dissipation in isotopic media requires two relaxation functions. These two relaxation functions are enough to describe anelastic characteristics of body waves because these modes decouple in a homogeneous medium with isotropic attenuation \cite{carcione1990}. In contrast, in anisotropic media, one has to decide the time (or frequency) dependence of 21 stiffness parameters. However, Mehrabadi and Cowin \cite{mehrabadi1990} and Helbig \cite{helbig1996} show that only six of the 21 stiffnesses parameters have an intrinsic physical meaning. 

Motivated by the works of Mehrabadi and Cowin \cite{mehrabadi1990} and Helbig \cite{helbig1996}, Carcione formulated a constitutive model and wave equations for linear viscoelastic anisotropic media \cite{carcione1995} . In three-dimensional anisotropic media, careful attention is required for modeling properties of the shear
modes, since the relaxation of a medium can be different for slow and fast waves . Carcione used a relaxation function to model the anelastic properties of the quasi-dilatational mode, whereas three relaxation functions are used to control the relaxation of the medium due to the shear waves along preferred directions. In this study, we use the constitutive model proposed by Caricone \cite{carcione1990} and pair it with equations of motion  described by Newton's second law of motion.  This results in a system of first order hyperbolic PDEs with stress-velocity field variables.  

Numerical simulations of seismic wave propagation solve the system of hyperbolic partial differential equations using various  numerical methods such as finite-differences, finite volumes and finite element methods. An overview of these methods is given by \cite{moczo2014, igel2017}. The most popular and simple method is the finite-difference (FD) method, and its application to the elastic wave equation has been studied by many researchers \cite{virieux1986, marfurt1984, moczo2007, bohlen2002, levander1988}. A detailed analysis of finite-difference methods is given in \cite{strikwerda2004}. Although the numerical representation of FD method is very simple,  it often comes with significant numerical dispersion, especially in the modeling of surface waves \cite{marfurt1984}. Additionally, the implementation of boundary conditions can require special treatment \cite{drainville2019}. FD methods are also difficult to apply to irregular geometries with out experiencing ``staircase effects" \cite{drainville2019}. To circumvent the effect of numerical dispersion and achieve high order accuracy, pseudo-spectral methods were first used by Tessmer and Kosloff \cite{tessmer}. The pseudo-spectral method uses global basis function for approximation of the solutions (e.g., Fourier or Chebyshev). The pseudo-spectral method requires few grid points per wavelength and produces a high order solution with less numerical dispersion. However, the choice of the global basis functions restricts the pseudo-spectral methods to smooth models, as it is difficult to represent materials with the discontinuities or sharp contrast. This can be addressed somewhat using domain decomposition, where  different meshes are used to represent the different domains. For example, Carcione \cite{carcione1991}  used Fourier basis functions along directions with smoothly varying materials properties and Chebyshev basis functions in directions with sharply varying medium properties.

The finite element method (FEM) method discretizes the domain using elements (e.g., triangles and quadrilaterals in 2D and tetrahedra and hexahedra in 3D). Since time-domain wave propagation is described by a hyperbolic system of partial differential equations, an explicit time integration can efficiently be applied. However, finite-element methods, when coupled with an explicit time integrator, require the inversion of a global mass matrix unless special techniques (such as diagonal mass lumping) are applied. Finite elements for elastic wave propagation were studied by Marfurt \cite{marfurt1984} and Bao \etal \cite{bao1998}.  In these works, FEM was shown to accurately represent sharp material properties and irregular geometries. However, the solution on each element is approximated using a low order polynomial, which results in significant numerical dispersion. In order to recover a more accurate solution, high order elements are required which results into large matrices to be inverted at each time step. To exploit the spectral properties in the finite element method Patera \cite{patera1984} proposed the spectral element method (SEM) to solve fluid flow problems. Subsequently, the SEM was successfully implemented by Seriani \etal \cite{seriani1994spectral} to solve the acoustic wave equation in a heterogeneous medium. Komatitsch and Vilotte \cite{komatitsch1998} used SEM to solve the elastic wave equation in a heterogeneous medium, described by a system of the second order PDEs. A detailed review of seismic modeling  is presented Carcione \etal \cite{carcione2002}.

In the present study, we introduce a high-order numerical scheme based on the discontinuous Galerkin method to solve the 3D viscoelastic  wave equation on unstructured tetrahedral meshes. High order methods provide one avenue towards improving fidelity in numerical simulations while maintaining reasonable computational costs, and methods which can accommodate unstructured meshes are desirable for problems with complex geometries. Among such methods, high order discontinuous Galerkin (DG) methods are particularly well-suited to the solution of time-dependent hyperbolic problems on modern computing architectures \cite{hesthaven2007, klockner2009}. The accuracy of high order methods can be attributed in part to their low numerical dissipation and dispersion compared to low order schemes \cite{ainsworth2004}. This accuracy has made them advantageous for the simulation of electro-magnetic and elastic wave propagation \cite{hesthaven2007, wilcox2010}. Spectral element methods avoid the inversion of the global mass matrix for quadrilateral and hexahedral elements by choosing nodal basis functions, which are discretely orthogonal with respect to an under-integrated $L^2$ inner product and result in a diagonal mass matrix \cite{komatitsch1998}. In contrast, high order DG methods produce block diagonal mass matrices, which are locally invertible. High order DG methods are often used for seismic simulation (elastic approximation) in combination with simplicial meshes \cite{kaser2007,de2007,ye2016}. 

DG methods impose inter-element continuity of approximate solutions between elements weakly through a numerical flux, of which the upwind flux  (solution of a Riemann problem)is more common. Käser \etal \cite{kaser2007} solved the 3D isotropic viscoelastic wave equation in a strain-velocity formulation using a local space-time DG method with an upwind flux by solving the exact Riemann problem on inter-element boundaries. In another study, Lambrecht \etal \cite{lambrecht2018} used a nodal DG method to solve the isotropic viscoelasltic wave equation using the same formulation proposed by Käser et al. \cite{kaser2007}. The solution of the Riemann problem requires diagonalization of Jacobian matrices into polarized waves constituents, which is a computationally intensive process for the viscoelastic system and does not extend naturally to anisotropic materials. Ye \etal\cite{ye2016} completely avoids the process of diagonalization for the coupled acoustic-elastic wave equation by using a penalty flux based on natural boundary conditions.  In this study, we use a similar energy-stable penalty flux for the anisotropic viscoelastic wave equations. 

The main new contributions of this paper are a new symmetric form of the anisotropic viscoelastic wave equation and its discretization using a high order DG method using penalty fluxes. The outline of the paper is as follows: Section \ref{sec:eq} will discuss the system of equations describing the viscoelastic wave equation. Section \ref{sec:ESdG} presents an energy stable formulation for the symmetric hyperbolic form of the viscoelastic wave equations. Finally, numerical results in Section \ref{numex} demonstrate the accuracy of this method for several problems in linear anisotropic viscoelasticity.

\section{Constitutive Relations}
\label{sec:eq}
\subsection{Notation and convention}
Let $f$ and $g$ be scalar and time-dependent functions. The Riemann convolution of $f$ with $g$ is defined as 

\begin{align} \label{eq1}
f * g =\begin{dcases} \int_0^t f(\tau) g(t-\tau)~d\tau: & \qquad t \ge 0  \\
                      \qquad \qquad \qquad \qquad  0: & \qquad   t <0,
\end{dcases}
\end{align}
where $t$ is the time variable.
Hooke's law is expressed in 3-D or 6-D space depending on whether stress or strain are considered as tensors or vectors, and  the convolution shown by (\ref{eq1}) is easily extended to include vectors and tensors as follows

\begin{align} \label{eq2}
\bm{\Psi} * \bm{A} = \int_0^t \bm{\Psi}(\tau) \cdot \bm{A}(t - \tau)~d\tau
\end{align}
where $\bm{A}$ and $\bm{\Psi}$ are $6 \times 6$ stress (strain) tensor and $\bm{\Psi}$ are relaxation matrices. 

If $f$ and $g$ are of Heaviside type, the Boltzman operation \cite{leitman1973} defines the time derivative of the convolution between $f$ and $g$ as 
\begin{align} \label{eq3}
f*\partial_t g=f \odot g = \mathring{f}g + \left( \dot{f} H \right)* g,
\end{align} 
where $\mathring{f} = f(t=0) = f(t=0^+)$ and $H(t)$ is the step function.
Sub-indices $i,~j, ~k$, and $m$ take values from 1 to 3 and correspond to the three Cartesian coordinates $x,~y$ and $z$.

\subsection{Boltzmann law}
The general constitutive relation for an anisotropic and linear viscoelastic medium can be expressed as \cite{carcione1990}
\begin{align} \label{eq4}
\bm{\sigma} = \dot{\bm{\Psi}} * \bm{\epsilon}
\end{align}
where $\bm{\sigma} =[\sigma_{11}, \sigma_{22}, \sigma_{22}, \sigma_{23}, \sigma_{13}, \sigma_{12}]^T$, $\bm{\epsilon} = [\epsilon_{11}, \epsilon_{22}, \epsilon_{33}, \epsilon_{23}, \epsilon_{13}, \epsilon_{12}]^T$  are stress and strain vectors with $\gamma_{ij}=2\epsilon_{ij}$, and $\bm{\Psi}$ is the symmetric relaxation matrix. 

The stress-strain relation in ({\ref{eq4}) is called the Boltzmann law and can be expressed as more explicitly by using the Einstein notation of summation over repeated indices
\begin{align}\label{eq5}
\sigma_I = \psi_{IJ}*\dot{\epsilon}_J, \qquad  I,J=1,...,6
\end{align}
To model the anelastic properties of shear waves, Carcione \cite{carcione1995} used one relaxation function for quasi dilatational mode and three relaxation functions for shear modes and expressed relaxation matrix $\bm{\Psi}$ as 
 \begin{align} \label{eq6}
 \bm{\Psi} = \left[ \begin{array}{cccccc}
 \psi_{11} & \psi_{12} & \psi_{13} & c_{14} & c_{15} & c_{16}\\ 
 ~&  \psi_{22} & \psi_{23} & c_{24} & c_{25} & c_{26} \\
 ~&~& \psi_{33} & c_{34} & c_{35} & c_{36}\\
 ~&~&~ & c_{44}\chi_2 & c_{45} & c_{46} \\
 ~&~&~&~& c_{55} \chi_{3} & c_{56} \\
 ~&~&~&~&~&c_{66} \chi_4
 \end{array}\right] H (t),
 \end{align}
where 
\begin{align}
\psi_{i(I)}&=c_{I(I)} -D + K \chi_1 + \myfrac{4}{3} G\chi_\delta \qquad \text{for~} I=1,2,3, \label{eq7}\\
\psi_{IJ}& =c_{IJ} - D + 2G + K\chi_1 -\myfrac{2}{3}G\chi_\delta\qquad \text{for~} I,J =1,2,3; I \neq J. \label{eq8}
\end{align}
The $C_{IJ}$ for $I,J = 1,....,6$ are the high-frequency limit (relaxed) elasticities i.e. $t  \rightarrow 0; \omega \rightarrow \infty$, and
\begin{align}\label{eq9}
K=D-\myfrac{4}{3}G,
\end{align} 
where 
\[
D=\myfrac{1}{3} (c_{11} + c_{22} + c_{33}), \qquad G=\myfrac{1}{3}(c_{44} + c_{55} + c_{66}).
\]
 The $\chi_\nu$ are dimensionless relaxation function with index $\nu=1$ representing the quasi-dilatational mode and indices $\nu=2,3,4$ corresponds to shear waves. In (\ref{eq7}), (\ref{eq8}) and (\ref{eq9}), $\chi_\delta$ is a shear relaxation function for $\delta=2,3,$ or $4$. $H(t)$ is the Heaviside function.
 
 The relaxation matrix shown in (\ref{eq6}) generalizes the anisotropic relaxation model (\ref{eq4}) as three relaxation functions are used to describe  anelastic properties of shear modes.This also allows for the control of the dissipation of energy along three preferred directions e.g. the principal axes of the anisotropic medium.
 
 The choice of the relaxation functions depends on the symmetry system of the material, since attenuation symmetries follow the symmetry of the crystallographic form of material \cite{carcione1990}. In this study, the following relaxation functions are used \cite{carcione2014}
 \begin{align} \label{eq10}
 \chi_{\nu}(t)= L_\nu \left( \sum_{l=1}^{L_\nu} \myfrac{\tau_{\epsilon l }^{(\nu)}}{\tau_{\sigma l}^{(\nu)}} \right)\left[1 - \myfrac{1}{L_{\nu}}\sum_{l=1}^{L_\nu} \left(1 -\myfrac{\tau_{\epsilon l}^{(\nu)}}{\tau_{\sigma l}^{(\nu)}} \right)\exp\left({-t/ \tau_{\sigma l}^{(\nu)}}\right)\right], \qquad  \nu=1,...,4,
 \end{align}
where $\tau_{\epsilon l}^{(\nu)}$ and $\tau_{\sigma l}^{(\nu)}$ are material relaxation times such that  $\tau_{\epsilon l}^{(\nu)} \geq \tau_{\sigma l}^{(\nu)}$. The pair $\tau_{\epsilon l}^{(\nu)}$ and $\tau_{\sigma l}^{(\nu)}$ define a dissipation mechanism.

Equation (\ref{eq10}) describes the relaxation function of generalized standard linear solid (also known as Zener model) consisting of $L_\nu$ elements \footnote{A mechanical system in which a spring and a parallel combination of a dashpot and a spring are connected in series.} connected in parallel. The complex modulus of the system is \cite{carcione2014}
\begin{align} \label{eq11}
M_\nu(\omega)=\mathcal{F}\left(\myfrac{d [\chi_\nu (t) H(t)]}{dt} \right),
\end{align}
where $\omega$ is the angular frequency and $\mathcal{F}(\cdot)$ represents the time Fourier transform of the variable.

$M_{\nu}(\omega)$ is expressed as
\begin{align} \label{eq12}
M_\nu(\omega)=\left(\sum_{l=1}^{L_\nu} \myfrac{\tau_{\epsilon l}^{(\nu)}}{\tau_{\sigma l}^{(\nu)}} \right)\sum_{l=1}^{L_\nu} \myfrac{1 + i\omega \tau_{\epsilon l} ^{(\nu)}}{1 + i \omega \tau_{\sigma l}^{(\nu)}}.
\end{align}
 From (\ref{eq12}) it can be easily seen that $M_\nu(\omega)=1$ as $\myfrac{\tau_{\epsilon l }^{(\nu)}}{\tau_{\sigma l}^{(\nu)}} \rightarrow 1$, which gives the low frequency limit. Thus, (\ref{eq12}) is a general relaxation function which can recover all possible type of frequency behavior of attenuation and velocity dispersion observed in subsurface materials. 

\subsection{Strain memory variables}
The time-domain stress–strain relation can be expressed as
\begin{align}\label{R1}
\sigma_I={\psi}_{IJ}*\partial_t e_J
\end{align}
Applying the Boltzmann operation (\ref{eq3}) to (\ref{R1}), we obtain
\begin{align} \label{R2}
\sigma_I = \mathring{\psi}_{IJ} e_J + \left( \dot{\psi}_{IJ} H\right)*e_J
\end{align}
Now we use 
\[\widecheck{\phi}_{\nu l} (t)=\dot{\chi}_{\nu}(t)=\myfrac{1}{\tau_{\sigma l}^{(\nu)}} \left(\sum_{l=1}^{L_\nu} \myfrac{\tau_{\epsilon l }^{(\nu)}}{\tau_{\sigma l}^{(\nu)}}\right)\left( 1- \myfrac{\tau_{\epsilon l} ^{(\nu)}}{\tau_{\sigma l}^{(\nu)}}\right) \exp\left(-t/\tau_{\sigma l}^{(\nu)}\right)\]
and write (\ref{R2}) in matrix form, which is expressed as
\begin{align}
\sigma_I=A_{IJ}^{(\nu)} e_J + B_{IJ}^{(\nu)} \sum_{l=1}^{L_\nu} e_{Jl}^{(\nu)},
\end{align}
where $A'$s and $B'$s are the matrices formed by the combination of elastic constants $c_{IJ}$ and 
\begin{align*}
e_{Jl}^{(v)}=\phi_{\nu l} (t)* e_J, \qquad J=1,..,6,\qquad l=1,..,L_{\nu},~~\nu=1,..,4,
\end{align*} 
where $\phi_{\nu l }=\widecheck{\phi}_{\nu l} (t)H(t)$
are the components of $6 \times 1$ strain memory array $\bm{e}_l^{(\nu)}$. 

In 3D, the symmetric strain memory tensor corresponding to the $l^{th}$ dissipation mechanism of the relaxation function $\chi_\nu$ is expressed as \cite{carcione1995}
\begin{align} 
\bm{e_l^{(\nu)}}& =\left[\begin{array}{ccc}
e_{11l}^{(\nu)} & e_{12l}^{(\nu)} & e_{13l}^{(\nu)}\\
~ & e_{22l}^{(\nu)} & e_{23l}^{(\nu)}\\
~&~& e_{33l}^{(\nu)}
\end{array}\right]=\phi_{\nu l } *\left[ \begin{array}{ccc}
\epsilon_{11} & \gamma_{12} & \gamma_{13}\\
~ & \epsilon_{22} & \gamma_{23} \\
~&~& \epsilon_{33} 
\end{array}
 \right] \label{eq13} \\
 &=\phi_{\nu l}(t) * e_J. \label{eq14}
\end{align}
The tensor $\bm{e_l^{(\nu)}}$ contains the past history of material due the dissipation mechanism defined in (\ref{eq10}). In the pure elastic case $\tau_{\epsilon l }^{(\nu)} \rightarrow \tau_{\sigma l}^{(\nu)},~\phi_{\nu l } \rightarrow 0$ and $\bm{e_l^{(\nu)}}$ vanishes.

Similar to the strain tensor, the memory strain variable can be decomposed as 
\begin{align}  \label{eq15}
\bm{e_l}^{(\nu)} = \bm{d}_{l}^{(\nu)} + \left(\myfrac{1}{3} \text{tr}(\bm{e}_l^{(\nu)}) \right)\bm{I}, \qquad \text{tr}\left(\bm{d}_{l}^{(\nu)}\right) = 0,
\end{align}
where $\bm{d}_{l}^{(\nu)}$ is the deviatoric strain memory tensor which is traceless and  $\bm{I}$ is $3\times3$ identity matrix.\\
Thus, the dilatation and shear memory variables are defined as 
\begin{align} \label{eq16}
e_{1 l} = \text{tr}\left(\bm{e}_l^{(1)}\right), \qquad~\text{and}~\qquad e_{ijl} = \left(\bm{d}_{l}^{(\nu)}\right)_{ij},
\end{align}
where $\nu=\delta$ for $i=j$, $\nu=2$ for $ij=23$, $\nu=3$ for $ij=13$ and $\nu=4$ for $ij=12$.

The stress-strain relations in terms of strain components and memory variables with one dissipation mechanism for each mode are \cite{carcione1995}
\begin{subequations} \label{eq17}
	\begin{align}
	\sigma_{11}&=c_{11} \epsilon_{11} + c_{12} \epsilon_{22} +  c_{13} \epsilon_{33} +  {c}_{14} \gamma_{23} +  {c}_{15} \gamma_{13} + {c}_{16} \gamma_{13} + Ke_{11} + 2Ge_{111}^{(\delta)} \\
	\sigma_{22}&=c_{12} \epsilon_{11} + c_{22} \epsilon_{22} +  c_{23} \epsilon_{33} +  {c}_{24} \gamma_{23} +  {c}_{25} \gamma_{13} +  {c}_{26} \gamma_{12} + Ke_{11} + 2Ge_{221}^{(\delta)} \\
	\sigma_{33}&=c_{13} \epsilon_{11} + c_{23} \epsilon_{22} +  c_{33} \epsilon_{33} +  {c}_{34} \gamma_{23} +  {c}_{35} \gamma_{13} +  {c}_{36} \gamma_{12} + Ke_{11} - 2G\left(e_{111}^{(\delta)} + e_{221}^{(\delta)} \right) \\
	\sigma_{23}&={c}_{14} \epsilon_{11} + {c}_{24} \epsilon_{22} +  {c}_{34} \epsilon_{33} +  c_{44} \gamma_{23} +  {c}_{44} e_{231}^{(2)}+  {c}_{45} \gamma_{13} + c_{46} \gamma_{13} \\
	\sigma_{13}&={c}_{15} \epsilon_{11} + {c}_{25} \epsilon_{22} +  {c}_{35} \epsilon_{33} +  {c}_{45} \gamma_{23} +  c_{55} \gamma_{23}+  {c}_{55} e_{131}^{(3)} + c_{56} \gamma_{12} \\
	\sigma_{12}&={c}_{16} \epsilon_{11} + {c}_{26} \epsilon_{22} +  {c}_{36} \epsilon_{33} +  {c}_{46} \gamma_{23} +  {c}_{56} \gamma_{13}+  c_{66}  \gamma_{13} + c_{66}e_{121}^{(4)} \\
	\end{align}
\end{subequations}
where $c_{IJ}=\psi_{IJ}(t=0^+)$ are unrelaxed elasticity constant at $\omega \rightarrow \infty$.

\subsection{Memory variable equation}
Applying the Boltzman equation to the deviatoric part of (\ref{eq14}) yields
\begin{align} \label{mem1}
\partial_t \bm{d}_l^{(\nu)}=\phi_{\nu l } (0) \bm{d} + (\partial_t \widecheck{\phi}_{\nu l }H)*\bm{d}.
\end{align}
Here, $\bm{d}$ denotes the deviatoric strain tensor with elements
\begin{align}
\bm{d}=\bm{\epsilon}- \myfrac{1}{3}\mathcal{V}\bm{I},
\end{align}
where the strain tensor $\bm{\epsilon}$ and $\mathcal{V}$ are
\[\bm{\epsilon}=\left[ \begin{array}{ccc}
\epsilon_{11} & \epsilon_{12} &\epsilon_{13} \\
\epsilon_{12} & \epsilon_{22} &\epsilon_{23} \\
\epsilon_{13} & \epsilon_{23} &\epsilon_{33} \\
\end{array}\right], \qquad \mathcal{V}=\epsilon_{11} + \epsilon_{22} + \epsilon_{33}.
\] 
Using $\partial_t\widecheck{ \phi}_{\nu l}=-\myfrac{\widecheck{\phi}_{\nu l}}{\tau_{\sigma l}^{(\nu)}}$ and substituting it in (\ref{mem1}), we recover
\begin{align}\label{mem2}
\partial_t \bm{d}_l^{(\nu)}=\phi_{\nu l}(0)\bm{d} - \myfrac{1}{\tau_{\sigma l}^{(\nu)}} \bm{d}_l^{(\nu)},
\end{align}
where $\bm{d}_l^{(\nu)}=\phi_{\nu l} (t)*\bm{d}$, with $\nu=2,~3,~\text{and}~4$.
Similarly applying the Boltzmann operation to the non-deviatoric part $\text{tr}(\bm{e}_l^{(1)})$, we get
\begin{align} \label{m4}
\partial_t \text{tr}(\bm{e}_l^{(1)})=\phi_{1l}(0) \text{tr}{(\bm{e})} - \myfrac{1}{\tau_{\sigma l}^{(1)}} \text{tr}(\bm{e}_l^{(1)})
\end{align}
Using (\ref{mem2}) and (\ref{m4}), The equations for the memory variables are expressed as 
\begin{subequations} \label{eq18}
	\begin{align}
	\partial_t e_{111}^{(\delta)}&= \phi_{\delta 1} (0) (\epsilon_{11}  - \bar{\epsilon}) - \myfrac{e_{111}^{(\delta)}}{\tau_{\sigma}^{(\delta)}}\\
	\partial_t e_{221}^{(\delta)} &= \phi_{\delta 1}(0) (\epsilon_{22} - \bar{\epsilon}) - \myfrac{e_{221}^{(\delta)}}{\tau_{\sigma}^{(\delta)}}\\
	\partial_t e_{231}&= \phi_{21}(0)\gamma_{23} - \myfrac{e_{231}}{\tau_{\sigma}^{(2)}}\\
	\partial_t e_{131}&= \phi_{31}(0)\gamma_{13} - \myfrac{e_{131}}{\tau_{\sigma}^{(3)}}\\
	\partial_t e_{121}&= \phi_{41}(0)\gamma_{12} - \myfrac{e_{121}}{\tau_{\sigma}^{(4)}}\\
		\partial_t e_{11} &= n \phi_{1l} (0) \bar{\epsilon} - e_{11}/ \tau_{\sigma}^{(1)} 
	\end{align}
\end{subequations}
where $\bar{\epsilon}=tr(\bm{S})/3$ and $n$ is taken as 2 for 2D and 3 for 3D.

\subsection{Equation of motion}
The conservation of momentum is expressed as
\begin{subequations}  \label{eq19}
	\begin{align}
	\myfrac{\partial \sigma_{11}}{\partial x_1} + \myfrac{\partial \sigma_{12}}{\partial x_2} + \myfrac{\partial \sigma_{13}}{\partial x_3}&=\rho \myfrac{\partial v_1}{\partial t} \\
	\myfrac{\partial \sigma_{12}}{\partial x_1} + \myfrac{\partial \sigma_{22}}{\partial x_2} + \myfrac{\partial \sigma_{23}}{\partial x_3}&=\rho \myfrac{\partial v_2}{\partial t} \\
	\myfrac{\partial \sigma_{13}}{\partial x_1} + \myfrac{\partial \sigma_{23}}{\partial x_2} + \myfrac{\partial \sigma_{33}}{\partial x_3}&=\rho \myfrac{\partial v_3}{\partial t} 
	\end{align}	
\end{subequations}

\subsection{System of equations in matrix form in 3D}
Let us consider the 3D-particle-velocity and stress equations for propagation in an  anisotropic medium. We assign one relaxation mechanism to both dilatational anelastic deformation $(\nu=1)$ and shear anelastic deformations $(\nu=2)$. The stress-strain relation is expressed as

\begin{align}  \label{eq20}
\begin{aligned}   
\myfrac{\partial{\sigma_{11}}} {\partial t} &= c_{11}\myfrac {\partial v_1}{\partial x_1} + c_{12} \myfrac{\partial v_2}{\partial x_2} +  c_{13}\myfrac {\partial v_3}{\partial x_3} + K e_1 + 2 G e_2,\\
\myfrac{\partial{\sigma_{22}}} {\partial t} &= c_{12}\myfrac {\partial v_1}{\partial x_1} + c_{11} \myfrac{\partial v_2}{\partial x_2} +  c_{13}\myfrac {\partial v_3}{\partial x_3} + K e_1 + 2 Ge_3,\\
\myfrac{\partial{\sigma_{33}}} {\partial t} &= c_{13}\myfrac {\partial v_1}{\partial x_1} + c_{13}\myfrac {\partial v_2}{\partial x_2} + c_{33}\myfrac {\partial v_3}{\partial x_3} + K e_1 - 2 G (e_2 + e_3)\\
\myfrac{\partial{\sigma_{23}}} {\partial t} &= c_{44}\left[\left(\myfrac {\partial v_2}{\partial x_3} +\myfrac {\partial v_3}{\partial x_2}\right) + e_4\right],\\
\myfrac{\partial{\sigma_{13}}} {\partial t} &= c_{55}\left[\left(\myfrac {\partial v_1}{\partial x_3} +\myfrac {\partial v_3}{\partial x_1}\right) + e_5\right],\\
\myfrac{\partial{\sigma_{12}}} {\partial t} &= c_{66}\left[\left(\myfrac {\partial v_1}{\partial x_2} +\myfrac {\partial v_2}{\partial x_1}\right) + e_6\right].
\end{aligned}
\end{align} 
Memory variables are expressed as
\begin{align} \label{eq21}
\begin{aligned}
\myfrac{\partial e_1}{\partial t}&=\myfrac{1}{\tau_{\sigma}^{(1)}}\left[\left( \myfrac{\tau_\sigma^{(1)}}{{\tau_{\epsilon}^{(1)}}}-1\right) \left(\myfrac{\partial v_1}{\partial x_1} + \myfrac{\partial v_2}{\partial x_2} + \myfrac{\partial v_3}{\partial x_3} \right) -e_1 \right], \\
\myfrac{\partial e_2}{\partial t}&=\myfrac{1}{3\tau_{\sigma}^{(2)}}\left[\left( \myfrac{\tau_\sigma^{(2)}}{{\tau_{\epsilon}^{(2)}}}-1\right) \left(2\myfrac{\partial v_1}{\partial x_1} - \myfrac{\partial v_2}{\partial x_2}  - \myfrac{\partial v_3}{\partial x_3} \right) -3e_2\right], \\
\myfrac{\partial e_3}{\partial t}&=\myfrac{1}{3\tau_{\sigma}^{(3)}}\left[\left( \myfrac{\tau_\sigma^{(3)}}{{\tau_{\epsilon}^{(3)}}}-1\right) \left(2\myfrac{\partial v_2}{\partial x_2} - \myfrac{\partial v_1}{\partial x_1}  - \myfrac{\partial v_3}{\partial x_3} \right) -3e_3\right], \\
\myfrac{\partial e_4}{\partial t}&=\myfrac{1}{\tau_{\sigma}^{(2)}}\left[\left( \myfrac{\tau_\sigma^{(2)}}{{\tau_{\epsilon}^{(2)}}}-1\right) \left(\myfrac{\partial v_2}{\partial x_3} + \myfrac{\partial v_3}{\partial x_2} \right) -e_4\right], \\
\myfrac{\partial e_5}{\partial t}&=\myfrac{1}{\tau_{\sigma}^{(3)}}\left[\left( \myfrac{\tau_\sigma^{(3)}}{{\tau_{\epsilon}^{(3)}}}-1\right) \left(\myfrac{\partial v_1}{\partial x_3} + \myfrac{\partial v_3}{\partial x_1} \right) -e_5\right], \\
\myfrac{\partial e_6}{\partial t}&=\myfrac{1}{\tau_{\sigma}^{(4)}}\left[\left( \myfrac{\tau_\sigma^{(4)}}{{\tau_{\epsilon}^{(4)}}}-1\right) \left(\myfrac{\partial v_1}{\partial x_2} + \myfrac{\partial v_2}{\partial x_1} \right) -e_6\right].\\
\end{aligned}
\end{align}
Combining (\ref{eq19})-(\ref{eq21}) in matrix form yields
\begin{align} \label{eq23}
\myfrac{\partial \bm{q}}{\partial t} + \bm{A(\bm{x})} \myfrac{\partial \bm{q}} {\partial x_1} + \bm{B}(\bm{x}) \myfrac{\partial \bm{q}}{\partial x_2} + \bm{C}(\bm{x}) \myfrac{\partial \bm{q}}{\partial x_3}= \bm{D(\bm{x})q} +\bm{f},
\end{align}
where
\[
\bm{q}=\left[\begin{array}{cccccccc} 
\sigma_{11},~\sigma_{22},~\sigma_{33},~\sigma_{23},~\sigma_{13},~\sigma_{12},~e_1,~e_2,~e_3,~e_4,~e_5,~e_6,~v_1,~v_2,~v_3
\end{array}
\right],
\].
\[
\bm{A}(\bm{x})=-\left[\begin{array}{ccccccccccccccc}
0 & 0 & 0 & 0 & 0 & 0 & 0 & 0 & 0  & 0 & 0 & 0 & c_{11}(\bm{x}) & 0 & 0\\
0 & 0 & 0 & 0 & 0 & 0 & 0 & 0 & 0  & 0 & 0 & 0 & c_{12}(\bm{x}) & 0 & 0\\
0 & 0 & 0 & 0 & 0 & 0 & 0 & 0 & 0  & 0 & 0 & 0 & c_{13}(\bm{x}) & 0 & 0\\
0 & 0 & 0 & 0 & 0 & 0 & 0 & 0 & 0  & 0 & 0 & 0 & 0 & 0 &  0 \\
0 & 0 & 0 & 0 & 0 & 0 & 0 & 0 & 0  & 0 & 0 & 0 & 0 & 0 &  c_{55}(\bm{x}) \\
0 & 0 & 0 & 0 & 0 & 0 & 0 & 0 & 0  & 0 & 0 & 0 & 0 & c_{66}(\bm{x}) &  0 \\
0 & 0 & 0 & 0 & 0 & 0 & 0 & 0 & 0  & 0 & 0 & 0 & T_1(\bm{x}) & 0 & 0\\
0 & 0 & 0 & 0 & 0 & 0 & 0 & 0 & 0  & 0 & 0 & 0 & \myfrac{2}{3}T_2(\bm{x}) & 0 & 0\\
0 & 0 & 0 & 0 & 0 & 0 & 0 & 0 & 0  & 0 & 0 & 0 & -\myfrac{1}{3}T_3(\bm{x}) & 0 & 0\\
0 & 0 & 0 & 0 & 0 & 0 & 0 & 0 & 0  & 0 & 0 & 0 & 0 & 0 & 0\\
0 & 0 & 0 & 0 & 0 & 0 & 0 & 0 & 0  & 0 & 0 & 0 & 0 & 0 & T_3(\bm{x})\\
0 & 0 & 0 & 0 & 0 & 0 & 0 & 0 & 0  & 0 & 0 & 0 &  0 & T_4(\bm{x}) & 0 \\
1/\rho(\bm{x}) & 0 & 0 & 0 & 0 & 0 & 0 & 0 & 0  & 0 & 0 & 0 & 0  & 0 & 0 \\
0 & 0 & 0 & 0 & 0 & 1/\rho(\bm{x}) & 0 & 0 & 0  & 0 & 0 & 0 & 0 & 0 & 0 \\
0 & 0 & 0 & 0 & 1/{\rho(\bm{x})} & 0 & 0 & 0 & 0  & 0 & 0 & 0 & 0 & 0 & 0
\end{array}
 \right],
\]
\[
\bm{B}(\bm{x})=-\left[\begin{array}{ccccccccccccccc}
0 & 0 & 0 & 0 & 0 & 0 & 0 & 0 & 0  & 0 & 0 & 0 & 0 & c_{12}(\bm{x}) & 0 \\
0 & 0 & 0 & 0 & 0 & 0 & 0 & 0 & 0  & 0 & 0 & 0 & 0 & c_{11}(\bm{x}) & 0 \\
0 & 0 & 0 & 0 & 0 & 0 & 0 & 0 & 0  & 0 & 0 & 0 & 0 & c_{13}(\bm{x}) & 0 \\
0 & 0 & 0 & 0 & 0 & 0 & 0 & 0 & 0  & 0 & 0 & 0 & 0 & 0 &  c_{44}(\bm{x}) \\
0 & 0 & 0 & 0 & 0 & 0 & 0 & 0 & 0  & 0 & 0 & 0 & 0 & 0 &  0 \\
0 & 0 & 0 & 0 & 0 & 0 & 0 & 0 & 0  & 0 & 0 & 0 &  c_{66}(\bm{x}) & 0  & 0 \\
0 & 0 & 0 & 0 & 0 & 0 & 0 & 0 & 0  & 0 & 0 & 0 & 0 & T_1(\bm{x}) & 0 \\
0 & 0 & 0 & 0 & 0 & 0 & 0 & 0 & 0  & 0 & 0 & 0 & 0 & -\myfrac{1}{3}T_2(\bm{x}) & 0\\
0 & 0 & 0 & 0 & 0 & 0 & 0 & 0 & 0  & 0 & 0 & 0 &  0 & \myfrac{2}{3}T_3(\bm{x}) & 0\\
0 & 0 & 0 & 0 & 0 & 0 & 0 & 0 & 0  & 0 & 0 & 0 & 0 & 0 & T_2(\bm{x})\\
0 & 0 & 0 & 0 & 0 & 0 & 0 & 0 & 0  & 0 & 0 & 0 & 0 & 0 & 0\\
0 & 0 & 0 & 0 & 0 & 0 & 0 & 0 & 0  & 0 & 0 & 0 &  T_4(\bm{x}) & 0 & 0 \\
0 & 0 & 0 & 0 & 0 & 1/\rho(\bm{x}) & 0 & 0 & 0  & 0 & 0 & 0 &  0 & 0  & 0\\
0 & 1/\rho(\bm{x}) & 0 & 0 & 0 & 0 & 0 & 0 & 0  & 0 & 0 & 0 & 0 & 0 & 0 \\
0 & 0 & 0 & 1/\rho(\bm{x}) & 0 & 0 & 0 & 0 & 0  & 0 & 0 & 0 & 0 & 0 & 0 
\end{array}
\right],
\]
\[
\bm{C}(\bm{x})=-\left[\begin{array}{ccccccccccccccc}
0 & 0 & 0 & 0 & 0 & 0 & 0 & 0 & 0  & 0 & 0 & 0 & 0 & 0 &  c_{13}(\bm{x}) \\
0 & 0 & 0 & 0 & 0 & 0 & 0 & 0 & 0  & 0 & 0 & 0 & 0 & 0 &  c_{13}(\bm{x}) \\
0 & 0 & 0 & 0 & 0 & 0 & 0 & 0 & 0  & 0 & 0 & 0 & 0 & 0 & c_{33}(\bm{x}) \\
0 & 0 & 0 & 0 & 0 & 0 & 0 & 0 & 0  & 0 & 0 & 0 & 0 &  c_{44}(\bm{x})  & 0 \\
0 & 0 & 0 & 0 & 0 & 0 & 0 & 0 & 0  & 0 & 0 & 0 & c_{55}(\bm{x}) & 0 &  0 \\
0 & 0 & 0 & 0 & 0 & 0 & 0 & 0 & 0  & 0 & 0 & 0 &  0 & 0  & 0 \\
0 & 0 & 0 & 0 & 0 & 0 & 0 & 0 & 0  & 0 & 0 & 0 & 0 & 0 & T_1(\bm{x})  \\
0 & 0 & 0 & 0 & 0 & 0 & 0 & 0 & 0  & 0 & 0 & 0 & 0 &  0 &-\myfrac{1}{3}T_2(\bm{x}) \\
0 & 0 & 0 & 0 & 0 & 0 & 0 & 0 & 0  & 0 & 0 & 0 &  0 & 0 & \myfrac{2}{3}T_3(\bm{x}) \\
0 & 0 & 0 & 0 & 0 & 0 & 0 & 0 & 0  & 0 & 0 & 0 & 0 & T_2(\bm{x}) & 0\\
0 & 0 & 0 & 0 & 0 & 0 & 0 & 0 & 0  & 0 & 0 & 0 & T_3(\bm{x}) & 0 & 0\\
0 & 0 & 0 & 0 & 0 & 0 & 0 & 0 & 0  & 0 & 0 & 0 &  0 & 0 & 0 \\
0 & 0 & 0 & 0 & 1/\rho(\bm{x}) & 0 & 0 & 0 & 0  & 0 & 0 & 0 &  0 & 0  & 0\\
0 & 0 & 0 & 1/\rho(\bm{x})  & 0 & 0 & 0 & 0 & 0  & 0 & 0 & 0 & 0 & 0 & 0 \\
0 & 0 & 1/\rho(\bm{x})  & 0 & 0 & 0 & 0 & 0 & 0  & 0 & 0 & 0 & 0 & 0 & 0  
\end{array}
\right],
\]
\[
\bm{D}(\bm{x})=-\left[\begin{array}{ccccccccccccccc}
0 & 0 & 0 & 0 & 0 & 0 & K(\bm{x}) & 2G(\bm{x}) & 0  & 0 & 0 & 0 & 0 & 0 &  0 \\
0 & 0 & 0 & 0 & 0 & 0 & K(\bm{x}) & 0 & 2G(\bm{x})  & 0 & 0 & 0 & 0 & 0 & 0 \\
0 & 0 & 0 & 0 & 0 & 0 & K(\bm{x}) & -2G(\bm{x}) & -2G(\bm{x})  & 0 & 0 & 0 & 0 & 0 & 0 \\
0 & 0 & 0 & 0 & 0 & 0 & 0 & 0 & 0  & c_{44}(\bm{x}) & 0 & 0 & 0 &  0  & 0 \\
0 & 0 & 0 & 0 & 0 & 0 & 0 & 0 & 0  & 0 & c_{55}(\bm{x}) & 0 & 0 &  0  & 0 \\
0 & 0 & 0 & 0 & 0 & 0 & 0 & 0 & 0  & 0 & 0 & c_{66}(\bm{x}) & 0 &  0  & 0 \\
0&0&0&0&0&0& -\myfrac{1}{\tau_\sigma^{(1)}(\bm{x})}&0&0&0&0&0&0&0&0\\
0&0&0&0&0&0& 0 &  -\myfrac{1}{\tau_\sigma^{(2)}(\bm{x})}&0&0&0&0&0&0&0\\
0&0&0&0&0&0& 0 & 0& -\myfrac{1}{\tau_\sigma^{(3)}(\bm{x})}&0&0&0&0&0&0\\
0&0&0&0&0&0& 0 & 0&0 & -\myfrac{1}{\tau_\sigma^{(2)}(\bm{x})}&0&0&0&0&0\\
0&0&0&0&0&0& 0 & 0&0 & 0&  -\myfrac{1}{\tau_\sigma^{(3)}(\bm{x})}&0&0&0&0\\
0&0&0&0&0&0& 0 & 0&0 & 0 & 0 & -\myfrac{1}{\tau_\sigma^{(4)}(\bm{x})}&0&0&0\\
0&0&0&0&0&0& 0 & 0&0 & 0 & 0 & 0&0&0&0\\
0&0&0&0&0&0& 0 & 0&0 & 0 & 0 & 0&0&0&0\\
0&0&0&0&0&0& 0 & 0&0 & 0 & 0 & 0&0&0&0\\
\end{array}
\right],
\]
where $T_i=\myfrac{1}{\tau_{\sigma}^{(i)}(\bm{x})}\left( \myfrac{\tau_\sigma^{(i)}(\bm{x})}{{\tau_{\epsilon}^{(i)}}(\bm{x})}-1\right).$

To prove stability of the scheme, we express (\ref{eq23}) in a form where spatially dependent material coefficient appear on left side of (\ref{eq24}). This will enable us to rewrite (\ref{eq21}) without terms involving the spatial derivatives.  From (\ref{eq20}), we compute $\myfrac{\partial v_1}{\partial x_1}, \myfrac{\partial v_2}{\partial x_2}, \myfrac{\partial v_3}{\partial x_3}, \myfrac{\partial v_2}{\partial x_3}, \myfrac{\partial v_3}{\partial x_2}, \myfrac{\partial v_1}{\partial x_3}, \myfrac{\partial v_3}{\partial x_1}, \myfrac{\partial v_1}{\partial x_2},~\text{and}~\myfrac{\partial v_2}{\partial x_1}$ and substitute in (\ref{eq21}), which yields
\begin{align} \label{RE1}
\begin{aligned}
\myfrac{\partial a_1}{\partial t} &=-w_1(a_1 + z_1) - w_2 (a_2 +z_2) -2 w_3 (a_3 +z_3) \\
\myfrac{\partial a_2 }{\partial t}&=w_4 (a_1 +z_1) + w_5 (a_2 +z_2) + 2w_6(a_3 +z_3)\\
\myfrac{\partial a_3 }{\partial t}&=w_7 (a_1 +z_1)+ w_8  (a_2 +z_2) + 2w_9 (a_3 + z_3)\\
\myfrac{\partial a_4}{\partial t} &= T_2 (a_4 +z_4) -\myfrac{(a_4 +z_4) }{\tau_{\sigma}^{(2)}}\\
\myfrac{\partial a_5}{\partial t} &= T_3 (a_5 + z_5) -\myfrac{(a_5 +z_5) }{\tau_{\sigma}^{(3)}} \\
\myfrac{\partial a_6}{\partial t} &= T_4 (a_6 + z_6) -\myfrac{(a_6 +z_6) }{\tau_{\sigma}^{(4)}},
\end{aligned}
\end{align}
where 
\begin{align*}
a_1 &= e_1 -z_1, \qquad a_2=e_2 - z_2, \qquad a_3=e_3 - z_3,\\
a_4&=e_4 - z_4, \qquad a_5=e_5 - z_5, \qquad a_6=e_6 - z_6, 
\end{align*}
and 
\[
d_1=r_{11} + r_{12} + r_{13},\qquad d_2 =r_{33} +  2r_{13},
\]
with $r_{ij}$ being the elements of the inverse of unrelaxed compliance matrix\[\bm{C}_p(\bm{x})=\left[\begin{array}{cccccc}
c_{11}(\bm{x}) & c_{12}(\bm{x}) & c_{13}(\bm{x}) & 0 & 0 & 0\\
c_{12}(\bm{x}) & c_{22}(\bm{x}) & c_{13}(\bm{x}) & 0 & 0 & 0 \\
c_{13}(\bm{x}) & c_{13}(\bm{x}) & c_{33}(\bm{x}) & 0 & 0 & 0 \\
0& 0& 0 & c_{44}(\bm{x}) & 0 & 0 \\
0 &0 & 0 & 0 & c_{55}(\bm{x}) & 0 \\
0 &0 & 0 & 0 &0 & c_{66}(\bm{x}) \\
\end{array}\right], \]
and inverse $\bm{C}_p$ is given in  \ref{A1}.
\begin{align*}
z_1&=T_1 (d_1 (\sigma_{11} + \sigma_{22}) + d_2 \sigma_{33}),\\
z_2 &= T_2\left[\left(r_{11} - \myfrac{1}{3}d_1 \right)\sigma_{11} + \left(r_{12} - \myfrac{1}{3}d_1 \right)\sigma_{22} + \left(r_{13} - \myfrac{1}{3}d_2 \right)\sigma_{33}\right],\\
z_3&=T_3\left[\left(r_{12} - \myfrac{1}{3}d_1 \right)\sigma_{11} + \left(r_{11} - \myfrac{1}{3}d_1 \right)\sigma_{22} + \left(r_{13} - \myfrac{1}{3}d_2 \right)\sigma_{33}\right],\\
z_4&=T_2 c_{44}^{-1}\sigma_{23},\qquad z_5=T_3 c_{55}^{-1}\sigma_{13},\qquad z_6=T_4 c_{66}^{-1}\sigma_{12},
\end{align*}
and 
\begin{align*}
w_1&=\left(T_1 p_\lambda +  \myfrac{1} {\tau_\sigma^{(1)}}\right),\qquad & w_2&=T_1 p_{\mu_1},\qquad & w_3&=T_1 p_{\mu_2},\\
w_4&=\left(\myfrac{T_2 p_\lambda}{3} -  \myfrac{1} {\tau_\sigma^{(2)} } - K({r_{11} + r_{12})} \right),\qquad & w_5&=  (T_2 p_{\mu_1} -G(r_{11} -r_{13})), \qquad & w_6&=(T_2 p_{\mu_2} -G(r_{12}- r_{13})),\\
w_7&=\left(\myfrac{T_3 p_\lambda}{3} -  \myfrac{1} {\tau_\sigma^{(3)} } - K({r_{12} + r_{11})} \right),\qquad& w_8&=(T_3 p_{\mu_1} -G(r_{12} -r_{13})),\qquad& w_9&=(T_3 p_{\mu_2} -G(r_{11}- r_{13})),
\end{align*}
with \[p_\lambda=K(2d_1 + d_2),\qquad p_{\mu_1}= G(d_1 - 2 d_2), \qquad p_{\mu_2}=2G(d_1 - d_2).\]
 We rewrite the system of equations (\ref{eq21}) with out spatial derivative using the set of equations in (\ref{eq20}), which yields
\begin{align} \label{eq24}
\begin{aligned}
\bm{Q}_{s}^{-1} (\bm{x}) \myfrac{\partial \bm{\bm{\sigma}}}{\partial t}&=\sum_{i=1}^{d} \bm{A}_i \myfrac{\partial \bm{v}}{\partial \bm{x}_i} + \bm{S\sigma}, \\
{\rho} \myfrac{\partial \bm{\bm{v}}}{\partial t}&=\sum_{i=1}^{d} \bm{A}_i^{T} \myfrac{\partial \bm{\sigma}}{\partial \bm{x}_i} +\bm{f}, \\
\end{aligned}	
\end{align}
where, $\bm{\sigma}=[\sigma_{11},~\sigma_{22},~\sigma_{33},~\sigma_{23}, \sigma_{13}, \sigma_{12},~a_1,~a_2,~a_3,~a_4,~a_5,~a_6]^T,~\text{and}~\bm{v}=[v_1,~v_2,~v_3]^T$ and 
\[\bm{Q}_{s}^{-1}(\bm{x})=\left[\begin{array}{c|c} 
\bm{C}^{-1}_{p}(\bm{x}) & \bm{0} \\
~\\
\hline~~\\
\bm{0} & \bm{I}\\\end{array}\right].\]
The matrices $\bm{A}_i$ and $\bm{S}$ are
\[
\bm{A}_1=\left[ \begin{array}{ccc} 
1&0&0\\
0&0&0\\
0&0&0\\
0&0&0\\
0&0&1\\
0&1&0\\
\hline
0&0&0\\
0&0&0\\
0&0&0\\
0&0&0\\
0&0&0\\
0&0&0\\
\end{array}
\right],\qquad \bm{A}_2=\left[ \begin{array}{ccc} 
0&0&0\\
0&1&0\\
0&0&0\\
0&0&1\\
0&0&0\\
1&0&0\\
\hline
0&0&0\\
0&0&0\\
0&0&0\\
0&0&0\\
0&0&0\\
0&0&0\\
\end{array}
\right],\qquad \bm{A}_3=\left[ \begin{array}{ccc} 
0&0&0\\
0&0&0\\
0&0&1\\
0&1&0\\
1&0&0\\
0&0&0\\
\hline
0&0&0\\
0&0&0\\
0&0&0\\
0&0&0\\
0&0&0\\
0&0&0\\
\end{array}
\right],
\]
and $\bm{S}=\bm{Q}_s^{-1} \bm{G}$ with\\
\[
\bm{G} = \left[\begin{array}{c|c} 
\bm{g}_{11_{6 \times 6}} & \bm{g}_{12_{6 \times 6}} \\
~\\
\hline
~\\
\bm{g}_{21_{6 \times 6}} & \bm{g}_{22_{6 \times 6}} \\
\end{array}\right]
\]
where
\begin{align*}
\bm{g}_{11}&=\left[\begin{array}{cccccccccccc} 
T_1d_1 + 2GT_2 g_{1}  &  K T_1d_1 + 2GT_2 g_{2} & K T_1 d_2 + 2GT_2 g_{3} & 0 & 0 & 0 \\
K T_1d_1 + 2GT_3 g_{4}  &  K T_1d_1 + 2GT_3 g_{5} & K T_1 d_2 + 2GT_3 g_{6} & 0 & 0 & 0 \\
K T_1d_1 - 2G g_{7} &  K T_1d_1 - 2G g_{8} & K T_1 d_2 - 2G g_{9} & 0 & 0 & 0 \\
0 & 0 & 0 & T_2 c_{44}^{-1} & 0 & 0 \\
0 & 0 & 0 & 0 & T_3 c_{55}^{-1} & 0 \\
0 & 0 & 0 & 0 & 0 & T_4 c_{66}^{-1}
\end{array}\right],\\
\bm{g}_{12}&=\left[\begin{array}{cccccc} 
K  & 2G  & 0 & 0 & 0 & 0\\
K  & 0 & 2G  & 0 & 0 & 0\\
K  & -2G  & -2G & 0 & 0 & 0\\
0 & 0& 0 & c_{44} & 0 & 0\\
0 & 0& 0 & 0 & c_{55}  & 0\\
0 & 0& 0 & 0 & 0  & c_{66}
\end{array}\right],\\
~ \bm{g}_{21}&=\left[ 
\begin{array}{cccccc}
t_1 & t_2 & t_3 & 0 & 0 & 0\\
t_4 & t_5 & t_6 & 0 & 0 & 0\\
t_7  & t_8 &t_9 & 0 & 0 & 0\\
0 & 0 & 0 & \left( T_2 - \myfrac{1}{\tau_{\sigma}^{(2)}}\right) T_2/c_{44} & 0 & 0 \\
0 & 0 & 0 & 0 & \left( T_3 - \myfrac{1}{\tau_{\sigma}^{(3)}}\right) T_3/c_{55} & 0 \\
0 & 0 & 0 & 0 & 0 & \left( T_4 - \myfrac{1}{\tau_{\sigma}^{(4)}}\right) T_4/c_{66} \\
\end{array}\right],
\end{align*}
\begin{align*}
\bm{g}_{22}&=\left[\begin{array}{cccccc}
-w_1& -w_2 & w_3 & 0 & 0 & 0\\
w_4 & w_5 & w_6 & 0 & 0 & 0\\
w_7 & w_8 & w_9 & 0 & 0 & 0\\
0  & 0 & 0 & \left( T_2 - \myfrac{1}{\tau_{\sigma}^{(2)}}\right) & 0 & 0\\
0  & 0 & 0 & 0 &  \left( T_3 - \myfrac{1}{\tau_{\sigma}^{(3)}}\right) & 0\\
0  & 0 & 0 & 0 & 0 &  \left( T_4 - \myfrac{1}{\tau_{\sigma}^{(4)}}\right)
\end{array}\right]
\end{align*}
where \begin{align*}
g_{1}&=\left( r_{11} - \myfrac{1}{3}d_1\right),\qquad & g_{2}&=\left(r_{12} - \myfrac{1}{3}d_1 \right),\qquad &g_{3}&=\left( r_{13} - \myfrac{1}{3}d_2\right),\\ g_4&=\left(r_{12} - \myfrac{1}{3}d_1 \right),\qquad & g_{5}&=\left(r_{11} - \myfrac{1}{3}d_1 \right), \qquad &g_{6}&=\left( r_{13} - \myfrac{1}{3}d_2\right),\\
 g_{7}&=T_2 g_{1} + T_3 g_{4},\qquad & g_{8}&=T_2 g_{2} + T_3 g_{5}, \qquad &g_{9}&=T_2 g_{3} + T_3 g_{6},\\
t_1&=-(w_1 T_1 d_1 + w_2 T_2 g_1 + 2w_3 T_3 g_4),\qquad &t_2&=-( w_1 T_1d_1  + w_2 T_2 g_2 + 2 w_3 g_5 ),\\
t_3&=-(w_1 T_1 d_2 + w_2 T_2 g_3 + 2 w_3 g_6),\qquad & t_4&=w_4 T_1 d_1 + w_5 T_2 g_1 + 2w_6 T_3 g_4,\\
t_5&=w_4 T_1d_1  + w_5 T_2 g_2 + 2 w_6 g_5, \qquad& t_6&=w_4 T_1 d_2 + w_5 T_2 g_3 + 2 w_6 g_6,\\
t_7&=w_7 T_1 d_1 + w_8 T_2 g_1 + 2w_9 T_3 g_4,\qquad& t_8&= w_7 T_1d_1  + w_8 T_2 g_2 + 2 w_9 g_5,\\
t_9&= w_7 T_1 d_2 + w_8 T_2 g_3 + 2 w_9 g_6 .
\end{align*}
Here elements of $\bm{G}$ and $\bm{g}_{ij}$ are space dependent.\\

It should be noted that the matrices $\bm{A}_i$ are spatially constant, while $\bm{Q}_{s}^{-1}$ and $\rho$ can vary spatially. We will also assume that $\rho$ , $\bm{Q}_{s}^{-1}$ and $\bm{Q}_{s}^{-1}$ are positive-definite and bounded pointwise such that 
\[
0 < \rho_{\min} \leq \rho \bm{(x)} \leq \rho_{\max}< \infty \\
\]
\[
0 < c_{\min} \leq \bm{u}^{T}\bm{Q}_s \bm{(x)} \bm{u} \leq c_{\max} < \infty \\
\]
\[
0 < \hat{c}_{\min} \leq \bm{u}^{T}\bm{Q}_s^{-1)} \bm{(x)} \bm{u} \leq \hat{c}_{\max} < \infty \\
\]
for all $\bm{x} \in \mathbb{R}^d$ and  $\forall~\bm{u} \in \mathbb{R}^{N_d}$.\\ Moreover, we assume that $\bm{S}$ is a semi negative-definite and bounded pointwise such that
\[
-\infty < s_{min} \leq \bm{s}^T \bm{S} \bm{s} \leq s_{max} < 0 \qquad \forall~\bm{x} \in \mathbb{R}^d~\text{and} ~\forall~\bm{u} \in \mathbb{R}^{N_d}. 
\]

\section{An energy stable discontinuous Galerkin formulation for the viscoelastic  wave equation}
\label{sec:ESdG}
Energy stable discontinuous Galerkin methods for viscoelastic wave propagation have been constructed based on  the formulations of the system shown in (\ref{eq24}). We assume that the domain $\Omega$ is exactly triangulated by a mesh $\Omega_h$ which consists of elements $D^k$ which are images of a reference element $\hat{D}$ under a local affine mapping.
\[
\bm{x}^{k}=\Phi^k \widehat{\bm{x}}.
\] 
Here $\bm{x}^k=\{x^k, y^k\}$ for $d=2$ and $\bm{x}^k=\{x^k, y^k, z^k\}$ for $d=3$ denote the physical coordinates on $D^k$ and $\hat{\bm{x}}=\{\hat{x}, \hat{y}\}$ for $d=2$ and $\widehat{\bm{x}}=\{\widehat{x}, \widehat{y}, \widehat{z}\}$ for $d=3$ denote coordinates on the reference element. We denote the determinant of the Jacobian of $\Phi^k$ as $J$. 

Solutions over each element $D^k$ are approximated from a local space $V_h(D^k)$, which is defined as composition of the mapping $\Phi^k$ and the reference approximation space $V_h(\widehat{D})$
\[
{V_h(D^k)=V_h(\widehat{D}) \circ \left(\Phi^k \right)^{-1} .}
\] 
Subsequently, the global approximation space $V_h(\Omega_h)$ is defined as 
\begin{align*}
V_h(\Omega_h)=\bigoplus_{k} V_h(D^k). 
\end{align*}
In this work, we will take $V_h(\widehat{D})=P^N(\widehat{D})$, with $P^N(\widehat{D})$ being the space of polynomials of total degree $N$ on the reference simplex. In two dimensions, $P^N$ on a triangle is 
\[
P^N(\widehat{D})=\{\widehat{x}^i\widehat{y}^j, 0 \le i +j \le N\},
\] 
and in three dimensions, $P^N$ on a tetrahedron is
\[
P^N(\widehat{D})=\{\widehat{x}^i\widehat{y}^j\widehat{x}^k, 0 \le i +j+k \le N\}.
\] 
The $L^2$ inner product and norm over $D^k$ is represented as
{
	\begin{align*}
	\LRp{\bm{g}, \bm{h}}=\int_{D^k} \bm{g} \cdot \bm{h}~d{\bm{x}} = \int_{\hat{D}} \left( \bm{g} \circ \Phi^k \right) \cdot \left(\bm{h} \circ \Phi^k\right) J~d\hat{\bm{x}}, \qquad ||\bm{g}||^2_{L^2{\left(\Omega\right)}} = (\bm{g}, \bm{g})_{L^2(D^k)},
	\end{align*}
}
where $\bm{g}$ and $\bm{h}$ are real vector-valued functions.  Global $L^2$ inner products and squared norms are defined as the sum of local $L^2$ inner products and squared norms over each elements. The $L^2$ inner product and norm over the boundary $\partial D^k$ of an element are similarly defined as
\[
\left< \bm{u}, \bm{v}  \right>_{L^2(\partial D^k)}=\int_{\partial D^k} \bm{u} \cdot \bm{v}~d\bm{x} = \sum_{f \in \partial D^k} \int_{\hat{f}} \bm{u} \cdot \bm{v} J^f~d\hat{\bm{x}}, \qquad ||\bm{u}||^2_{L^2(\partial D^k)}=\left<\bm{u}, \bm{u}\right>, 
\]
where $J^f$ is the Jacobian of the mapping from a reference face $\hat{f}$ to a physical face $f$ of an element.

Let $f$ be a face of an element $D^k$ with neighboring element $D^{k,+}$ and unit outward normal $\bm{n}$. Let $u$ be a function with discontinuities across element interfaces. We define the interior value $u^-$ and exterior value $u^+$ on face $f$ of $D^k$
\[
u^- = u|_{f \cap \partial D^k}, \qquad u^+ = u|_{f\cap \partial D^{k,+}}.
\]
The jump and average of a scalar function $u \in V_h(\Omega_h)$ over $f$ are then defined as
\[
\jump{u}=u^+ - u^-, \qquad \avg{u}=\myfrac{u^+ +  u^-}{2}.
\]
Jumps and averages of vector-valued functions $\bm{u} \in \mathbb{R}^{m}$ and and matrix-valued functions $\tilde{\bm{S}} \in \mathbb{R}^{m \times n}$  are defined component-wise.
\[
\left( \jump{\bm{u}}\right)_i = \jump{\bm{u}_i}, \qquad~1 \le i \le m\qquad~\left(\jump{\tilde{\bm{S}}}\right)_{ij}=\jump{\tilde{\bm{S}}}
\]
We can now specify a DG formulation for the viscoelastic wave equation (\ref{eq24}) which readily admits a DG formulation based on a penalty flux \cite{chan2018}. For the first order viscoelastic wave equation in (\ref{eq24}), the DG formulation in strong form  expressed as
\begin{equation}
\begin{aligned}
\label{eq:scehem}
\sum_{D^k \in \Omega_h} \left( \bm{Q}_{s}^{-1} \myfrac{\partial \bm{\sigma}}{\partial t} , \bm{h} \right)_{L^2(D^k)}=&\sum_{{D^k \in \Omega_h}} \Biggl (\left(\sum_{i=1}^{d}\bm{A}_i \myfrac{\partial \bm{v}}{\partial \bm{x}_i} , \bm{h}   \right)_{L^2(D^k)} + \left \langle \myfrac{1}{2} \bm{A_n}\jump{\bm{v}} + \myfrac{\alpha_{\bm{\sigma}}}{2} \bm{A_n}\bm{A}_n^T \jump{\bm{\sigma}}, \bm{h}  \right \rangle_{L^2(\partial D^k)}  \\ 
& + \left(\bm{S}\bm{\sigma}, \bm{g}\right)_{L^2 \left( {D^k} \right)} \Biggr) \\
\sum_{D^k \in \Omega_h} \left( \rho \myfrac{\partial \bm{v}}{\partial t} , \bm{g} \right)_{L^2(D^k)}=&\sum_{{D^k \in \Omega_h}} \Biggl(\left(\sum_{i=1}^{d}{\bm{A}_i}^T \myfrac{\partial \bm{\sigma}}{\partial \bm{x}_i} + \bm{f} ,\bm{g}   \right)_{L^2(D^k)}  + \left \langle \myfrac{1}{2} \bm{A_n}^T \jump{\bm{\sigma}} + \myfrac{\alpha_{\bm{v}}}{2} \bm{A}_n^T \bm{A_n} \jump{\bm{v}}, \bm{g}  \right \rangle_{L^2(\partial D^k)} \Biggr),
\end{aligned}
\end{equation}
for all $\bm{h},~\bm{g} \in V_h(\Omega_h)$.  Here, $\bm{A}_n$ is the normal matrix defined on a face $f$ of an element 
\[
\bm{A}_n=\sum_{i=1}^d n_i \bm{A}_i=\left[\begin{array}{ccc}
n_x &  0 & 0  \\
0 & n_y &  0 \\
0 & 0 & n_z\\
0 & n_z & n_y \\
n_z & 0 & n_x\\
n_y & n_x & 0\\
\hline
0 & 0 & 0\\
0 & 0 & 0\\
0 & 0 & 0\\
0 & 0 & 0\\
0 & 0 & 0\\
0 & 0 & 0\\
\end{array}\right].
\]
The factors $\alpha_{\tau},~\alpha_{v}$ are penalty parameters and defined on element interfaces. We assume that $\alpha_{\tau}~\text{and}~\alpha_{v} \ge 0$ and are piecewise constant over each shared face between two elements. These penalty constants can be taken to be zero, which results in a non-dissipative central flux, while $\alpha_{\tau},~\alpha_{v} > 0$  results in energy dissipation similar to the upwind flux \cite{hesthaven2007, chan2017weight, Shukla2020}. The stability of DG formulations are independent of the magnitude of these penalty parameters. However, a naive choice of these parameters can result in a stiffer semi-discrete system of ODEs and necessitates a smaller time under explicit time integration schemes \cite{chan2017weight, Shukla2020}.  In this work, we take $\alpha= 1/2$ unless stated otherwise. 

In most of DG formulations, material parameters are present in the numerical flux and in the penalty parameters $\alpha$.  However, for the presented formulation, the scheme is stable and high order accurate even when the penalty parameters are zero. The difference in the presented formulation is that the material data has been factored out into the mass matrix multiplying the time derivative.  Multiplying by the inverse mass matrix incorporates material parameters through an appropriate combination and scaling of the flux terms \cite{chan2017weight}.

\subsection{Boundary Conditions}
In many applications, the  top surface of a domain is a free surface boundary (stress free), with remaining surfaces taken to be absorbing boundaries. In our DG formulation, the boundary conditions are imposed by choosing appropriate exterior values which result in modified boundary numerical fluxes.  Boundary conditions on the normal components of the stress can be imposed by modifying the term $\jump{\bm{A}_n^T \bm{\tau}}=\jump{\widetilde{\bm{S}} \bm{n}}$ in numerical flux, where $\bm{\widetilde{S}}$ is
\[
\bm{\widetilde{S}}=\left[\begin{array}{ccccccc}
\sigma_{11} & \sigma_{12} & \sigma_{13} \\
\sigma_{12} & \sigma_{22} & \sigma_{23} \\
\sigma_{13} & \sigma_{23} & \sigma_{33} \\
\end{array}
\right].
\]
For a face of the element lying on the top surface of the domain, free surface boundary or zero traction conditions can be imposed by setting
\[
\jump{\bm{A}_n^T \bm{\sigma}}=\jump{\widetilde{\bm{S}} \bm{n}}=-2\widetilde{\bm{S}}^- \bm{n}=-2\bm{A}_n^T \bm{\sigma}^-,\qquad \bm{v}^+=\bm{v}^-\implies \jump{\bm{v}}=0.
\]
For problems which require the truncation of infinite or large domains, basic ``extrapolation'' absorbing boundary conditions can be imposed by setting 
\[
\jump{\bm{A}_n^T \bm{\sigma}}=\jump{\widetilde{\bm{S}} \bm{n}}=-\widetilde{\bm{S}}^- \bm{n}=-\bm{A}_n^T \bm{\sigma}^-,\qquad \bm{v}^+=0 \implies\jump{\bm{v}}=-\bm{v}^-.
\]
In addition to the above boundary conditions, more accurate absorbing boundary conditions can be also imposed using perfectly matching layers (PML) \cite{berenger1994} or high order absorbing boundary conditions (HABC) \cite{hagstrom2004}. However, in all cases, the boundary conditions are imposed by computing the numerical fluxes based on the modified jumps, ensuring energy stability for free surface and absorbing boundary conditions.  

\subsection{Energy stability}
The DG formulation in (\ref{eq:scehem}) can be proven to be energy stable in the absence of external forces $(\bm{f}=0)$ and for free-surface and absorbing boundary conditions. Integrating by parts the velocity equation in (\ref{eq:scehem}) gives
\begin{equation} \label{eq:es1}
\begin{aligned}
\sum_{D^k \in \Omega_h} \left( \bm{Q}_{s}^{-1} \myfrac{\partial \bm{\sigma}}{\partial t} , \bm{h} \right)_{L^2(D^k)}& =\sum_{{D^k \in \Omega_h}} \Biggl(\left(\sum_{i=1}^{d}\bm{A}_i \myfrac{\partial \bm{v}}{\partial \bm{x}_i} , \bm{h}   \right)_{L^2(D^k)} + \left \langle \myfrac{1}{2} \bm{A_n}\jump{\bm{v}} + \myfrac{\alpha_{\bm{\sigma}}}{2} \bm{A_n}\bm{A}_n^T \jump{\bm{\sigma}}, \bm{h}  \right \rangle_{L^2(\partial D^k)} \\
& + \left(\bm{S}\bm{\sigma}, \bm{\sigma}\right)_{L^2(D^k)}\Biggr) \\
\sum_{D^k \in \Omega_h} \left( \rho \myfrac{\partial \bm{v}}{\partial t} , \bm{g} \right)_{L^2(D^k)}=&\sum_{{D^k \in \Omega_h}} \Biggl(-\left(\sum_{i=1}^{d} \bm{\sigma},{\bm{A}_i} \myfrac{\partial \bm{g}}{\partial \bm{x}_i}\right)_{L^2(D^k)} \\ & + \left \langle  \bm{A_n}^T \avg{\bm{\sigma}} + \myfrac{\alpha_{\bm{v}}}{2} \bm{A}_n^T \bm{A_n} \jump{\bm{v}}, \bm{g}  \right. \rangle_{L^2(\partial D^k)}  
\Biggr)
\end{aligned}
\end{equation}
Taking $\left(\bm{h}, \bm{g}\right)=(\bm{\sigma}, \bm{v})$  in (\ref{eq:es1}) and adding both equations together yields
\begin{align*}
\sum_{D^k \in \Omega_h} \myfrac{1}{2}&\myfrac{\partial}{\partial t}\left( (\bm{Q}_s^{-1} \bm{\sigma}, \bm{\sigma})_{L^2(D^k)} + (\rho \bm{v}, \bm{v})_{L^2(D^k)}\right) \\
=&\sum_{D^k \in \Omega_h} \left \langle \myfrac{1}{2} \bm{A_n}\jump{\bm{v}} + \myfrac{\alpha_{\bm{\tau}}}{2} \bm{A_n}\bm{A}_n^T \jump{\bm{\sigma}}, \bm{h}  \right \rangle_{L^2(\partial D^k)}+ \left \langle  \bm{A_n}^T \avg{\bm{\sigma}} + \myfrac{\alpha_{\bm{v}}}{2} \bm{A}_n^T \bm{A_n} \jump{\bm{v}}, \bm{g}  \right \rangle_{L^2(\partial D^k)}\\ & + \left(\bm{S}\bm{\sigma}, \bm{\sigma}\right)_{L^2(D^k)}\\
=&\sum_{D^k \in \Omega_h}  \sum_{f \in \partial D^k}  \int_f\left(\myfrac{1}{2} \bm{\sigma}^T \bm{A}_n \jump{\bm{v}} + \myfrac{\alpha_{\bm{\tau}}}{2} \bm{\sigma}^T \bm{A_n} \bm{A}_n^T\jump{\bm{\sigma}}+\bm{v}^T \bm{A}_n^T\avg{\bm{\sigma}} + \myfrac{\alpha_{\bm{v}}}{2} \bm{v}^T \bm{A}_n^T \bm{A}_n \jump{\bm{v}}\right) \text{d}\bm{x} \\ \nonumber
& + \sum_{D^k \in \Omega_h} \int_{D^k} \bm{\sigma}^T \bm{S} \bm{\sigma}~\text{d}\bm{x},
\end{align*}
where the term 
\[\sum_{D^k \in \Omega_h} \myfrac{1}{2} \myfrac{\partial}{\partial t}\left( (\bm{Q}_s^{-1} \bm{\sigma}, \bm{\sigma})_{L^2(D^k)} + (\rho \bm{v}, \bm{v})_{L^2(D^k)}\right)\]
is the total energy of the system. 

Let $\Gamma_h$ be the set of unique faces in $\Omega_h$ and let $\Gamma_{\bm{\sigma}}$, $\Gamma_{\text{abc}}$ denote boundaries where free-surface and absorbing boundary conditions are imposed, respectively. We split surface terms into contributions from interior shared faces and from boundary faces. On interior shared faces, we sum the contributions from the two adjacent elements, which gives
\begin{align*}
&\sum_{f \in \Gamma_h \setminus \partial \Omega} \int_f\left(\myfrac{1}{2} \bm{\tau}^T \bm{A}_n \jump{\bm{v}} + \myfrac{\alpha_{\bm{\sigma}}}{2} \bm{\sigma}^T \bm{A_n} \bm{A}_n^T\jump{\bm{\sigma}}+\bm{v}^T \bm{A}_n^T\avg{\bm{\sigma}} + \myfrac{\alpha_{\bm{v}}}{2} \bm{v}^T \bm{A}_n^T \bm{A}_n \jump{\bm{v}}\right) \text{d}\bm{x} + \sum_{D^k \in \Omega_h} \int_{D^k} \bm{\sigma}^T \bm{D} \bm{\sigma}~\text{d}\bm{x}\\
&=-\sum_{f \in \Gamma_h \setminus \partial \Omega} \int_f \left( \myfrac{\alpha_{\bm{\sigma}}}{2} \left |\bm{A}_n^T \jump{\bm{\sigma}}\right|^2 + \myfrac{\alpha_{\bm{v}}}{2} \left|\bm{A}_n \jump{\bm{v}} \right|^2 \right)\text{d}\bm{x} + \sum_{D^k \in \Omega_h} \int_{D^k} \bm{v}^T \bm{S} \bm{v}~\text{d}\bm{x},
\end{align*}
where $\bm{v}^T \bm{S} \bm{v} < 0 $, since $\bm{S}$ is a negative semi-definite matrix.
For faces which lie on $\Gamma_{\bm{\sigma}}$ , $\bm{A}_n^T=-2\bm{A}_{n}^T \bm{\sigma}^-,~\bm{A}_n^T \avg{\bm{\sigma}}=0$ and $\jump{\bm{v}}=0$ yielding 
\begin{align*}
&\sum_{f \in \Gamma_{\bm{\sigma}}} \int_f\left(\myfrac{1}{2} \bm{\sigma}^T \bm{A}_n \jump{\bm{v}} + \myfrac{\alpha_{\bm{\sigma}}}{2} \bm{\sigma}^T \bm{A_n} \bm{A}_n^T\jump{\bm{\tau}}+\bm{v}^T \bm{A}_n^T\avg{\bm{\sigma}} + \myfrac{\alpha_{\bm{v}}}{2} \bm{v}^T \bm{A}_n^T \bm{A}_n \jump{\bm{v}}\right) \text{d}\bm{x} \\
&=-\sum_{f \in \Gamma_{\bm{\tau}}} \int_f \left ( \alpha_{\bm{\sigma}} \left| \bm{A}_n^T \bm{\sigma}^- \right|^2\right)~\text{d} \bm{x}.
\end{align*}

Finally, for faces in $\Gamma_{\text{abc}}$, we have $\bm{A}_n^T \avg{\bm{\sigma}}=\myfrac{1}{2} \bm{A}_n^T \bm{\sigma}{^-}$, 
$\bm{A}_n^T\jump{\bm{\sigma}}=-\bm{A}_n^T \bm{\sigma}^- $ and $\jump{\bm{v}}=-\bm{v}^-$, yielding
\begin{align*}
&\sum_{f \in \Gamma_h \setminus \partial \Omega} \int_f\left(\myfrac{1}{2} \bm{\sigma}^T \bm{A}_n \jump{\bm{v}} + \myfrac{\alpha_{\bm{\sigma}}}{2} \bm{\sigma}^T \bm{A_n} \bm{A}_n^T\jump{\bm{\sigma}}+\bm{v}^T \bm{A}_n^T\avg{\bm{\sigma}} + \myfrac{\alpha_{\bm{v}}}{2} \bm{v}^T \bm{A}_n^T \bm{A}_n \jump{\bm{v}}\right) \text{d}\bm{x}\\
&=-\sum_{f \in \Gamma_{\text{abc}}} \int_f \left( \myfrac{\alpha_{\bm{\sigma}}}{2} \left |\bm{A}_n^T \bm{\sigma}^-\right|^2 + \myfrac{\alpha_{\bm{v}}}{2} \left|\bm{A}_n \bm{v}^- \right|^2 \right)\text{d}\bm{x} ,
\end{align*}
Combining contributions from all faces and dissipation in the system yields the following result:
\begin{theorem}
	The DG formulation in (\ref{eq:scehem}) is energy stable for $\alpha_{\bm{\sigma}}, \alpha_{\bm{v}} \geq 0$ such that
	\begin{align}
	\label{eq:th1}
	\sum_{D^k \in \Omega_h} & \myfrac{1}{2} \myfrac{\partial}{\partial t}\left( (\bm{Q}_s^{-1} \bm{\sigma}, \bm{\sigma})_{L^2(D^k)} + (\bm{\rho} \bm{v}, \bm{v})_{L^2(D^k)}\right)=-\sum_{f \in \Gamma_h \setminus \partial \Omega} \int_f \left( \myfrac{\alpha_{\bm{\sigma}}}{2} \left |\bm{A}_n^T \jump{\bm{\sigma}}\right|^2 + \myfrac{\alpha_{\bm{v}}}{2} \left|\bm{A}_n \jump{\bm{v}} \right|^2 \right)\text{d}\bm{x} \nonumber \\
	&-\sum_{f \in \Gamma_{\bm{\sigma}}} \int_f \left ( \alpha_{\bm{\sigma}} \left| \bm{A}_n^T \bm{\sigma}^- \right|^2\right)~\text{d} \bm{x}-\sum_{f \in \Gamma_{\text{abc}}} \int_f \left( \myfrac{\alpha_{\bm{\sigma}}}{2} \left |\bm{A}_n^T \bm{\sigma}^-\right|^2 + \myfrac{\alpha_{\bm{v}}}{2} \left|\bm{A}_n \bm{v}^- \right|^2 \right)\text{d}\bm{x} \nonumber \\ & + \sum_{D^k \in \Omega_h} \int_{D^k} \bm{v}^T \bm{S} \bm{v}~\text{d}\bm{x} \le 0. 
	\end{align}
\end{theorem}

 The left hand side of (\ref{eq:th1}) is an $L^2$-equivalent norm on $(\bm{\tau}, \bm{v})$ as  $\bm{Q}_s^{-1}$ and $\rho$ are positive definite. Theorem 1 implies that magnitude of the DG solution is non-increasing in time dissipation is present for penalization parameters ${\alpha_{\bm{\tau}}}, {\alpha_{\bm{v}}} \ge 0$.

\subsection{The semi-discrete matrix system for DG}

Let ${\{\phi_i\}}_{i=1}^{N_p}$ be a nodal basis function for $P^N\left(\widehat{D}\right)$ located at Warp and Blend interpolation points \cite{hesthaven2007}. These basis functions  are defined implicitly using an orthogonal polynomial basis on the reference simplex. We define the reference mass matrix $\widehat{\bm{M}}$ and the physical mass matrix $\bm{M}$ for an element $D^k$ as
\[
\left(\widehat{\bm{M}}\right)_{ij}=\int_{\widehat{D}} \phi_j \phi_i~\text{d}\bm{x}, \qquad (\bm{M}_{ij})=\int_{D^k} \phi_j \phi_i~\text{d}\bm{x}=\int_{\widehat{D}} \phi_j \phi_i J~\text{d}\bm{\widehat{x}}.
\]
 $J$ is constant for affine mappings, and $\bm{M}=J\bm{\widehat{M}}$. We also define weak differentiation matrices $\bm{S}_{k}$ and face mass matrices $\bm{M}_f$ such that
\[
\left(\bm{S}_{k}\right)_{ij}=\int_{D^k} \myfrac{\partial \phi_j}{\partial \bm{x}_k} \phi_i~\text{d} \bm{x}, ~~~ (\bm{M}_f)_{ij}=\int_f \phi_j \phi_i \text{d}\bm{x}=\int_{\hat{f}} \phi_j \phi_i J^f \text{d}\hat{\bm{x}},
\]
where $J^f$ is the Jacobian of the mapping from the reference face $\widehat{f}$ to $f$. For affinely mapped simplices, $J^f$ is also constant and $M_f=J^f\widehat{\bm{M}}_f$, where the definition of the reference face mass matrix $\widehat{\bm{M}}_f$ is analogous to the definition of the reference mass matrix $\widehat{\bm{M}}$.

Finally, we introduce weighted mass matrices.  Let $w(\bm{x}) \in \mathbb{R}$ and $\bm{W(x)} \in \mathbb{R}^{m \times n}$. Then, scalar and matrix-weighted mass matrices $\bm{M}_w$ and $\bm{M_W}$ are defined as 
\begin{align}
(\bm{M}_w)_{ij}=\int_{D^k} w(\bm{x}) \phi_j(\bm{x}) \phi_i(\bm{x})~ \text{d}\bm{x},\qquad\bm{M_W}=\begin{bmatrix}\bm{M_{W_{1,1}}} & \dots & \bm{M_{W_{1,n}}}\\
\vdots &\ddots & \vdots\\
\bm{M_{W_{m,1}}} & \dots & \bm{M_{W_{m,n}}},
\end{bmatrix}
\label{eq:matweight}
\end{align}
where $\bm{M_{W_{i,j}}}$ is the scalar weighted mass matrix weighted by the $(i,j)^{\text{th}}$ element of $\bm{W}$.  Note that $\bm{M}_w, \bm{M}_{\bm{W}}$ are positive definite if $w(x), \bm{W}$ are pointwise positive definite.

Local contributions to the DG variational form may be evaluated in a quadrature-free manner using matrix-weighted mass matrices as defined above. Let $\bm{\Sigma_i},~\bm{V}_i$ denote vectors containing degrees of freedom for solutions components $\bm{\sigma}_i$ and $\bm{v}_i$, such that
\begin{align*}
\bm{v}_i (\bm{x}, t)&=\sum_{j=1}^{N_p} (\bm{V}_i(t))_j\phi_j(\bm{x}), \qquad\qquad1 \le i \le 6 \\
\bm{\sigma}_i (\bm{x}, t)&=\sum_{j=1}^{N_p} (\bm{\Sigma_i(t)})_j \phi_j(\bm{x}), \qquad\qquad1 \le i \le 7 
\end{align*}
Then, the local DG formulation can be written as a block system of ordinary differential equations by concatenating $\bm{\Sigma}_i,~\bm{V}_i$ into single vectors $\bm{\Sigma}$ and $\bm{V}$ and using the Kronecker product $\otimes$
\begin{align} 
\bm{M}_{\bm{Q}_s^{-1}} \myfrac{\partial \bm{\Sigma}}{\partial t} &= \sum_{i=1}^{d} \left(\bm{A}_i \otimes \bm{S_i}\right)\bm{V} + \sum_{f \in \partial D^k} \left(\bm{I} \otimes \bm{M}_f\right) \bm{F}_\sigma +  \bm{M}_{\bm{S}}\bm{\Sigma}  \label{eqsd1}\\
\bm{M}_\rho \myfrac{\partial \bm{V}}{\partial t}&=\sum_{i=1}^{d}\left({\bm{A}_i}^T \otimes \bm{S_i}\right)\bm{\Sigma} + \sum_{f \in \partial D^k} \left(\bm{I} \otimes \bm{M}_f\right) \bm{F}_v \label{eqsd2},
\end{align}
where $\bm{F}_v$ and $\bm{F}_\sigma$ denote the degrees of freedom for the velocity and stress numerical fluxes.

In order to apply a time integrator, we must invert $\bm{M}_{\bm{Q}_s}$ and $\bm{M}_{\bm{\rho}}$.  While the inversion of $\bm{M}_{\bm{Q}_s^{-1}}$ and  $\bm{M}_{\bm{\rho}}$. can be parallelized from element to element, doing so typically requires either the precomputation and storage of the dense matrix inverses or on-the-fly construction and solution of a large dense matrix system at every time step. The former option requires a large amount of storage, while the latter option is computationally expensive and difficult to parallelize among degrees of freedom. This cost can be avoided when ${\bm{Q}_s^{-1}}$ and $\rho$ are constant over an element $D^k$. In this case, $\bm{M}_{\bm{Q}_s}$ reduces to
\begin{align*} 
\bm{M}_{\bm{Q}_s^{-1}}^{-1}=\begin{bmatrix}{\bm{Q}_s^{-1}}_{(1,1)}\bm{M} & \dots & {\bm{Q}_s^{-1}}_{(1,N_d)}\bm{M}\\
\vdots &\ddots & \vdots\\
{\bm{Q}_s^{-1}}_{(N_d,1)}\bm{M} & \dots & {\bm{Q}_s^{-1}}_{(N_d,N_d)}\bm{M}
\end{bmatrix}^{-1} =\left({\bm{Q}_s} \otimes \bm{M}  \right)^{-1}=\bm{Q}_s \otimes \left(\myfrac{1}{J} \bm{\widehat{M}^{-1}}\right).
\end{align*} 
Similarly ${\bm{M}^{-1}_{\bm{S}}}$  and $\bm{M}^{-1}_\rho$ can be expressed as $\bm{M}_{\bm{S}}={\bm{S}} \otimes \left(\myfrac{1}{J} \bm{\widehat{M}^{-1}}\right)$  and $\bm{M}^{-1}_{\rho}=\rho^{-1} \otimes \left(\myfrac{1}{J} \bm{\widehat{M}^{-1}}\right)$, and respectively.  Applying these observations to (\ref{eqsd1}) and (\ref{eqsd2}) yields the following sets of local ODEs over each element
\begin{align}
\myfrac{\partial \bm{\Sigma}}{\partial t} &= \sum_{i=1}^{d} \left({\bm{Q}_s}\bm{A} \otimes \bm{D_i}\right)\bm{V} + \sum_{f \in \partial D^k} \left({\bm{Q}_s} \otimes \bm{M}_f\right) \bm{F}_\sigma + {\bm{M}^{-1}_{\bm{Q}_s^{-1}}}\bm{M}_{\bm{D}}\bm{\Sigma},,\\
\myfrac{\partial \bm{V}}{\partial t}&=\sum_{i=1}^{d}\left({\rho} ^{-1}{\bm{A}}^T \otimes \bm{D_i}\right)\bm{\Sigma} + \sum_{f \in \partial D^k} \left({\rho}^{-1} \otimes \bm{M}_f\right) \bm{F}_v 
\end{align}
where we have introduced the differentiation matrix $\bm{D}_i=\bm{M}{^{-1}}\bm{S}_i$ and lift matrix $\bm{L}_f=\bm{M}^{-1}\bm{M}_f$. For affine elements, both derivative and lift matrices are applied using products of geometric factors and reference derivative and lift matrices.

If $\bm{Q}_s^{-1}$ and $\rho$ vary spatially within the element, then the above approach can no longer be used to invert $\bm{Q}_s$ and $\bm{Q}_v$.  This case can be treated using the approach of \cite{chan2018}, where $\bm{M}_{\bm{Q}_s},{\rho}$ are replaced with  weight-adjusted approximations. These approximations are low storage, simple to invert, and yield an energy stable and high order accurate DG method to approximate the matrix-weighted $L^2$ inner product (and corresponding matrix-weighted mass-matrices $\bm{Q}_s$ and $\rho$).  

We also note that, material coefficients $\bm{Q_s}^{-1},~\rho$ appear only on the left hand side of (\ref{eq24}). The right hand side of (\ref{eq24}) is equivalent to the discretization of a constant coefficient system.  This provides additional advantages in that the right hand side can be evaluated using efficient techniques for DG discretizations of constant-coefficient problems \cite{chan2017gpu, guo2018}.

\begin{figure}
	\centering
	\subfloat[$\alpha=0$, central flux]{
		\includegraphics[width=0.5\textwidth]{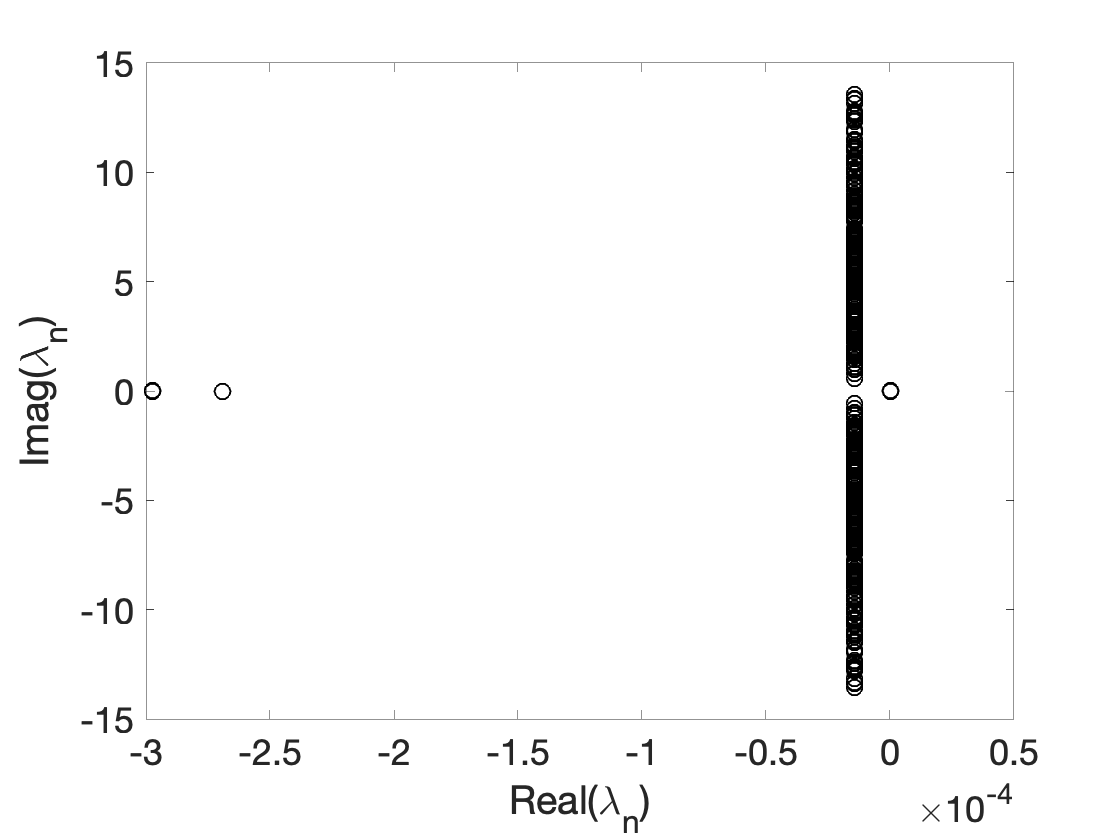}
	}
	\subfloat[$\alpha=1$, penalty flux]{
			\includegraphics[width=0.5\textwidth]{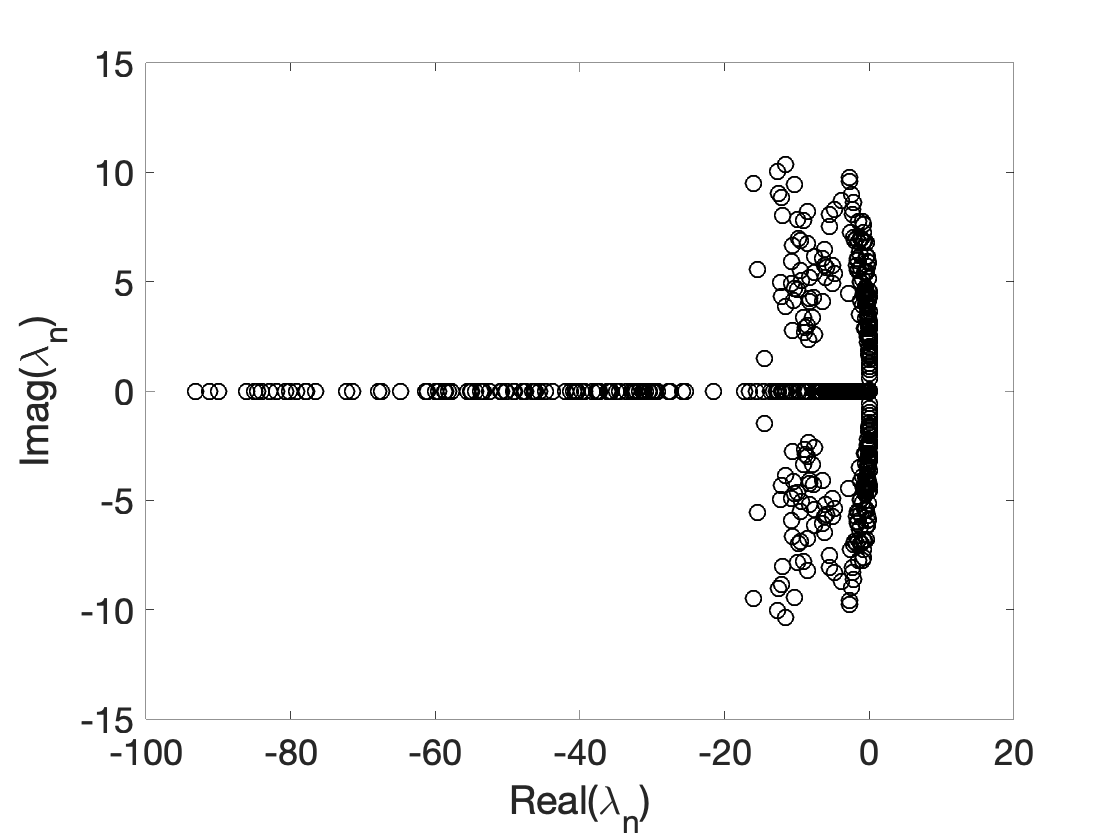}
	}
	\caption{Spectra for $N = 3$ and $h = 1/2$ with a material property of isotropic Sandstone (Column 3 of Table 1).  For  $\alpha_{\tau}=\alpha_{v}=0$  and $\alpha_{\tau}=\alpha_{v}=1$ , the largest real part of spectra are $1.83519\text{e}-14$ and $2.18232\text{e}-14$, respectively. }
\end{figure}
\begin{table}
	\caption{Material properties of anisotropic-viscoelastic media \cite{carcione1995} }
	\centering
	\medskip
	\begin{tabular}{c c c  c} 
		\hline
		Properties  & Clay shale  & Phenolic  & Isotropic Sandstone\\ [0.5ex] 
		\hline
		\hspace*{-2cm} \textbf{Elasticities} & &  & \\
		$\rho_s~$ $($kg/m$^3$$)$ & 2590\ & 1364 & 2500 \\
		$c_{11}~$(GPa) & 66.6 &  11.7 &  25.6  \\ 
		$c_{12}~$(GPa) & 19.7  & 6.7 & 9.4 \\ 
		$c_{13}~$(GPa) & 39.4  & 7.0  &  9.4\\ 
		$c_{22}~$(GPa) & 66.6 &  15.4 & 25.6 \\ 
		$c_{23}~$(GPa) & 39.4  & 7.0 & 9.4\\ 
		$c_{33}~$(GPa) & 39.9 &  17.4 & 25.6 \\ 
		$c_{44}~$(GPa) & 10.9 & 3.8 &  16.2\\ 
		$c_{55}~$(GPa) & 10.9 & 3.5 & 16.2\\ 
		$c_{66}~$(GPa) & 23.4 & 3.1 & 16.2\\ 
		\hspace*{-2cm} \textbf{Relaxation time~(s)} & ~ & ~ \\
		$\tau_\epsilon^{(1)}$ & $8.00 \times 10^{-3}$ & $6.4 \times 10^{-3}$  & $3.72 \times 10^-3$ \\
		$\tau_\sigma^{(1)}$ & $7.49 \times 10^{-3}$ & $6.00 \times 10^{-3}$  & $3.36 \times 10^-3$ \\
		$\tau_\epsilon^{(2)}$ & $8.00 \times 10^{-3}$ & $6.4 \times 10^{-3}$ & $3.78 \times 10^-3$  \\
		$\tau_\sigma^{(2)}$ & $7.25 \times 10^{-3}$ & $5.80 \times 10^{-3}$  & $3.30 \times 10^-3$ \\
		$\tau_\epsilon^{(3)}$ & $8.00 \times 10^{-3} $ & $6.4 \times 10^{-3}$ & $3.78 \times 10^-3$  \\
		$\tau_\sigma^{(3)}$ & $7.25\times 10^{-3}$ & $5.60 \times 10^{-3}$ & $3.30 \times 10^-3$   \\
		$\tau_\epsilon^{(4)}$ & $8.00 \times 10^{-3}$ & $6.4 \times 10^{-3}$ & $3.78 \times 10^-3$  \\
		$\tau_\sigma^{(4)}$ & $7.25 \times 10^{-3}$ & $5.30 \times 10^{-3}$ & $3.30 \times 10^-3$   \\
		\hline
	\end{tabular}
	\label{tab:props}
\end{table}

\section{Numerical experiments}
\label{numex}
In this section, we present several numerical experiments to validate the stability and accuracy of the proposed method in two and three dimensions. The convergence of the new DG formulation in piecewise constant isotropic viscoelastic media is confirmed. Finally, the method is applied to solve the anisotropic viscoelastic wave equation in various heterogeous media.

Time integration is performed using the low-storage 4$^\text{th}$ order five-stage Runge-Kutta scheme of Carpenter and Kennedy \cite{carpenter1994fourth}, and the time step is chosen based on the global estimate
\begin{align}
\label{eqdt}
dt=\min_k\myfrac{C_{CFL}}{ \max{\left(\lambda_i\right)} C_N \left \Vert J^f \right \Vert_{L^{\infty}(\partial D^k)} \left \Vert J^{-1} \right \Vert_{L^{\infty}(D^k)} }
\end{align}
where $\lambda_i$ are wave speeds of the system \cite{carcione2002}, $C_N=O(N^2)$ is the order-dependent constant in the surface polynomial trace inequality \cite{chan2016gpu}, and $C_{\rm CFL}$ is a tunable global CFL constant. This estimate is derived by bounding the eigenvalues of the spatial DG discretization matrix appearing in the semi-discrete system of ODEs.  This choice of $dt$ is very conservative as it is derived based on an upper bound on the spectral radius.

\subsection{Spectra and choice of penalty parameter}
We first verify the energy stability of proposed DG formulation. Let $\bm{A}_h$ denote the matrix induced by the global semi-discrete DG formulation, such that time evolution of the solution $\bm{q}$ is governed by 
\[
\myfrac{\partial \bm{q}}{\partial t}=\bm{A}_h \bm{q},
\]
where $\bm{q}$ denotes a vector of degrees of freedom. We show in Figure 1 eigenvalues of $\bm{A}_h$ for $\alpha=0$ and $\alpha=1$ with material parameters of isotropic sandstone (given in Column 3 of Table 1). The discretization parameters are $N=3$ and $h=1/2$. In both cases, the largest real part of any eigenvalues is $O(10^{-14})$, which suggests that the semi-discrete scheme is indeed energy stable. It is to be also noted that some eigenvalues for $\alpha=0$ have purely negative real part, corresponding to the dissipation present in viscoelastic system. 

For practical simulations, the choice of the penalty parameter $\alpha$ remains to be specified. Taking $\alpha >0$ results in damping of under-resolved spurious components of the solutions. However, a naive selection of $\alpha$ can result in an overly restrictive time-step restriction for stability.  A guiding principle for determining appropriate values of the penalty parameters $\alpha$ is to ensure that the spectral radius is the same magnitude as the case when $\alpha= 0$.   For example, the spectral radius of $\bm{A}_h$, $\rho(\bm{A}_h)$ is $13.5653$ for $\alpha=0$ which is $O(N^2/h)$.  The spectral radius $\rho(\bm{A}_h)$ is $45.6388$ for $\alpha_{\bm{\tau}},\alpha_{\bm{v}}=0.5$, while the spectral radius for $\alpha_{\bm{\tau}},\alpha_{\bm{v}}=1$ is $\rho(\bm{A}_h) = 93.1184$.  Since the maximum stable time step is proportional to the spectral radius, taking $\alpha=1$ in this case results in a more restrictive CFL condition.  This phenomena is related to observations in \cite{chan2017penalty} that large penalty parameters result in extremal eigenvalues of $\bm{A}_h$ with very large negative real parts. 
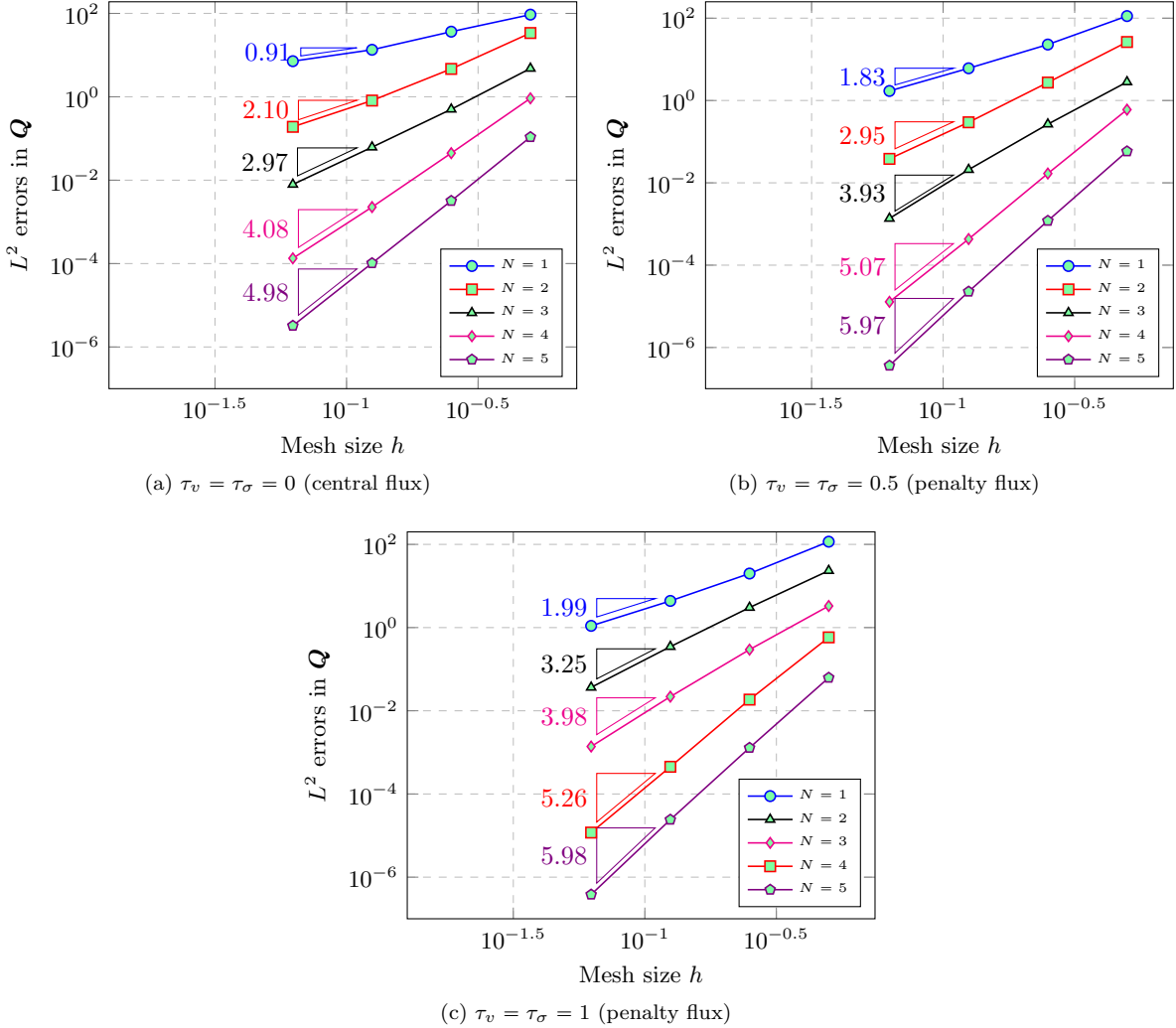
\begin{figure}
	\centering
	\subfloat[$\tau_v=\tau_\sigma=0$ (central flux)]{
		\begin{tikzpicture}
		\begin{loglogaxis}[
		legend cell align=left,
		width=.475\textwidth,
		xlabel={Mesh size $h$},
		ylabel={$L^2$ errors in $\bm{Q}$}, 
		xmin=.0125, xmax=0.75,
		ymin=1e-7, ymax=190,
		legend pos=south east, legend cell align=left, legend style={font=\tiny},
		xmajorgrids=true,
		ymajorgrids=true,
		grid style=dashed,
		] 	
		\addplot[color=blue,mark=*,semithick, mark options={fill=markercolor}]
		coordinates{(0.5,93.416652199577868)(0.25,36.461055091130305)(0.125,13.419807385475686)(0.0625,7.144012450937141)};
		\logLogSlopeTriangleFlip{0.53}{0.120}{0.860}{0.91}{blue}
		
		\addplot[color=red,mark=square*,semithick, mark options={fill=markercolor}]
		coordinates{(0.5,33.847322601986292)(0.25,4.664858070496184)(0.125,0.812481931025341)(0.0625, 0.190352686003959)};
		\logLogSlopeTriangleFlip{0.53}{0.125}{0.695}{2.10}{red}
		
		\addplot[color=black,mark=triangle*,semithick, mark options={fill=markercolor}]
		coordinates{(0.5, 4.786112396198761)(0.25,0.501918675729696)(0.125,0.061832282648908)(0.0625, 0.007904911009118)};
		\logLogSlopeTriangleFlip{0.53}{0.127}{0.550}{2.97}{black}
		
		\addplot[color=magenta,mark=diamond*,semithick, mark options={fill=markercolor}]
		coordinates{(0.5,0.919465420828410)(0.25, 0.044526272950184)(0.125,0.002266723445930)(0.0625,0.000133825573583)};
		\logLogSlopeTriangleFlip{0.53}{0.125}{0.365}{4.08}{magenta}
		
		\addplot[color=violet,mark=pentagon*,semithick, mark options={fill=markercolor}]
		coordinates{(0.5,0.108224424061850)(0.25,0.003193194110986)(0.125,0.000102373940223)(0.0625,0.000003234259164)};
		\logLogSlopeTriangleFlip{0.53}{0.125}{0.190}{ 4.98}{violet}
		
		\legend{$N=1$,$N=2$,$N=3$,$N=4$,$N=5$}
		\end{loglogaxis}
		\end{tikzpicture}
	}
	\subfloat[$\tau_v=\tau_\sigma=0.5$ (penalty flux)]{
		\begin{tikzpicture}
		\begin{loglogaxis}[
		legend cell align=left,
		width=.475\textwidth,
		xlabel={Mesh size $h$},
		ylabel={$L^2$ errors in $\bm{Q}$}, 
		xmin=.0125, xmax=0.75,
		ymin=1e-7, ymax=250,
		legend pos=south east, legend cell align=left, legend style={font=\tiny},
		xmajorgrids=true,
		ymajorgrids=true,
		grid style=dashed,
		] 
		
		\addplot[color=blue,mark=*,semithick, mark options={fill=markercolor}]
		coordinates{(0.5,1.124130635617782e+02)(0.25,0.228703290452522e+02)(0.125,0.060579644716072e+02)(0.0625,0.017040904280704e+02)};
		\logLogSlopeTriangleFlip{0.53}{0.125}{0.785}{1.83}{blue}
		
		\addplot[color=red,mark=square*,semithick, mark options={fill=markercolor}]
		coordinates{(0.5,0.263619359948516e+02)(0.25,0.027293429824478e+02)(0.125,0.002953557596423e+02)(0.0625, 0.000382111919095e+02)};
		\logLogSlopeTriangleFlip{0.53}{0.125}{0.620}{2.95}{red}
		
		\addplot[color=black,mark=triangle*,semithick, mark options={fill=markercolor}]
		coordinates{(0.5, 0.028167350684454e+02)(0.25,0.002656635936889e+02)(0.125,0.000208832149096 e+02)(0.0625, 0.000013625036149e+02)};
		\logLogSlopeTriangleFlip{0.53}{0.125}{0.459}{3.93}{black}
		
		\addplot[color=magenta,mark=diamond*,semithick, mark options={fill=markercolor}]
		coordinates{(0.5,0.005993972510410e+02)(0.25, 0.000168206603431e+02)(0.125,0.000004305361829e+02)(0.0625,0.000000128659892e+02)};
		\logLogSlopeTriangleFlip{0.53}{0.125}{0.255}{5.07}{magenta}
		
		\addplot[color=violet,mark=pentagon*,semithick, mark options={fill=markercolor}]
		coordinates{(0.5,0.000582583540833e+02)(0.25,0.000011961312972e+02)(0.125,0.000000228088767e+02)(0.0625,0.000000003632606e+02)};
		\logLogSlopeTriangleFlip{0.53}{0.125}{0.092}{ 5.97}{violet}
		
		\legend{$N=1$,$N=2$,$N=3$,$N=4$,$N=5$}
		\end{loglogaxis}
		\end{tikzpicture}
	}\\
	\subfloat[$\tau_v=\tau_\sigma=1$ (penalty flux)]{
		\begin{tikzpicture}
		\begin{loglogaxis}[
		legend cell align=left,
		width=.475\textwidth,
		xlabel={Mesh size $h$},
		ylabel={$L^2$ errors in $\bm{Q}$}, 
		xmin=.0125, xmax=0.75,
		ymin=1e-7, ymax=200,
		legend pos=south east, legend cell align=left, legend style={font=\tiny},
		xmajorgrids=true,
		ymajorgrids=true,
		grid style=dashed,
		] 
		
		\addplot[color=blue,mark=*,semithick, mark options={fill=markercolor}]
		coordinates{(0.5,1.160878895681427e+02)(0.25,0.199202420066448e+02)(0.125,0.043754293552746e+02)(0.0625,0.011005364295797e+02)};
		\logLogSlopeTriangleFlip{0.53}{0.125}{0.780}{1.99}{blue}
		
		\addplot[color=black,mark=triangle*,semithick, mark options={fill=markercolor}]
		coordinates{(0.5, 0.233367197675592e+02)(0.25,0.030323878949011e+02)(0.125,0.003478094120710e+02)(0.0625, 0.00036516082598e+02)};
		\logLogSlopeTriangleFlip{0.53}{0.125}{0.620}{3.25}{black}
		
		\addplot[color=magenta,mark=diamond*,semithick, mark options={fill=markercolor}]
		coordinates{(0.5,0.033056519318894 e+02)(0.25, 0.002953163077176 e+02)(0.125,0.000218486038857e+02)(0.0625,0.000013823345125e+02)};
		\logLogSlopeTriangleFlip{0.53}{0.125}{0.476}{3.98}{magenta}
		\addplot[color=red,mark=square*,semithick, mark options={fill=markercolor}]
		coordinates{(0.5,0.005815901969312e+02)(0.25,0.000185458445416e+02)(0.125,0.000004503632810e+02)(0.0625, 0.000000117994542e+02)};
		\logLogSlopeTriangleFlip{0.53}{0.125}{0.250}{5.26}{red}
		
		\addplot[color=violet,mark=pentagon*,semithick, mark options={fill=markercolor}]
		coordinates{(0.5,0.000623845962079e+02)(0.25,0.000012886724665e+02)(0.125,0.000000241280057 e+02)(0.0625,0.000000003823192e+02)};
		\logLogSlopeTriangleFlip{0.53}{0.125}{0.092}{ 5.98}{violet}

		\legend{$N=1$,$N=2$,$N=3$,$N=4$,$N=5$}
		\end{loglogaxis}
		\end{tikzpicture}
	}
	\caption{Convergence of $L^2$ error for plane wave in a viscoelastic media }
\end{figure}

\subsection{Convergence for a plane wave in viscoelastic medium}
The analytical solution to (\ref{eq23}) for a plane wave is given as
\begin{align}
\label{eq25}
\bm{q}_n(\bm{x},t)=\bm{q}_n^0 \exp[\text{i} \cdot (\omega t- \bm{k}\cdot \bm{x})], \qquad n=1...15,
\end{align}
where $\bm{Q}_n^0$ is the initial amplitude vector of stress and velocity components; $\omega$ are wave frequencies; $\bm{k}=(k_x, k_y, k_z)$ is the wave-number vector.
To achieve realistic viscoelastic behavior, we superimpose three plane waves, of the form given by (\ref{eq25}), corresponding to a  P-wave, and S-wave.

Now, we briefly describe how we determine the wave frequencies $\omega$. Substituting (\ref{eq25}) into (\ref{eq23}) yields
\begin{align}
\label{eq26}
\omega \bm{q}_n^0=(\bm{A} (\bm{x})k_x + \bm{B}(\bm{x}) k_y + \bm{C}(\bm{x}) k_z -\text{i} \bm{D}(\bm{x}))\bm{q}_n^0
\end{align}

Solving the three eigenvalues problem for in (\ref{eq26}) for each wave mode $l$ yields in matrix of right eigenvectors $(R^{(l)}_{mn})$ and eigenvalue $(\omega_{l})$.  Following \cite{toro2009,de2008}, the solution  of (\ref{eq23}) for the plane wave can be constructed as
\begin{align}
\label{eq33}
\bm{q}_n(\bm{x},t)=\sum_{l=1}^2 R^{(l)}_{mn} \gamma_{n}^{(l)} \exp[\text{i} \cdot (\omega^{(l)} t- \bm{k^{(l)}}\cdot \bm{x})],
\end{align}
where $\gamma_{n}^{(l)}$ is a amplitude coefficient with $\gamma_{1}^{(1)}=\gamma_{2}^{(2)}=1$.

Next, we study the accuracy and convergence of our DG method for a plane wave propagating in an isotropic porous sandstone with material properties given in Table 1 (Column 3). Unless otherwise stated, we report relative $L^2$ errors for all components of the solution $\bm{U}$
	\[
	\myfrac{\left \Vert \bm{U} -\ \bm{U}_h \right \Vert_{L^2(\Omega)}}{\left \Vert \bm{U} \right \Vert_{L^2(\Omega)}}=\myfrac{\left( \sum_{i=1}^m \left \Vert \bm{U_i} -\bm{U}_{\bm{i},\bm{h}} \right \Vert ^2_{L^2(\Omega)}\right)^{1/2}}{\left(\sum_{i=1}^m \left \Vert \bm{U}_i \right \Vert^2_{L^2(\Omega)} \right)^{1/2}}.
	\]
The error is computed for an  isotropic sandstone.  In Figure 2, we show the $L^2$ errors computed at T=1, using uniform triangular meshes constructed by bisecting an uniform mesh of quadrilaterals along the diagonal. Figure 2a shows error plots using the central flux with $\alpha=0~(\text{CFL}=1)$. We observe convergence rate of $O(h^{N})$ or $O(h^{N+1})$ for odd-even $N$. Figure 2b and 2c show errors for penalty fluxes with  $\alpha=0.5~(\text{CFL}=0.5),~\text{and}~1~(\text{CFL}=0.8)$, respectively. For $N=1,...,5$, $O(h^{N+1})$, rates of convergence are observed. We note that for $N=4$ and $N=5$, we observe results which are better than the $4^{\text{th}}$ order accuracy of our time-stepping scheme. This is most likely due to the benign nature of the solution in time and the choice of time step (\ref{eqdt}) which scales as $O(h/N^2)$. 
\begin{figure}
	\centering
	\subfloat[Orthotropic shale, $v_1$]{
		\begin{tikzpicture}
		\begin{axis}[enlargelimits=false, axis on top, axis equal image, xlabel=x (km), ylabel=y (km), width=0.6\textwidth]
		
		\addplot graphics [xmin=0,xmax=2,ymin=0,ymax=2] {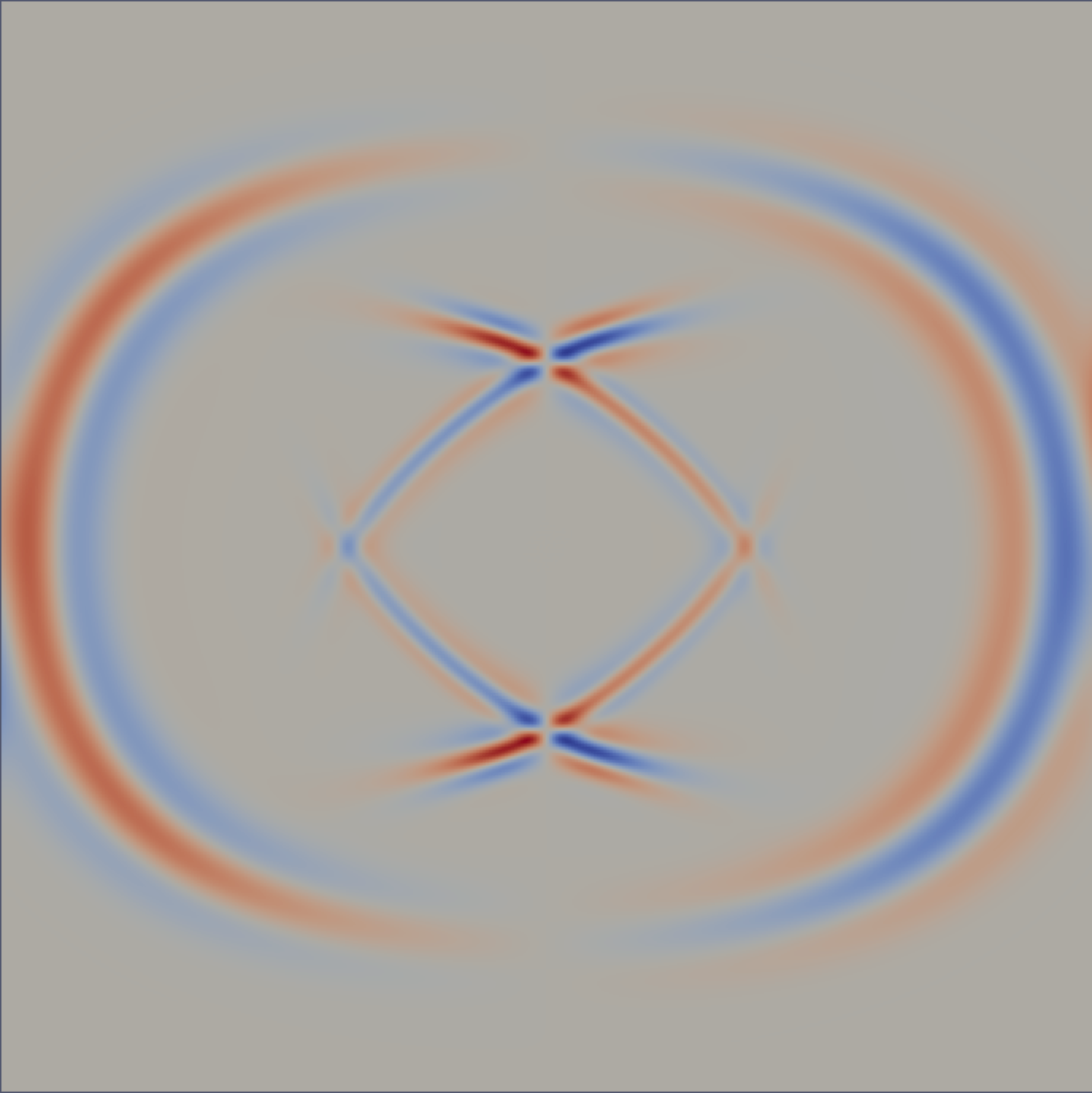};
		\end{axis}
		\end{tikzpicture}
	}
	\subfloat[Orthotropic shale, $v_3$]{
		\begin{tikzpicture}
		\begin{axis}[enlargelimits=false, axis on top, axis equal image, xlabel=x (km), ylabel=y (km), width=0.6\textwidth]
		
		\addplot graphics [xmin=0,xmax=2,ymin=0,ymax=2] {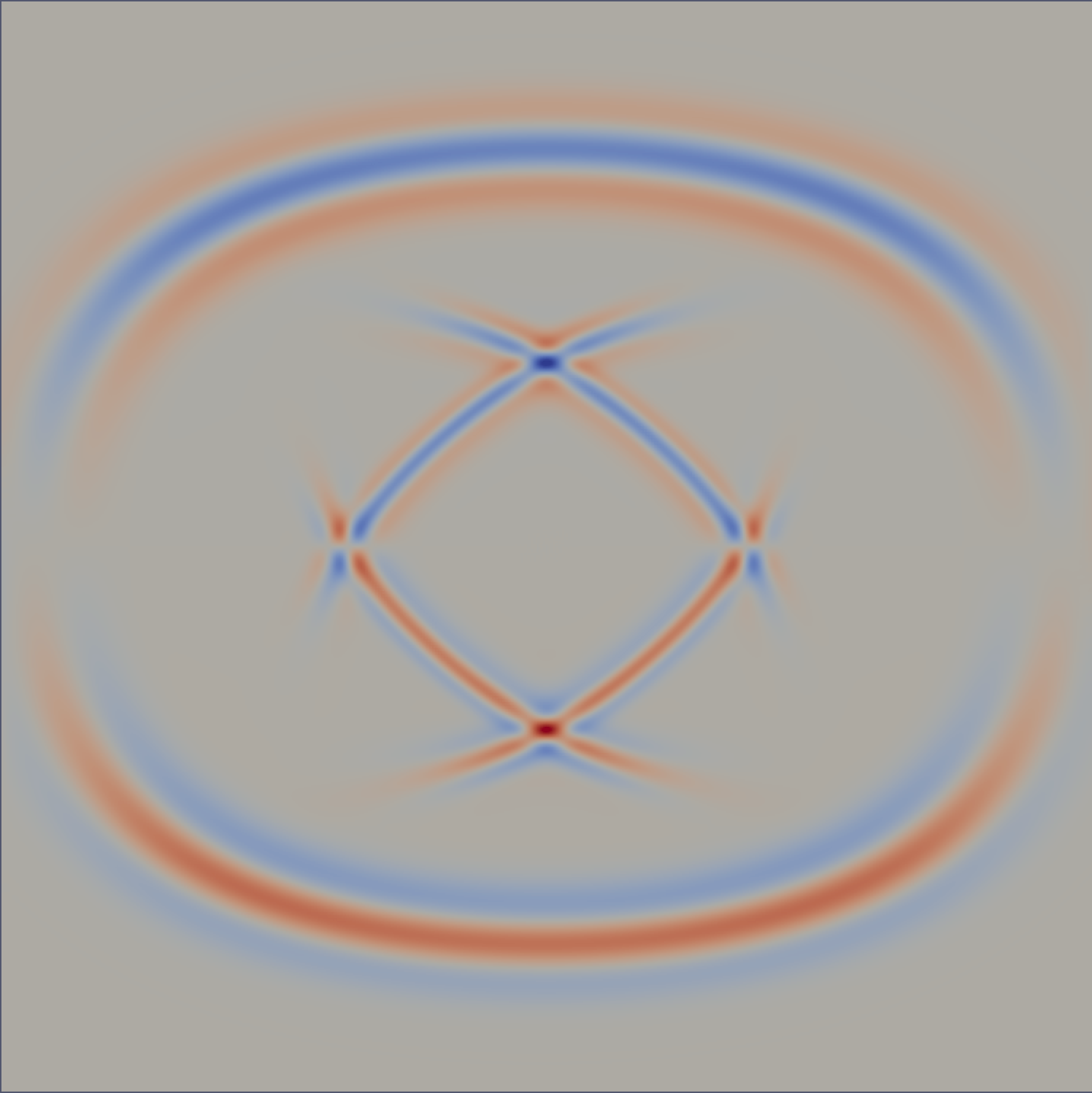};
		\end{axis}
		\end{tikzpicture}
	}\\
	\caption{Snapshots of  particle velocities in orthotropic shale (column 1 of Table 1 ), computed at $t=0.22~\text{s}$, where (a) and (b) corresponds to $v_1$ and $v_3$ components. The central frequency of the forcing function is  20 Hz corresponding to the relaxation peak of materials. The point source is located in the center of the domain. The solution is computed using polynomials of degree $N=3$ and $K=32,768$. }
\end{figure}
\subsection{Application examples}
We next demonstrate the accuracy and flexibility of the proposed DG method for several application-based problems in linear viscoelasticity with anisotropy.  All computations are done using penalty parameters $\alpha=0.5$ unless specified otherwise. In subsequent simulations, the forcing is applied to both the  $x-$ nd $z-$ components of stress i.e. $(\sigma_{11}, \sigma_{33})$.
\begin{align}
\label{Ricker}
f(\bm{x}, t)=(1-2(\pi f_0(t-t_0))^2)\exp[-(\pi f_0(t-t_0))^2] \delta(\bm{x}-\bm{x}_0),
\end{align}
where $\bm{x}_0$ is the position of the point source and $f_0$ is the central frequency.

In the following simulations, two types of elastic waves are observed: a  P wave, and an S wave with an anisotropic dissipative phenomena.
\begin{figure}
	\centering
	\subfloat[Phenolic type material, $v_1$]{
		\begin{tikzpicture}
		\begin{axis}[enlargelimits=false, axis on top, axis equal image, xlabel=x (cm), ylabel=y (cm), width=0.6\textwidth]
		
		\addplot graphics [xmin=0,xmax=40,ymin=0,ymax=40] {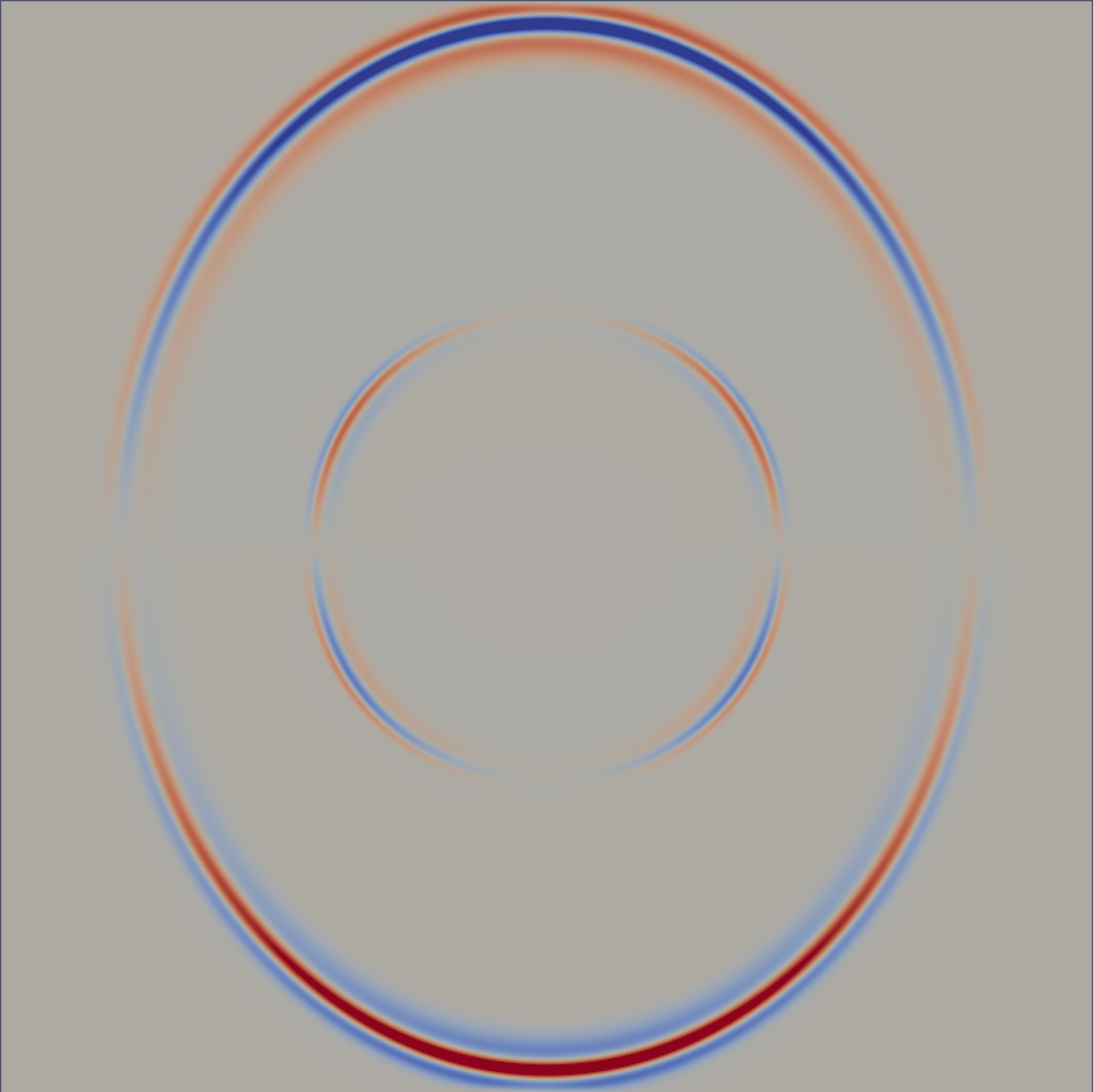};
		\end{axis}
		\end{tikzpicture}
	}
	\subfloat[Phenolic type material, $v_2$]{
		\begin{tikzpicture}
		\begin{axis}[enlargelimits=false, axis on top, axis equal image, xlabel=x (cm), ylabel=y (cm), width=0.6\textwidth]
		\addplot graphics [xmin=0,xmax=40,ymin=0,ymax=40] {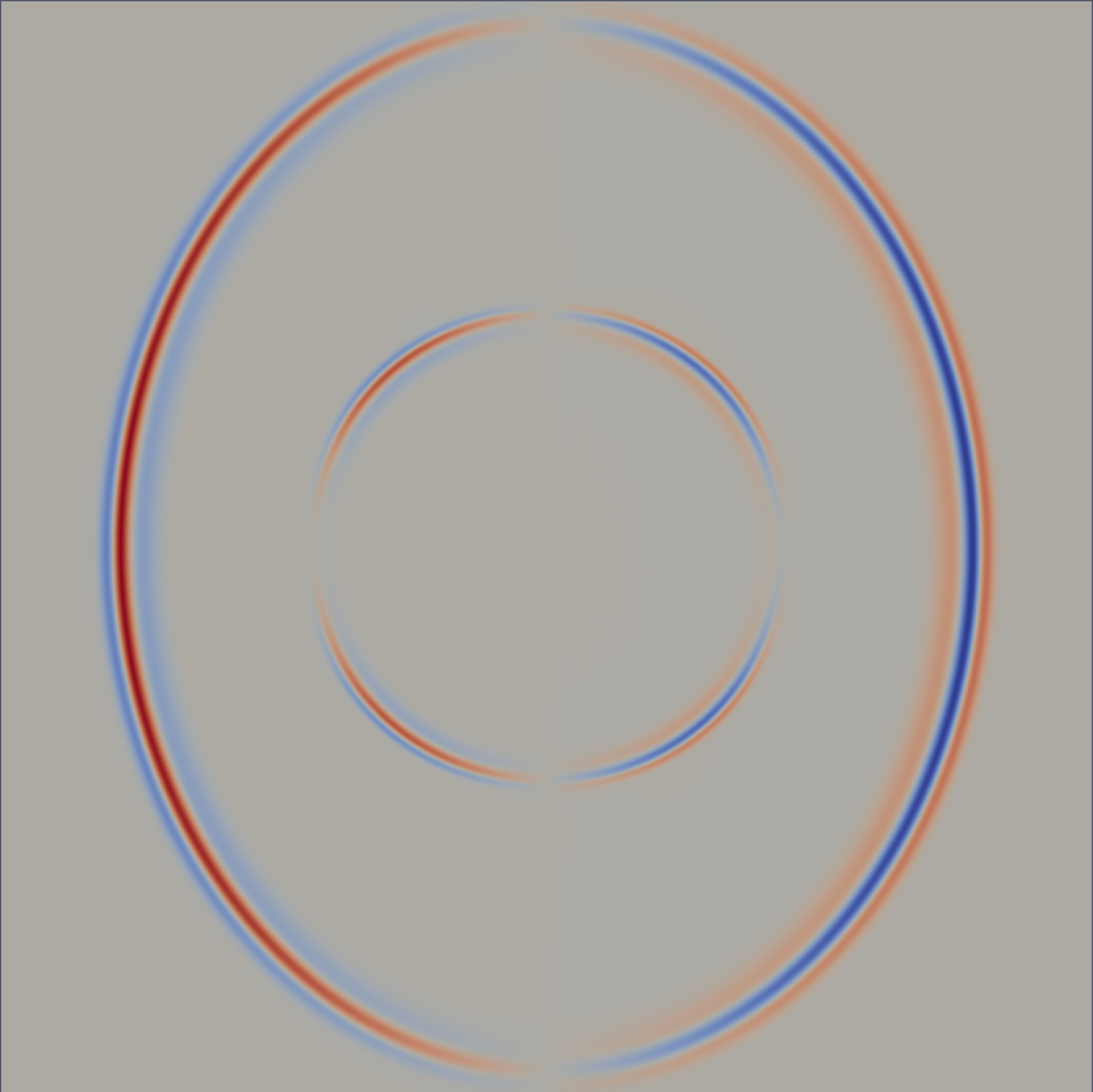};
		\end{axis}
		\end{tikzpicture}
	}
	\caption{Snapshots of  particle velocities in Phenolic material (column 2 of Table 1), computed at $t=53.2~\mu$\text{s}, where (a) and (b) corresponds to $v_1$ and $v_2$ components. The central frequency of the forcing function is  250 kHz which also corresponds to relaxation peak of the material. The point source is located in the center of the domain. The solution is computed using polynomials of degree $N=3$ and $K=131,072$.}
\end{figure}
\subsubsection{2D Orthotropic shale}
To illustrate the effect of anisotropic dissipation in a viscoelastic medium, we perform a computational experiment in orthotropic shale with material properties given in Table 1 (Column 1).  The size of the computational domain is $2~\text{km} \times 2~\text{km} $. The domain is discretized with uniform triangular elements with a minimum edge length of $15.625~\text{m}$. Figures 3(a)-(b) show the $x-$ and $z-$ components of the  particle velocity of the orthotropic shale, respectively. The central frequency of the forcing function is $f_0=20~\text{Hz}$, which is also the frequency for relaxation peak of the material. Polynomials of degree $N=3$ are used for the simulation, and the propagation time is $0.22~\text{s}$. Both wave modes can be observed: the  P mode  and the shear mode (S, inner wavefront). A shear wave cusp is clearly observed in Figure 3a and 3b.
\begin{figure}
	\centering
	\subfloat[$v_1$ for elastic and vosocelastic approximation]{
		\includegraphics[width=0.5\textwidth]{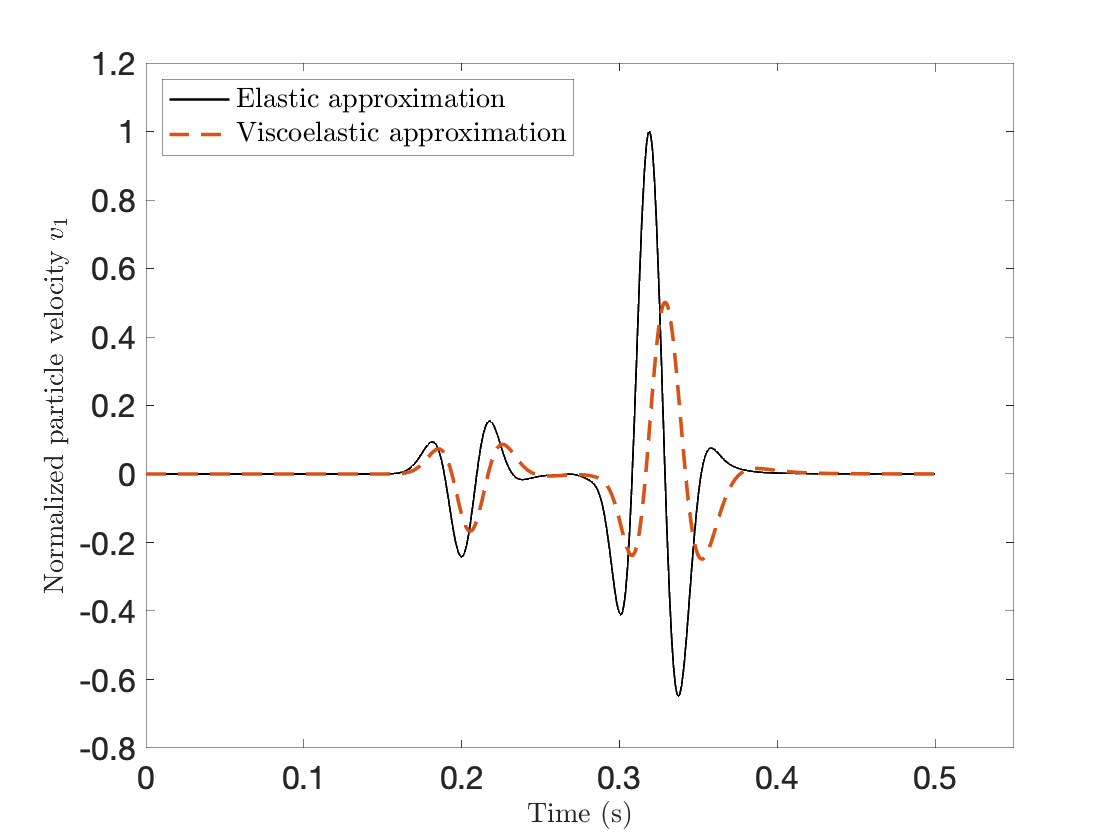}
	}
	\subfloat[$v_3$ for elastic and vosocelastic approximation]{
		\includegraphics[width=0.5\textwidth]{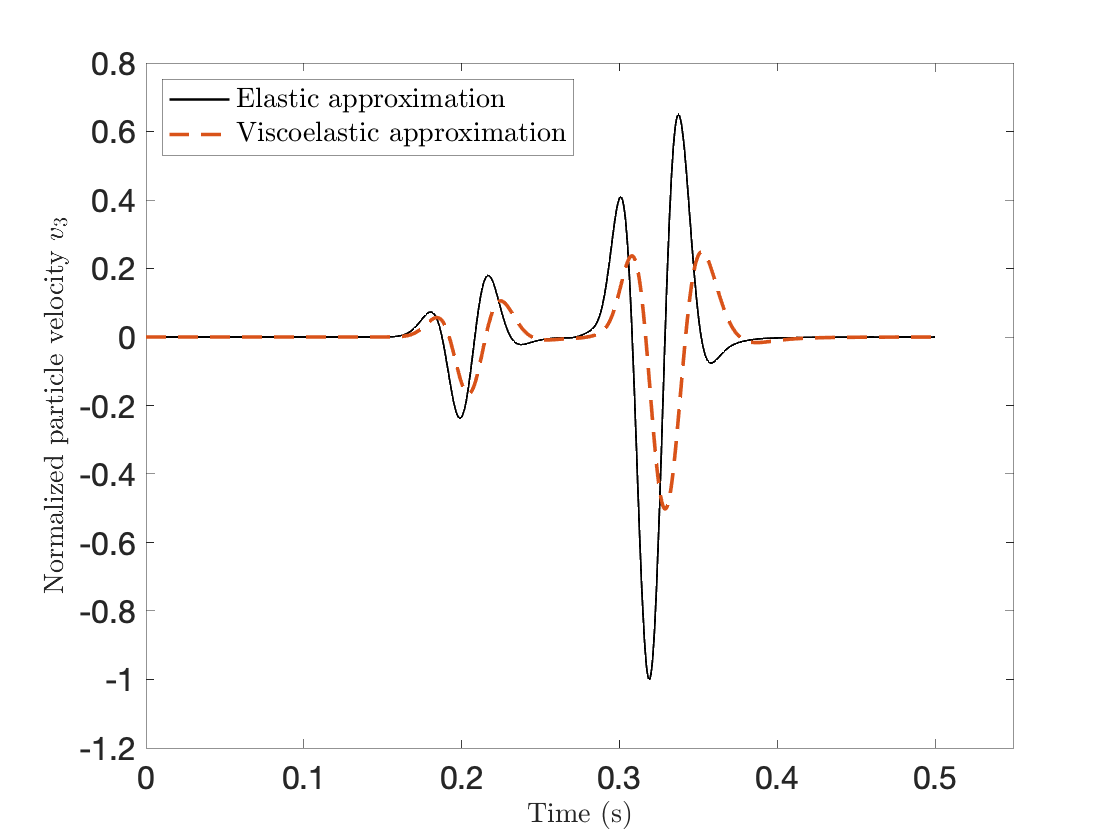}
	}
	\caption{A comparison of the time history of particle velocities for both elastic and viscoelastic  approximations. Subfigures (a) and (b) represent the horizontal and vertical particle velocities, respectively. A difference between elastic and viscoelastic approximation is clearly visible in both by phase and amplitude difference between the traces.}
	
\end{figure}
\begin{figure}
	\centering
	\subfloat[Analytical vs numerical simulation for $v_1$]{
		\includegraphics[width=0.5\textwidth]{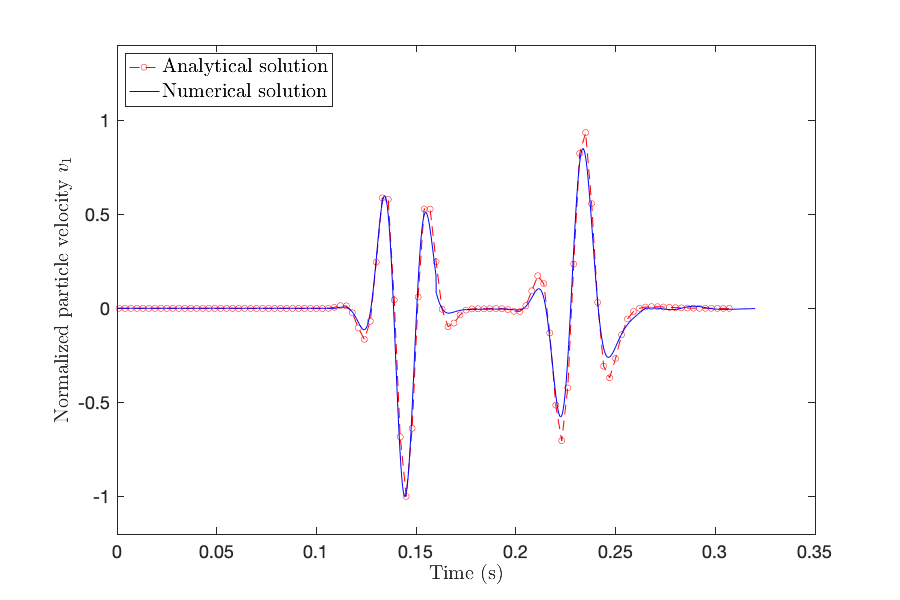}
	}
	\subfloat[Analytical vs numerical simulation for $v_2$]{
		\includegraphics[width=0.5\textwidth]{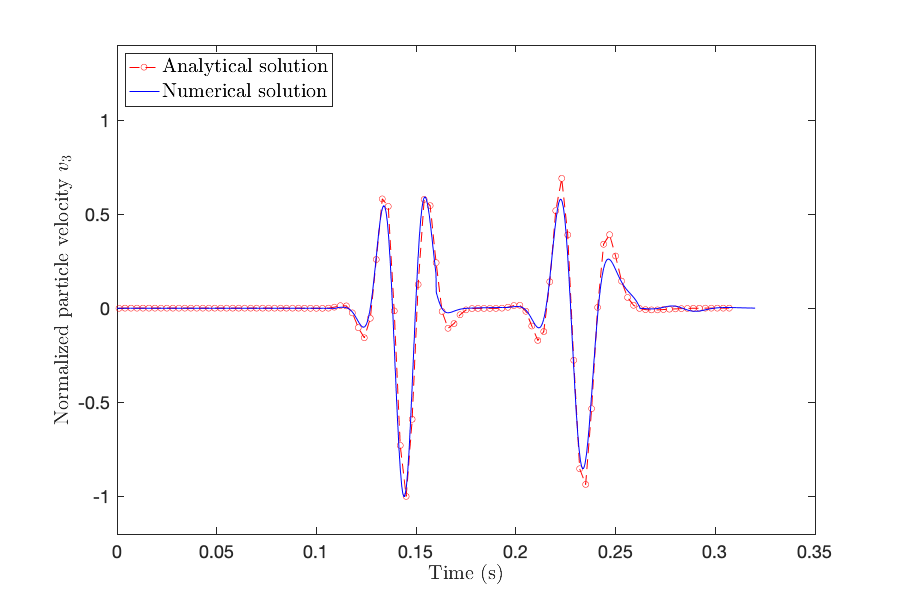}
	}
	\caption{A comparison between numerical and analytical solution of viscoelastic  wave equation in a homogeneous media. Subfigures (a) and (b) represent the horizontal and vertical particle velocities, respectively. Numerical solution is computed in 2D and for polynomials of degree $N=3$.}
	
\end{figure}

\subsubsection{2D Phenolic  material}
To further validate our numerical scheme, we perform a computational experiment in Phenolic material which has a high relaxation frequency of 250~kHz with  material properties given in Table 1 (Column 2). The size of the computational domain is $40~\text{cm} \times 40~\text{cm} $. The domain is discretized with uniform triangular element with a minimum edge length of $0.1562~\text{cm}$. Figures 4(a)-(b) represent the $x-$ and $z-$ components of the  particle velocity of the Phenolic material, respectively. The central frequency of the forcing function is $f_0=250~\text{kHz}~\text{(the frequency for relaxation peak)}$. Polynomials of degree $N=3$ are used for the simulation. The propagation time is $53.2~\text{$\mu$s}$. Both modes of waves can be observed: the  P mode  and the shear mode (S, inner wavefront). 
\begin{figure}
	\centering
	\subfloat[ $v_1$]{
		\begin{tikzpicture}
		\begin{axis}[enlargelimits=false, axis on top, axis equal image, xlabel=x (km), ylabel=y (km), width=0.6\textwidth]
		
		\addplot graphics [xmin=0,xmax=3,ymin=0,ymax=3] {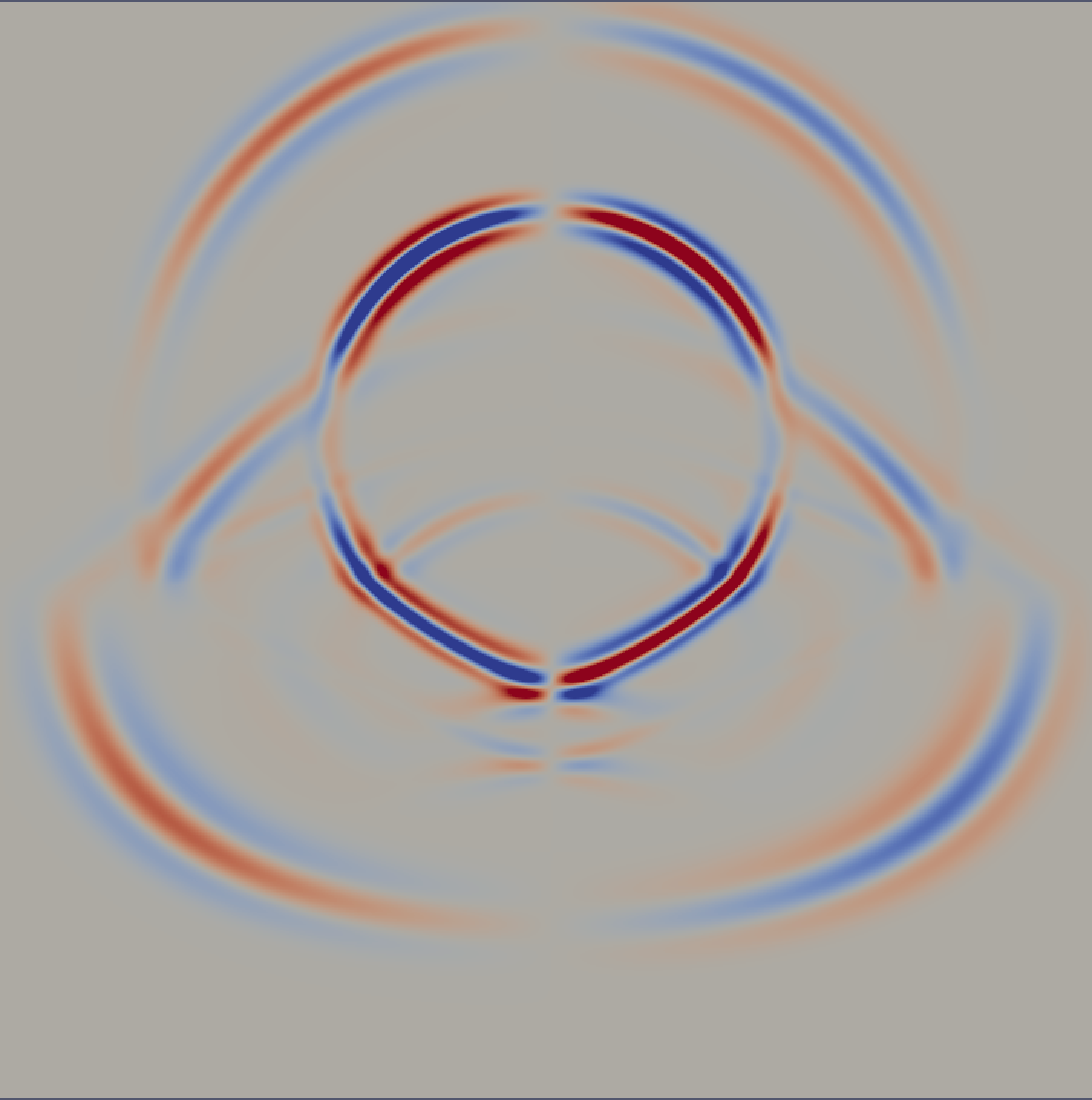};
		\draw[black, thin] (0, 140) -- (300, 140);
		\end{axis}
		\end{tikzpicture}
	}
	\subfloat[$v_2$]{
		\begin{tikzpicture}
		\begin{axis}[enlargelimits=false, axis on top, axis equal image, xlabel=x (km), ylabel=y (km), width=0.6\textwidth]		
		\addplot graphics [xmin=0,xmax=3, ymin=0,ymax=3] {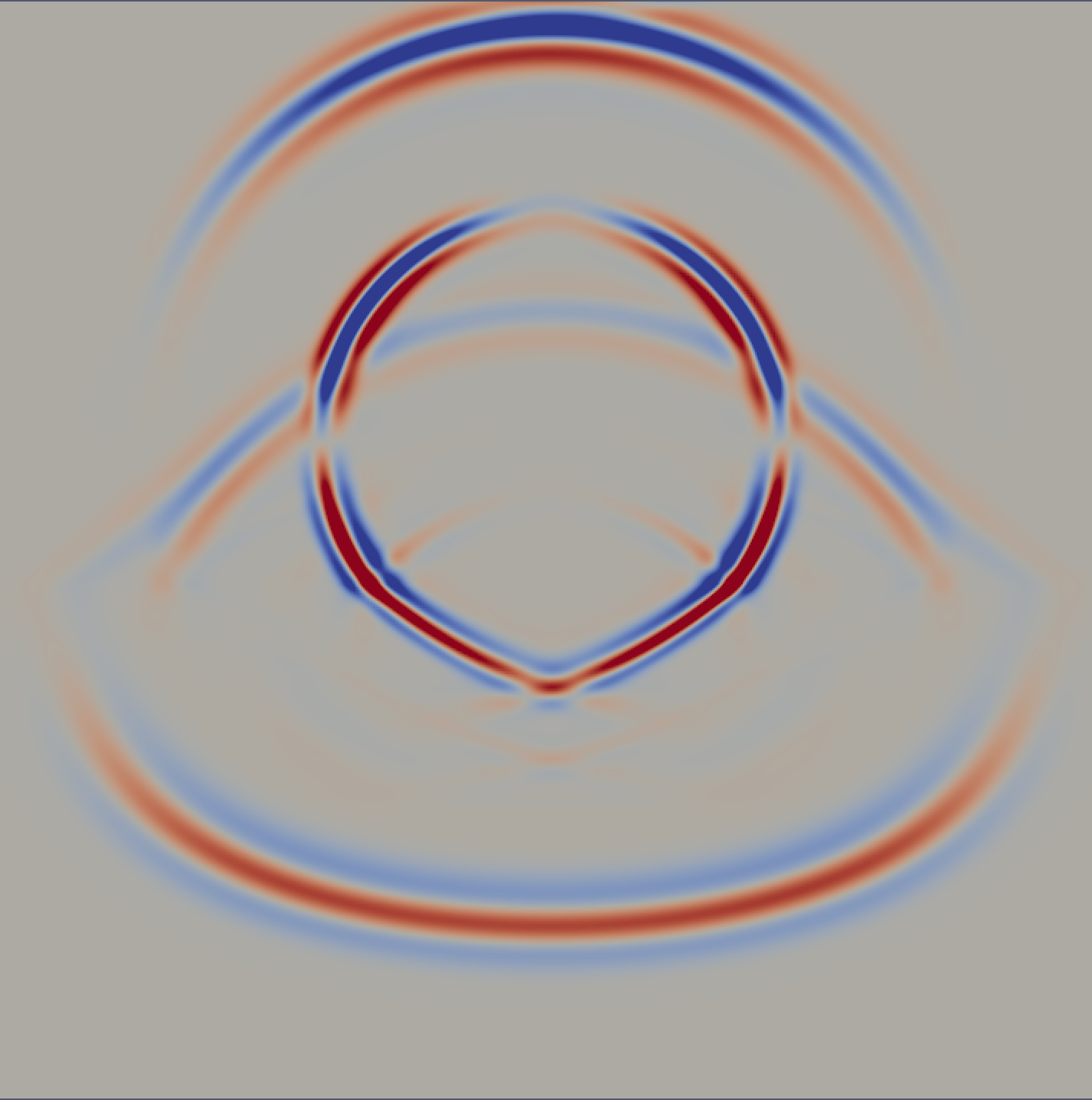};
		\draw[black, thin] (0, 140) -- (300, 140);
		\end{axis}
		\end{tikzpicture}
	}
	\caption{Snapshots of particle velocities in the layered model with (a) and (b) showing $v_1$ and $v_2$ components at $t=0.4~\text{s}$. The central frequency of the forcing function is 20 Hz. The point source is located at $(1.5~\text{km}, 1.8~\text{km})$.The solution is computed using polynomials of degree $N=3$ and $K=32,768$. }
\end{figure}

\subsubsection{Comparison of elastic and viscoelastic models}
To show the effect of the attenuation on wave propagation, we compare numerical solutions of elastic and viscoelastic wave equation in an isotropic sandstone with material properties given in Table 1 (Column 3). Numerical solution are computed in a domain of dimension  $[-1~\text{km}, 1~\text{km}]\times [-1~\text{km}, 1~\text{km}]$, and discretized with uniform triangular elements with a minimum edge length of $20.833~\text{m}$. Figure 5 shows a comparison between the numerical solutions of particle velocities for elastic and viscoelastic equation with $x-$ and $z-$ components represented in Figure 5a and 5b, respectively. The central frequency of the forcing function is $f_0=20~\text{Hz}$. Polynomials of degree $N=3$ are used for the simulation. The solution are stored at receiver position $(250~\text{m}, 250~\text{m})$ with source located at $(0~\text{m},0~\text{m})$. A difference between the amplitude of  elastic and viscoelastic solutions in Figure 5 is due to the attenuation brought in to the system due to the relaxation. The relaxation also results into decreasing the velocity of waves (when compared against the pure elastic or lossless case), which is clearly reflected by the difference in phases between the elastic and viscoelastic solutions, shown in Figure 5.

\subsubsection{Comparisons of analytical and numerical solutions}
Now, we compare the analytical and numerical solution, computed from our DG method, of 2D viscoelastic wave equation. (\ref{eq24}). The analytical solution of isotropic viscoelastic wave equation is computed by Carcione \cite{carcione2014} using correspondence principle \cite{blandtheory}. The derivation of analytical solution in a homogeneous and isotropic medium is given in \ref{A2}. The following forcing function is used to compute the analytical solution 
\begin{align} \label{cfst}
f(\bm{x}, t)=\exp{\left[- \myfrac{\Delta \omega^2 (t-t_0)^2}{4}\right]}\cos[\bar{\omega}(t-t_0)]\delta(\bm{x} - \bm{x}_0),
\end{align}
where $\bar{\omega}=2\pi f_0$ is central angular frequency with $\Delta \omega = \myfrac{\bar{\omega}}{2}$. \\
To compute the analytical solution one required the frequency spectrum of (\ref{cfst}), which is expressed as
\begin{align} \label{cfsf}
F(\omega)=\myfrac{\sqrt{\pi}}{\Delta\omega}\left(\exp\left[-\left(\myfrac{\omega + \bar{\omega}}{\Delta \omega}\right)\right] + \exp\left[-\left(\myfrac{\omega - \bar{\omega}}{\Delta \omega} \right) \right] \right)\exp(-\text{i}\omega t_0).
\end{align}
Equation (\ref{cfsf}) also satisfies the condition $F(\bar{\omega} + \Delta \omega)=F(\bar{\omega})/e$, which is the requirement to compute the analytical solution.
Figures 6 shows a comparison between time histories of numerical and analytical solutions of viscoelastic wave equation with $x-$ and $z-$ components being represented in Figure 6a and 6b, respectively. Numerical solution are computed in a domain of dimension  $[-0.5~\text{km}, 0.5~\text{km}]\times [-0.5~\text{km}, 0.5~\text{km}]$, which is discretized with uniform triangular elements with a minimum edge length of $4~\text{m}$. The material properties of isotropic sandstone (Table 1, Column 3) is used to compute the solutions. The central frequency $f_0$ of forcing function is $45~\text{Hz}$, which is located at $(0~\text{m}, 0~\text{m})$. The forcing function is added to the force corresponding to $\sigma_{22} $. Polynomials of degree $N=3$ are used for the simulation. The solution is stored at the node with coordinate $(250~\text{m}, 250~\text{m})$. Figure 6a and 6b show a very good agreement between and analytical and numerical solutions and thus validating the accuracy of the proposed numerical method.

\begin{figure}
	\centering
	\subfloat[$v_1$]{
		\includegraphics[trim={9cm, 0cm, 0cm, 0cm}, clip, width=0.5\textwidth]{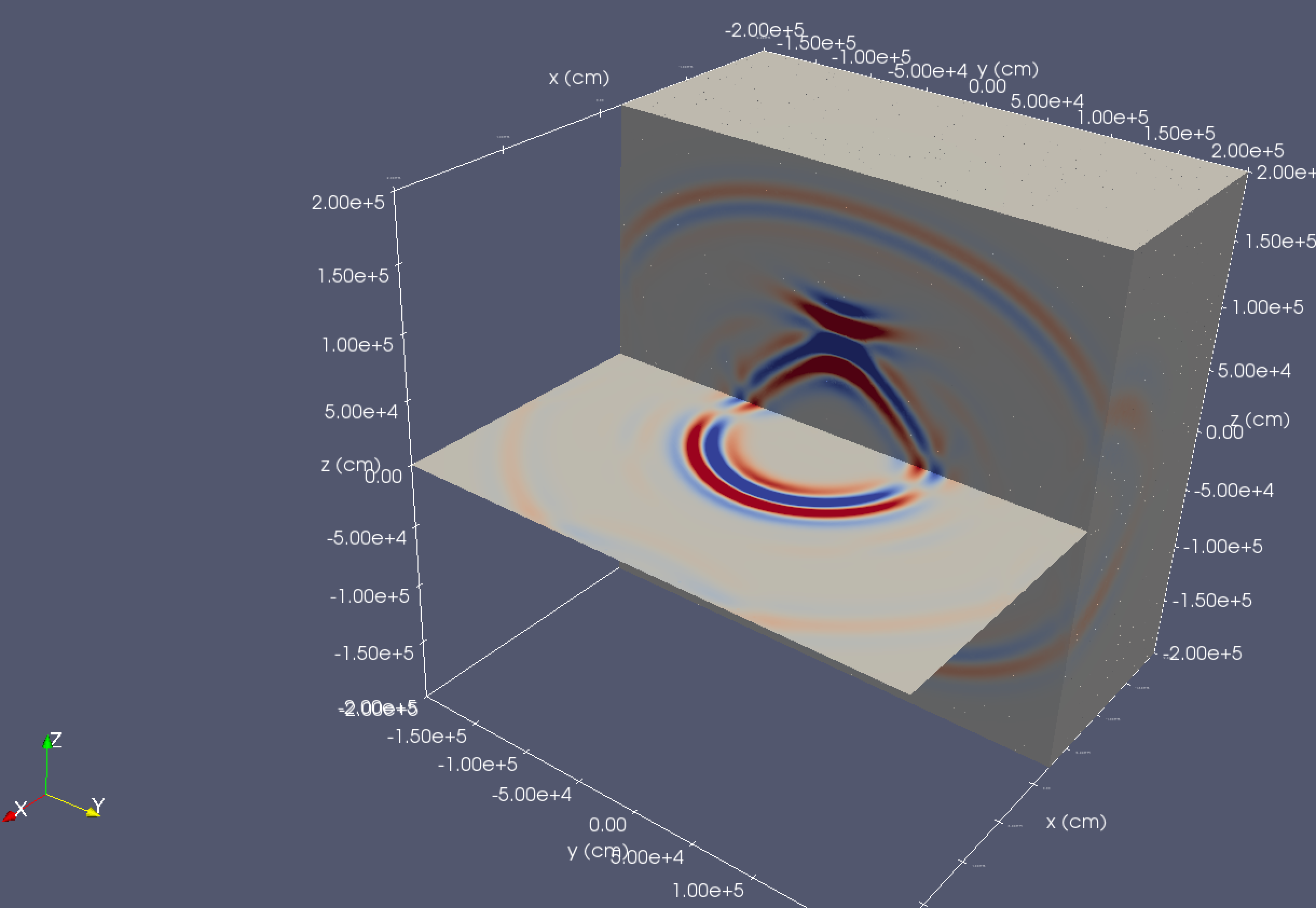}
	}
	\subfloat[$v_2$]{
		\includegraphics[trim={9cm, 0cm, 0cm, 0cm}, clip, width=0.5\textwidth]{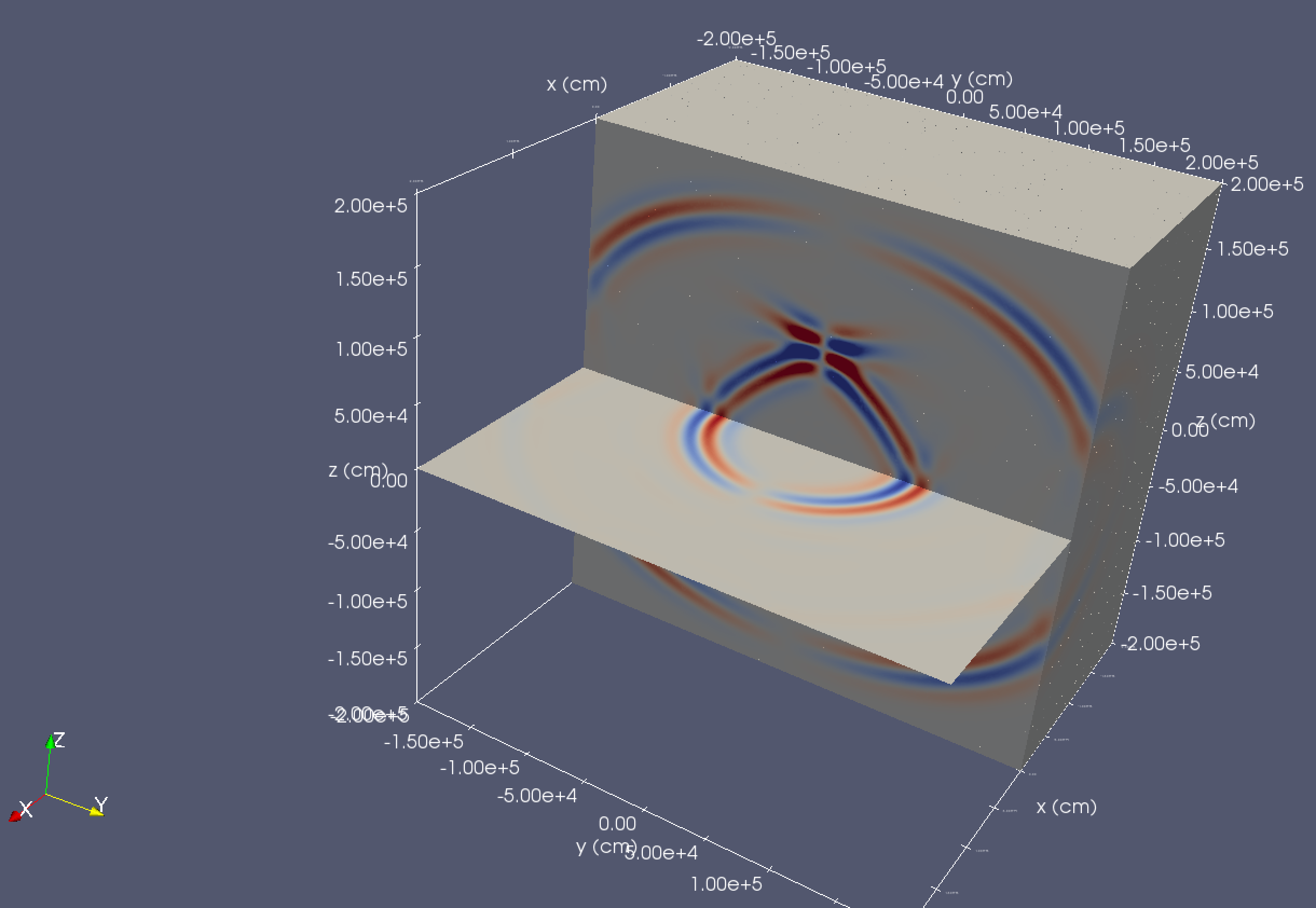}
	}\\
	\subfloat[$v_3$]{
		\includegraphics[trim={9cm, 0cm, 0cm, 0cm}, clip, width=0.5\textwidth]{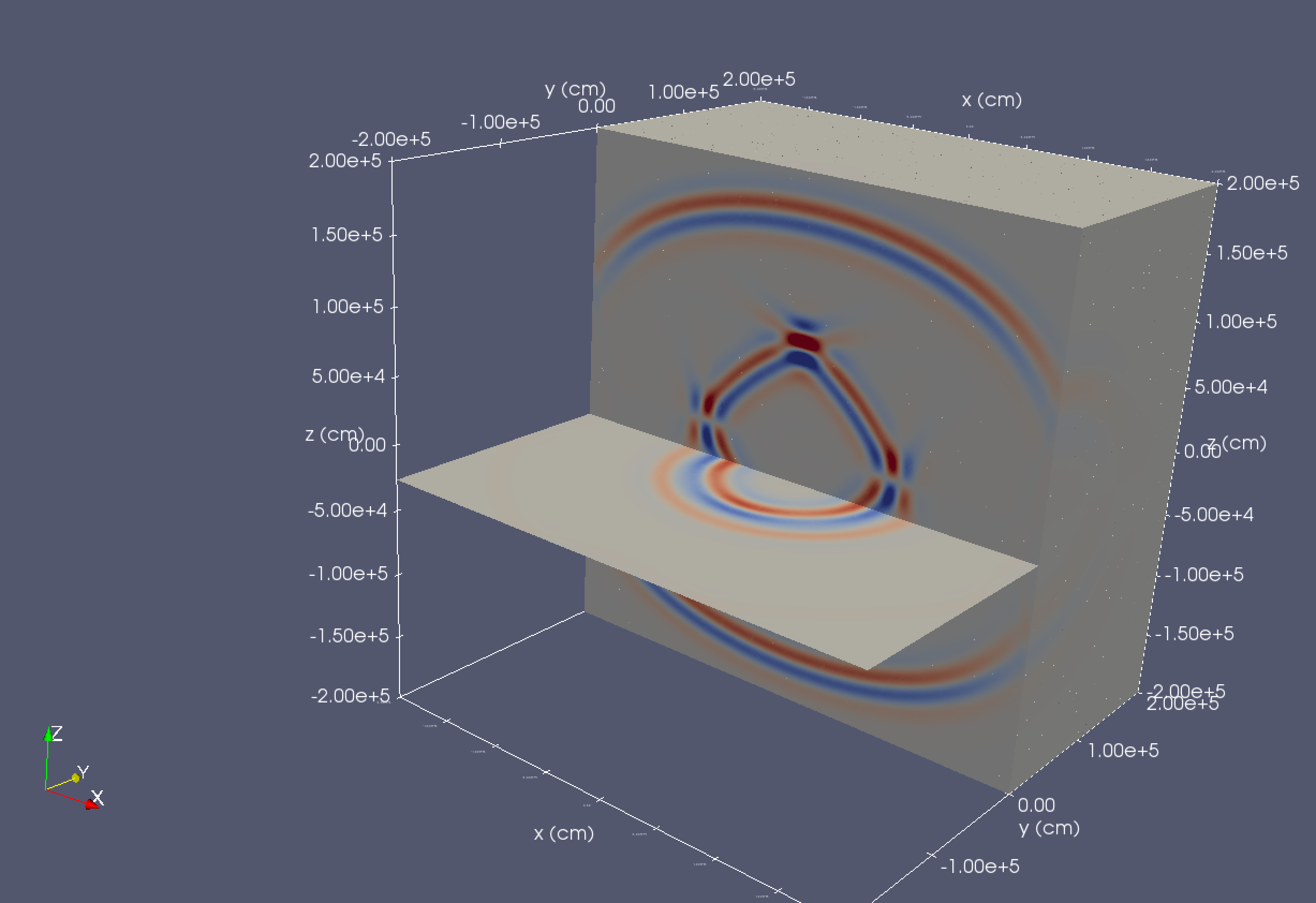}
	}
	
	\caption{Snapshots of particle velocities in a 3D homogeneous Orthotropic shale  with material properties shown in columns 1 of Table 1. Sub figures (a), (b), and (c) are showing $v_1,~v_2$ and $v_3$ components at $t=0.48~\text{s}$. The central frequency of the forcing function is 20 Hz. The point source is located at the center of the domain. The solution is computed using polynomials of degree $N=3$ and $h=32.5~\text{m}$.}
\end{figure}

\begin{figure}
	\centering
	\subfloat[$v_2$]{
		\includegraphics[trim={0cm, 0cm, 0cm, 0cm}, clip, width=0.5\textwidth]{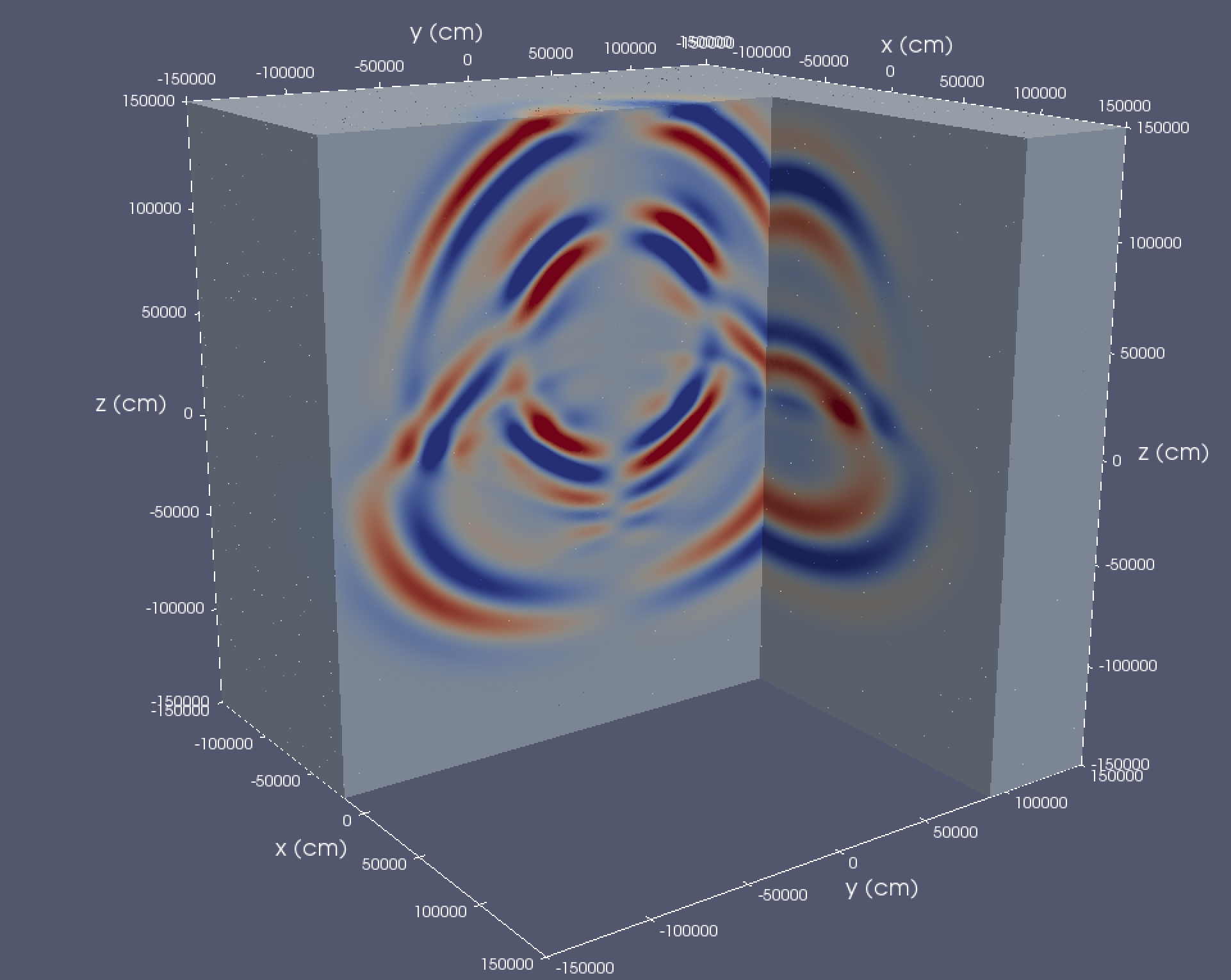}
	}
	\subfloat[$v_3$]{
		\includegraphics[trim={0cm, 0cm, 0cm, 0cm}, clip, width=0.5\textwidth]{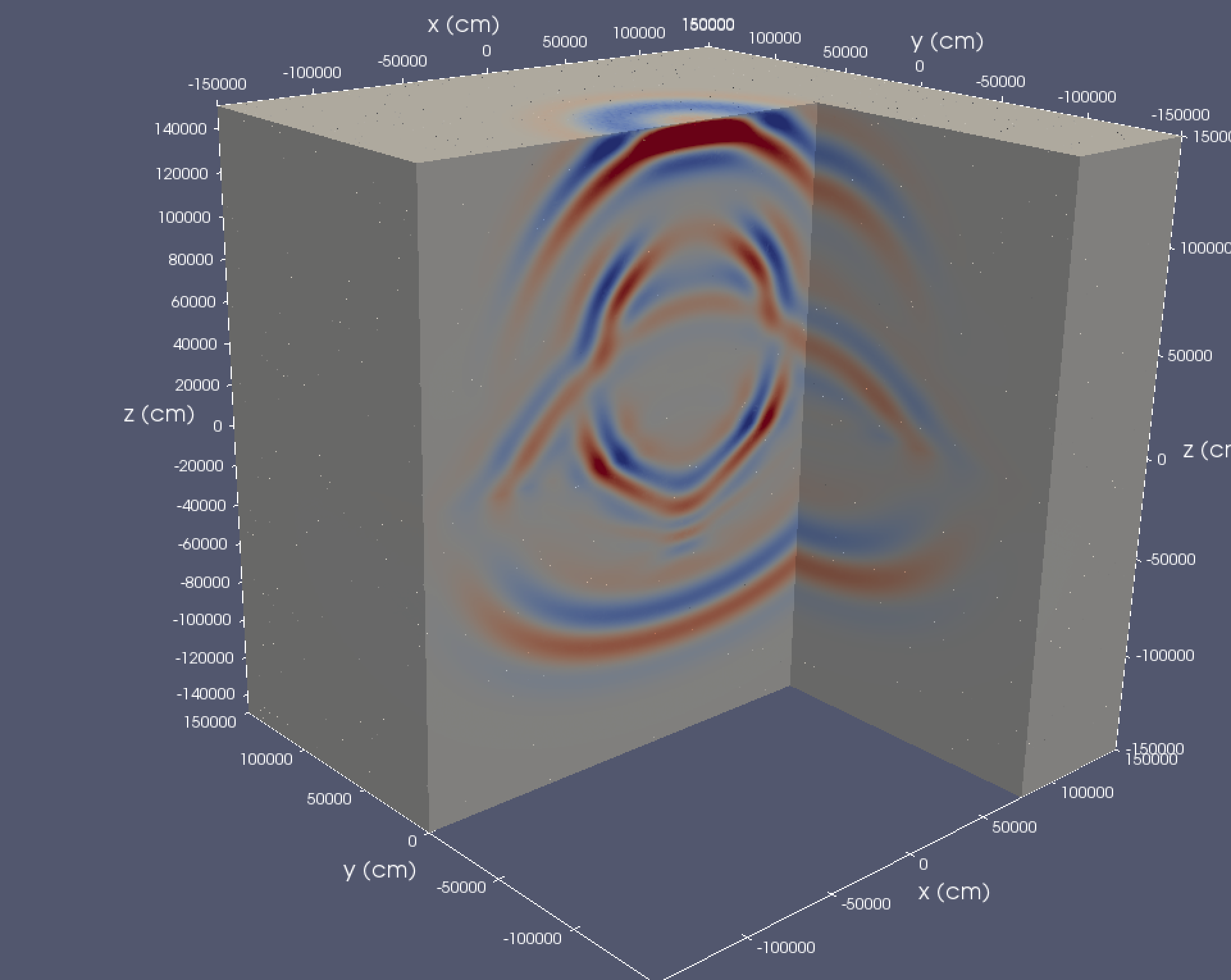}
	}
	\caption{Snapshots of particle velocities in a 3D anisotropic heterogeneous two layer model of size $3~\text{km}\times3~\text{km}\times3~\text{km}$. Sub figures (a), and (b) are represnt  $v_2$ and $v_3$ components at $t=0.4~\text{s}$. The central frequency of the forcing function is 20 Hz. The point source is located in the domain at $(0~\text{km}, 0~\text{km}, 0.3~\text{km})$. The solution is computed using polynomials of degree $N=3$ and $h=32.5~\text{m}$.}
\end{figure}

\subsubsection{2D Isotropic-anisotropic layered model}
In this example, we illustrate the effect of an interface between two layers of viscoelastic media. 
In the layered model, the top and bottom layers correspond to isotropic sandstone and orthotropic shale, with material properties are given in Table 1 (Column 3 and Column 1). The size of the computational domain is $3~\text{km} \times 3~\text{km}$ in the $x$ and $z$ directions, respectively. The minimum edge size of the triangular elements used to mesh the domain is $23.4375~\text{m}$. The point source is located at $(1.5~\text{km}, 1.8~\text{km})$ with a Ricker wavelet of frequency 20 Hz. The propagation time is $0.4~\text{s}$. The simulation is performed using polynomials of degree $N=3$. Snapshots of the $x$ and $z$ components of the particle velocity are shown in Figures 7a and 7b, respectively. Figure 7 clearly shows the direct, reflected, and transmitted wavefronts, corresponding to both P and S wave modes.  The effect of anisotropy on all three modes is clearly seen as wavefronts move with different phase velocities.

\subsubsection{3D Orthotropic  material}
Now, we  validate our numerical scheme for a 3D model. First, we perform a 3D computational experiment for orthotropic shale   with the material properties given in Table 1 (Column 1). The size of the computational domain is $4~\text{km} \times 4~\text{km} \times 4~\text{km} $ and is discretized by tetrahedral element with a minimum edge length of $32.5~\text{m}$. The central frequency of the forcing function is $f_0=20~\text{Hz}~\text{(the frequency for relaxation peak)}$. Polynomials of degree $N=3$ are used for the simulation. The propagation time is $0.48~\text{$\mu$s}$.  Figures 8(a), (b) and (c) represent the $x-,~y-$ and $z-$ components of the  particle velocity of the orthotropic material, respectively.  Both modes of waves can be observed: the  P mode  and the shear mode (S, inner wavefront). Figure 8 shows that the plane perpendicular to $z-$ direction is  a plane of symmetry, as wave propagation in the plane is isotropic. 

\subsubsection{3D Isotropic-anisotropic layered model}
In this example, we illustrate the effect of a two dimensional interface between two layers of the 3D viscoelastic media. The top and bottom layer of 3D model are comprised of isotropic sandstone and orthotropic shale, respectively.  The size of the computational domain is $3~\text{km} \times 3~\text{km} \time 3~\text{km}$ in the $x, ~y$ and $z$ directions, respectively. The minimum edge size of the triangular elements used to mesh the domain is $32.5~\text{m}$. The point source is located at $(0~\text{km}, 0~\text{km}~0.4~\text{km})$ with a Ricker wavelet of frequency 20 Hz. The propagation time is $0.4~\text{s}$. The simulation is performed using polynomials of degree $N=3$. Snapshots of the $y$ and $z$ components of the particle velocity are shown in Figures 9a and 9b, respectively. Figure  9 clearly demonstrate the effect the interface responsible for the direct, reflected, and transmitted wavefronts of P and shear waves present in the system.

\subsection{A large 3D heterogeneous subsurface model}
We use a 3D reservoir model from Shukla \etal \cite{Shukla2020}. The model is characterized by rock layers, discontinuity, and a surface with undulated topography. The discretized model is shown in Figure 10a. The dimension of the model is  $(22.8\times \text{km}\times 17.4~\text{km} \times 8.0~\text{km})$ in x, y and z directions, respectively. The domain is discretized with tetrahedral elements with a minimum edge length of $125~\text{m}$. The top surface of the model is perturbed so that the effects of the topography, assumed as a free surface, could be incorporated into numerical simulations. Figure 10(b) represent the z- component of the particle velocity at $3.5~\text{s}$. The central frequency of the forcing function is $20~\text{Hz}$. Polynomials of degree $N=3$  are used for simulation. The various modes of transmissions, reflections and scattering can be clearly seen in Figure 10b.

\begin{figure}
	\centering
	\subfloat[3D model constructed with topography with minimum edge length of element $h=125~\text{m}$]{
		\includegraphics[width=0.9\textwidth]{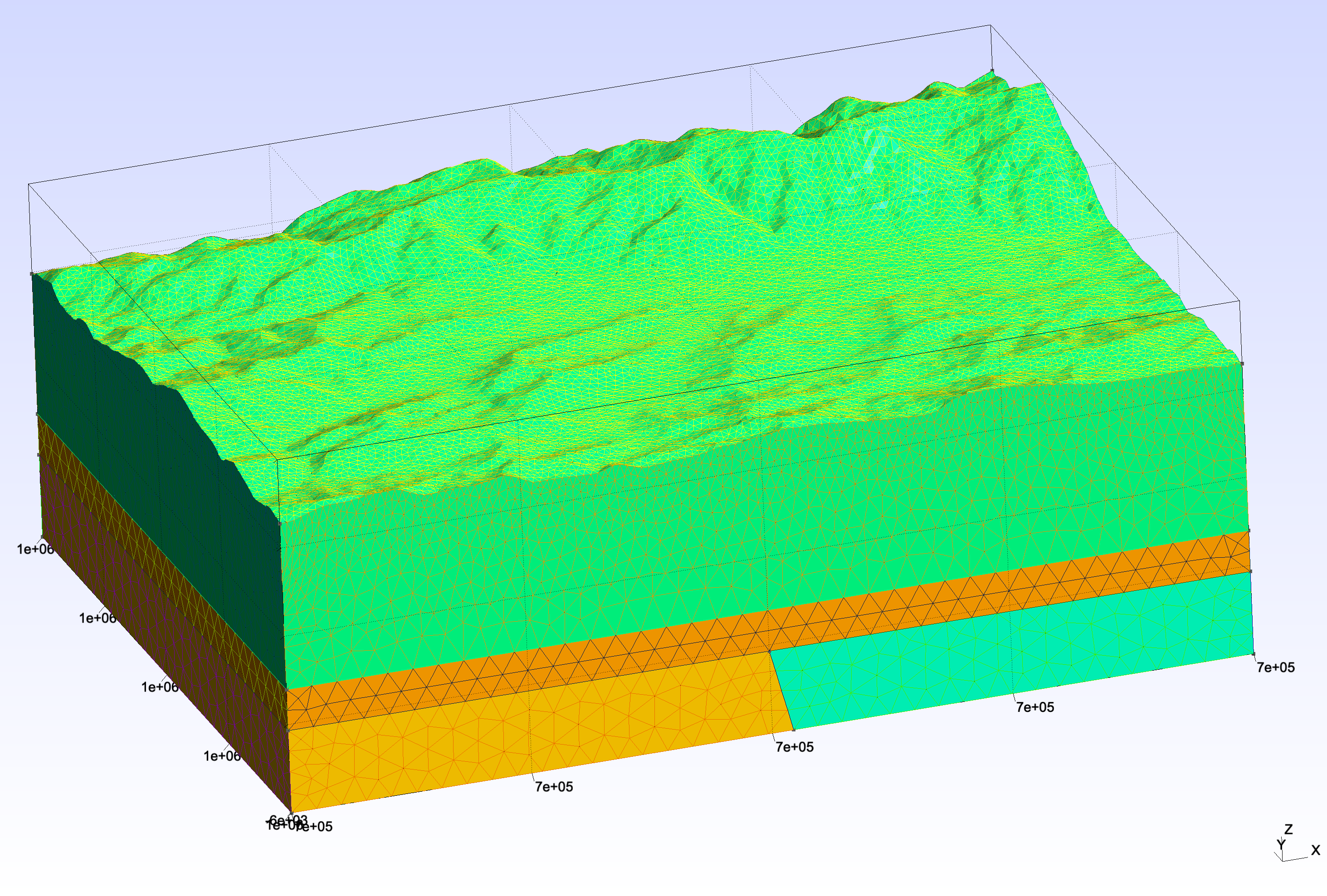}
	}\\
	\subfloat[Snapshot of $v_3$ at $t=3.5~\text{s}$]{
		\includegraphics[width=0.9\textwidth]{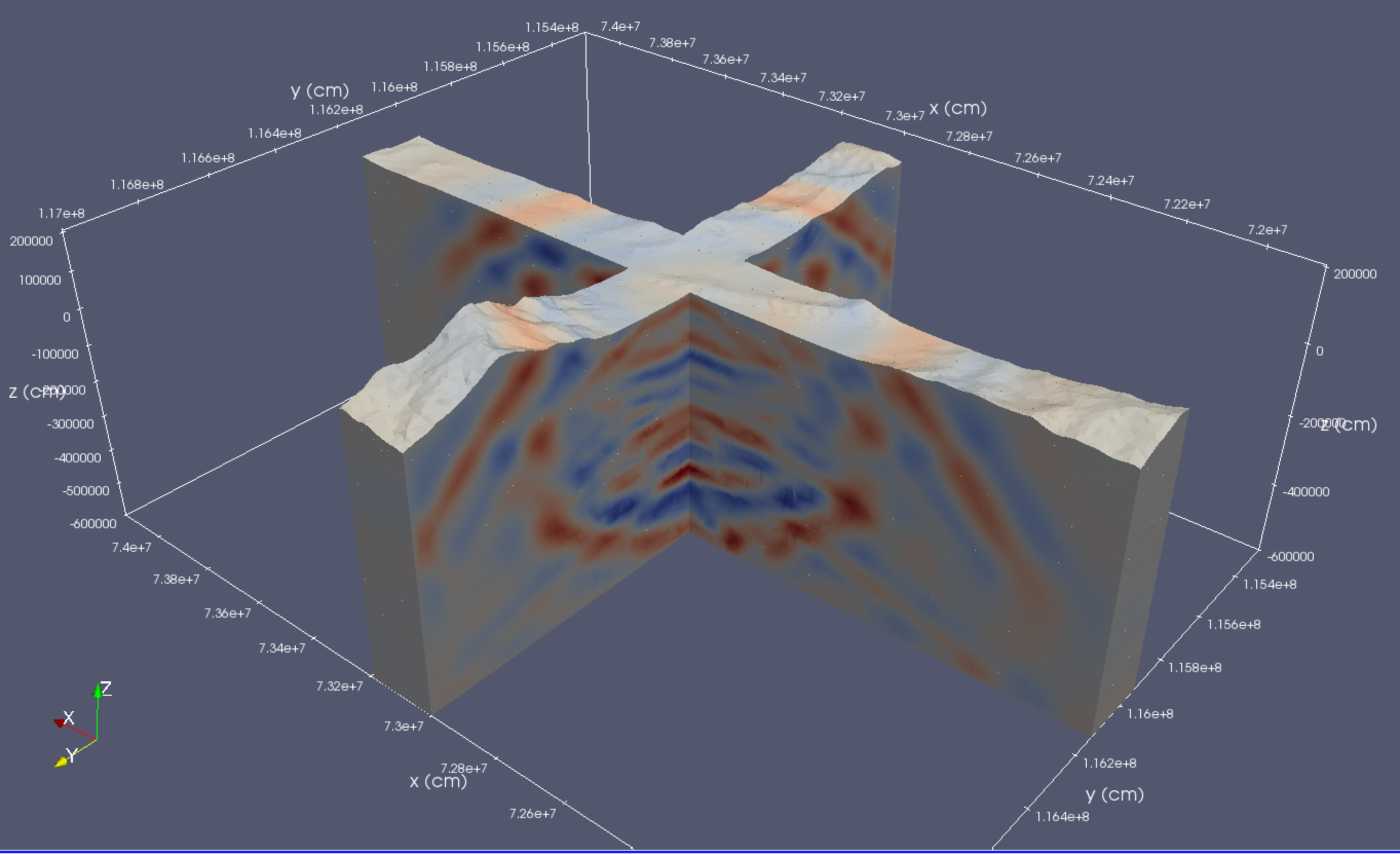}
	}
	\caption{Snapshot of $v_3$ for a 3D heterogeneous model constructed with topography on the top surface and a discontinuity in lower most layer. Subfigure (a) shows the 3D model of domain size $(22.8~\text{km},~17.4~\text{km}~8.0~\text{km})$ and discretized with tetrahedral element with a minimum edge length of element being $h=125~\text{m}$. Subfigure (b) shows the snapshot of $v_3$ at $3.5~\text{s}$. The point source is located in the domain at $(11.4~\text{km}, -8.7~\text{km},-50~\text{m})$ . The solution is computed using polynomials of degree $N=3$.}
\end{figure}

\section{Conclusions}
This work presents a high order discontinuous Galerkin method for a new symmetric form of the linear anisotropic viscoelastic wave equations. The method is energy stable and high order accurate for arbitrary stiffness tensors. We confirm the high-order accuracy of the numerical method using an analytic plane wave solution in a viscoelastic media. Finally, we provide computational results for various combinations of homogeneous and heterogeneous medium.

\section{Acknowledgments}
The authors gratefully thank the sponsors of the Geo-Mathematical Imaging Group at Rice University for providing the resources to carry out this work. Jesse Chan gratefully acknowledges support from the NSF
under awards DMS-1719818 and DMS-1712639. MVdH gratefully acknowledges support from the Simons Foundation under the MATH + X program and the NSF under grant DMS-1815143. The authors gratefully acknowledge Dr. Jos{\'e}  M Carcione of INOGS Italy, for his help in deriving and implementing the analytical solution in a 2D viscoelastic medium. 

\bibliographystyle{unsrt}
\bibliography{dg}

\begin{thebibliography}{10}

\bibitem{hosten1987}
B~Hosten, M~Deschamps, and Bernhard~R Tittmann.
\newblock Inhomogeneous wave generation and propagation in lossy anisotropic
  solids. application to the characterization of viscoelastic composite
  materials.
\newblock {\em The Journal of the Acoustical Society of America},
  82(5):1763--1770, 1987.

\bibitem{arts1992}
Rob~J Arts and Patrick~NJ Rasolofosaon.
\newblock Approximation of velocity and attenuation in general anisotropic
  rocks.
\newblock In {\em SEG Technical Program Expanded Abstracts 1992}, pages
  640--643. Society of Exploration Geophysicists, 1992.

\bibitem{carcione2014}
JM~Carcione.
\newblock {\em Wave Fields in Real Media: Wave Propagation in Anisotropic,
  Anelastic, Porous and Electromagnetic Media}.
\newblock Elsevier Science, 2014.

\bibitem{carcione1990}
JM~Carcione.
\newblock Wave propagation in anisotropic linear viscoelastic media: {T}heory
  and simulated wavefields.
\newblock {\em Geophysical Journal International}, 101(3):739--750, 1990.

\bibitem{mehrabadi1990}
Morteza~M Mehrabadi and Stephen~C Cowin.
\newblock Eigentensors of linear anisotropic elastic materials.
\newblock {\em The Quarterly Journal of Mechanics and Applied Mathematics},
  43(1):15--41, 1990.

\bibitem{helbig1996}
K~Helbig.
\newblock Foundations of anisotropy for exploration seismics.
\newblock In {\em International Journal of Rock Mechanics and Mining Sciences
  and Geomechanics Abstracts}, volume~1, pages 19A--20A, 1996.

\bibitem{carcione1995}
Jose~M Carcione.
\newblock Constitutive model and wave equations for linear, viscoelastic,
  anisotropic media.
\newblock {\em Geophysics}, 60(2):537--548, 1995.

\bibitem{moczo2014}
Peter Moczo, Jozef Kristek, and Martin G{\'a}lis.
\newblock {\em The finite-difference modelling of earthquake motions: Waves and
  ruptures}.
\newblock Cambridge University Press, 2014.

\bibitem{igel2017}
Heiner Igel.
\newblock {\em Computational {S}eismology: {A} {P}ractical {I}ntroduction}.
\newblock Oxford University Press, 2017.

\bibitem{virieux1986}
Jean Virieux.
\newblock P-{SV} wave propagation in heterogeneous media: Velocity-stress
  finite-difference method.
\newblock {\em Geophysics}, 51(4):889--901, 1986.

\bibitem{marfurt1984}
Kurt~J Marfurt.
\newblock Accuracy of finite-difference and finite-element modeling of the
  scalar and elastic wave equations.
\newblock {\em Geophysics}, 49(5):533--549, 1984.

\bibitem{moczo2007}
Peter Moczo, Johan~OA Robertsson, and Leo Eisner.
\newblock The finite-difference time-domain method for modeling of seismic wave
  propagation.
\newblock {\em Advances in {G}eophysics}, 48:421--516, 2007.

\bibitem{bohlen2002}
Thomas Bohlen.
\newblock Parallel {3-D} viscoelastic finite difference seismic modelling.
\newblock {\em Computers \& Geosciences}, 28(8):887--899, 2002.

\bibitem{levander1988}
Alan~R Levander.
\newblock Fourth-order finite-difference {P-SV} seismograms.
\newblock {\em Geophysics}, 53(11):1425--1436, 1988.

\bibitem{strikwerda2004}
John~C Strikwerda.
\newblock {\em Finite difference schemes and partial differential equations},
  volume~88.
\newblock Siam, 2004.

\bibitem{drainville2019}
Robert~Andrew Drainville, Laura Curiel, and Samuel Pichardo.
\newblock Superposition method for modelling boundaries between media in
  viscoelastic finite difference time domain simulations.
\newblock {\em The Journal of the Acoustical Society of America},
  146(6):4382--4401, 2019.

\bibitem{tessmer}
Ekkehart Tessmer and Dan Kosloff.
\newblock {3-D} elastic modeling with surface topography by a {C}hebychev
  spectral method.
\newblock {\em Geophysics}, 59(3):464--473, 1994.

\bibitem{carcione1991}
Jos{\'e}~M Carcione.
\newblock Domain decomposition for wave propagation problems.
\newblock {\em Journal of Scientific Computing}, 6(4):453--472, 1991.

\bibitem{bao1998}
Hesheng Bao, Jacobo Bielak, Omar Ghattas, Loukas~F Kallivokas, David~R
  O'Hallaron, Jonathan~R Shewchuk, and Jifeng Xu.
\newblock Large-scale simulation of elastic wave propagation in heterogeneous
  media on parallel computers.
\newblock {\em Computer {M}ethods in {A}pplied {M}echanics and {E}ngineering},
  152(1-2):85--102, 1998.

\bibitem{patera1984}
Anthony~T Patera.
\newblock A spectral element method for fluid dynamics: laminar flow in a
  channel expansion.
\newblock {\em Journal of Computational Physics}, 54(3):468--488, 1984.

\bibitem{seriani1994spectral}
G~za Seriani and Enrico Priolo.
\newblock Spectral element method for acoustic wave simulation in heterogeneous
  media.
\newblock {\em Finite {E}lements in {A}nalysis and {D}esign}, 16(3):337--348,
  1994.

\bibitem{komatitsch1998}
Dimitri Komatitsch and Jean-Pierre Vilotte.
\newblock The spectral element method: An efficient tool to simulate the
  seismic response of {2D} and {3D} geological structures.
\newblock {\em Bulletin of the Seismological Society of America},
  88(2):368--392, 1998.

\bibitem{carcione2002}
Jose~M Carcione, G{\'e}rard~C Herman, and APE Ten~Kroode.
\newblock Seismic {M}odeling.
\newblock {\em Geophysics}, 67(4):1304--1325, 2002.

\bibitem{hesthaven2007}
J.S. Hesthaven and T.~Warburton.
\newblock {\em Nodal discontinuous {G}alerkin methods: algorithms, analysis,
  and applications}, volume~54.
\newblock Springer, 2007.

\bibitem{klockner2009}
A.~Kl{\"o}ckner, T.~Warburton, J.~Bridge, and J.S. Hesthaven.
\newblock Nodal discontinuous {G}alerkin methods on graphics processors.
\newblock {\em Journal of Computational Physics}, 228(21):7863--7882, 2009.

\bibitem{ainsworth2004}
M.~Ainsworth.
\newblock {Dispersive and dissipative behaviour of high order discontinuous
  {G}alerkin finite element methods}.
\newblock {\em Journal of Computational Physics}, 198(1):106--130, 2004.

\bibitem{wilcox2010}
L.C. Wilcox, G.~Stadler, C.~Burstedde, and O.~Ghattas.
\newblock A high-order discontinuous {G}alerkin method for wave propagation
  through coupled elastic--acoustic media.
\newblock {\em Journal of Computational Physics}, 229(24):9373--9396, 2010.

\bibitem{kaser2007}
M.~K{\"a}ser, M.~Dumbser, J.~De~La~Puente, and H.~Igel.
\newblock {An arbitrary high-order discontinuous {G}alerkin method for elastic
  waves on unstructured meshes---III. Viscoelastic attenuation}.
\newblock {\em Geophysical Journal International}, 168(1):224--242, 2007.

\bibitem{de2007}
J.~de~la Puente, M.~K{\"a}ser, M.~Dumbser, and H.~Igel.
\newblock {An arbitrary high-order discontinuous {G}alerkin method for elastic
  waves on unstructured meshes---IV. Anisotropy}.
\newblock {\em Geophysical Journal International}, 169(3):1210--1228, 2007.

\bibitem{ye2016}
R.~Ye, M.V. de~Hoop, C.L. Petrovitch, L.J. Pyrak-Nolte, and L.C. Wilcox.
\newblock {A discontinuous {G}alerkin method with a modified penalty flux for
  the propagation and scattering of acousto-elastic waves}.
\newblock {\em Geophysical Journal International}, 205(2):1267--1289, 2016.

\bibitem{lambrecht2018}
L~Lambrecht, A~Lamert, W~Friederich, T~M{\"o}ller, and MS~Boxberg.
\newblock A nodal discontinuous {G}alerkin approach to 3-{D} viscoelastic wave
  propagation in complex geological media.
\newblock {\em Geophysical Journal International}, 212(3):1570--1587, 2018.

\bibitem{leitman1973}
Marshall~J Leitman and George~MC Fisher.
\newblock The linear theory of viscoelasticity (constitutive equations, creep
  laws, stress functions, variational principles and differential operators in
  dynamic and static linear viscoelasticity theory).
\newblock {\em Solid-state mechanics 3.(A 73-45495 24-32) Berlin,
  Springer-Verlag, 1973,}, pages 1--123, 1973.

\bibitem{chan2018}
J.~Chan.
\newblock Weight-adjusted discontinuous {G}alerkin methods: Matrix-valued
  weights and elastic wave propagation in heterogeneous media.
\newblock {\em International Journal for Numerical Methods in Engineering},
  113(12):1779--1809, 2018.

\bibitem{chan2017weight}
Jesse Chan, Russell~J Hewett, and Timothy Warburton.
\newblock Weight-adjusted discontinuous {G}alerkin methods: wave propagation in
  heterogeneous media.
\newblock {\em SIAM Journal on Scientific Computing}, 39(6):A2935--A2961, 2017.

\bibitem{Shukla2020}
Khemraj Shukla, Jesse Chan, V~Maarten, and Priyank Jaiswal.
\newblock A weight-adjusted discontinuous {G}alerkin method for the poroelastic
  wave equation: penalty fluxes and micro-heterogeneities.
\newblock {\em Journal of Computational Physics}, 403:109061, 2020.

\bibitem{berenger1994}
Jean-Pierre Berenger.
\newblock A perfectly matched layer for the absorption of electromagnetic
  waves.
\newblock {\em Journal of {C}omputational {P}hysics}, 114(2):185--200, 1994.

\bibitem{hagstrom2004}
Thomas Hagstrom and Timothy Warburton.
\newblock A new auxiliary variable formulation of high-order local radiation
  boundary conditions: corner compatibility conditions and extensions to
  first-order systems.
\newblock {\em Wave {M}otion}, 39(4):327--338, 2004.

\bibitem{chan2017gpu}
J.~Chan and T.~Warburton.
\newblock {{GPU}-Accelerated Bernstein--B\'{e}zier Discontinuous {G}alerkin
  Methods for Wave Problems}.
\newblock {\em SIAM Journal on Scientific Computing}, 39(2):A628--A654, 2017.

\bibitem{guo2018}
K.~Guo and J.~Chan.
\newblock {Bernstein-B\'{e}zier weight-adjusted discontinuous Galerkin methods
  for wave propagation in heterogeneous media}.
\newblock {\em arXiv preprint arXiv:1808.08645}, 2018.

\bibitem{carpenter1994fourth}
Mark~H Carpenter and Christopher~A Kennedy.
\newblock Fourth-order 2{N}-storage {R}unge-{K}utta schemes.
\newblock Technical Report NASA-TM-109112, NASA Langley Research Center, 1994.

\bibitem{chan2016gpu}
Jesse Chan, Zheng Wang, Axel Modave, Jean-Francois Remacle, and Tim Warburton.
\newblock {GPU}-accelerated discontinuous galerkin methods on hybrid meshes.
\newblock {\em Journal of Computational Physics}, 318:142--168, 2016.

\bibitem{chan2017penalty}
Jesse Chan and T~Warburton.
\newblock On the penalty stabilization mechanism for upwind discontinuous
  galerkin formulations of first order hyperbolic systems.
\newblock {\em Computers \& Mathematics with Applications}, 74(12):3099--3110,
  2017.

\bibitem{toro2009}
Eleuterio~F Toro.
\newblock The {HLL} and {HLLC} riemann solvers.
\newblock In {\em Riemann solvers and numerical methods for fluid dynamics},
  pages 315--344. Springer, 2009.

\bibitem{de2008}
J.~de~la Puente, M.~Dumbser, M.~K{\"a}ser, and H.~Igel.
\newblock Discontinuous {G}alerkin methods for wave propagation in poroelastic
  media.
\newblock {\em Geophysics}, 73(5):T77--T97, 2008.

\bibitem{blandtheory}
David~Russell Bland.
\newblock {\em The theory of linear viscoelasticity}.
\newblock Courier Dover Publications, 2016.

\bibitem{eason1956}
G~Eason, J~Fulton, and Ian~Naismith Sneddon.
\newblock The generation of waves in an infinite elastic solid by variable body
  forces.
\newblock {\em Philosophical Transactions of the Royal Society of London.
  Series A, Mathematical and Physical Sciences}, 248(955):575--607, 1956.

\end{thebibliography}

\appendix
\section{Inverse of compliance matrix $C$}
\label{A1}
The expressions for $r_{ij}$ in (\ref{RE1}) are
\begin{align}
r_{11}&= -\myfrac{( c_{11}c_{33} - c_{13}^2)}{(c_{11} - c_{12})(c_{11}c_{33} + c_{12}c_{33}- 2c_{13}^2 )}\\
r_{12}&= -\myfrac{(c_{12}c_{33}- c_{13}^2)}{(c_{11} - c_{12})(c_{11}c_{33} + c_{12}c_{33}- 2c_{13}^2 )}\\
r_{13}&= -\myfrac{c_{13}}{(c_{11}c_{33} + c_{12}c_{33}- 2c_{13}^2)} \\
r_{33}&= \myfrac{(c_{11} + c_{12})}{(c_{11}c_{33} + c_{12}c_{33}-2c_{13}^2)}
\end{align}
\section{Analytic solution in a  homogeneous viscoelastic media}
\label{A2}
The solution of the elastic wave equation in an 2-D isotropic  medium  for an impulsive point force  is given by Eason \etal \cite{eason1956}. For a force acting in the positive $x_3$-direction, displacement solutions are expressed as \cite{carcione2014}
\begin{align} \label{B1}
\begin{aligned}
u_1(r,t) &=\left(\myfrac{F_0}{2\pi\rho}\right)\myfrac{xz}{r^2}\left[G_1(r,t) + G_3(r,t) \right],\\
u_3(r,t) &=\left(\myfrac{F_0}{2\pi\rho}\right)\myfrac{1}{r^2}\left[z^2G_1(r,t) - x^2G_3(r,t) \right],
\end{aligned}
\end{align}
where $F_0$ is a the magnitude of the force and $r^2 = x^2 + z^2$, and $G_1(r, t)$ and $G_3(r,t)$ are Green's function expresses as
\begin{align} \label{B2}
\begin{aligned}
G_1(r,t) &= \myfrac{1}{c_p^2}(t^2 - \tau_p^2)^{-1/2} H(t-\tau_p) + \myfrac{1}{r^2}(t^2 - \tau_p^2)^{1/2} H(t-\tau_p) -\myfrac{1}{r^2}(t^2 - \tau_s^2)^{1/2} H(t - \tau_s), \\
G_3(r,t) &= -\myfrac{1}{c_s^2}(t^2 - \tau_s^2)^{-1/2} H(t-\tau_s) + \myfrac{1}{r^2}(t^2 - \tau_p^2)^{1/2} H(t-\tau_p) -\myfrac{1}{r^2}(t^2 - \tau_s^2)^{1/2} H(t - \tau_s),
\end{aligned}
\end{align}

where $\tau_p = \myfrac{r}{c_p},~\tau_s=\myfrac{r}{c_s}$ with $c_p$ and $c_s$ being phase velocities of the compressional and shear waves. $H(t)$ is Heaviside function. To recover the anelastic solution the correspondence principle \cite{carcione2014} is applied on frequency domain representation of (\ref{B2}). We also use following identities of transform pairs of zero- and first-order Hankel function of the second kind
\begin{align} \label{B3}
\begin{aligned}
\int_{-\infty}^{\infty} \myfrac{1}{\tau^2} \left(t^2 - \tau^2\right)^{1/2} H(t - \tau)\exp(i\omega t)dt &= \myfrac{\text{i}\pi}{2\omega\tau} H_1^{(2)}(\omega \tau), \\
\int_{-\infty}^{\infty} \myfrac{1}{\tau^2} \left(t^2 - \tau^2\right)^{-1/2} H(t - \tau)\exp(i\omega t)dt &= -\myfrac{\text{i}\pi}{2\omega\tau} H_0^{(2)}(\omega \tau).
\end{aligned}
\end{align}
Now, Fourier transform of (\ref{B2}) with respect to time yields

\begin{align}
\widehat{G}_1(r, \omega, c_p, c_s)&=-\myfrac{\text{i}\pi}{2}\left[\myfrac{1}{(c_p(\omega))^2} H_0^{(2)}\left(\myfrac{ \omega r}{c_p}\right) + \myfrac{1}{\omega r c_s(\omega)} H_1^{(2)}\left(\myfrac{\omega r}{c_s(\omega)} \right)  - \myfrac{1}{\omega r c_p(\omega)} H_1^{(2)}\left(\myfrac{\omega r}{c_p(\omega)} \right) \right], \label{B4} \\
\widehat{G}_3(r, \omega, c_p, c_s)&=\myfrac{\text{i}\pi}{2}\left[\myfrac{1}{(c_s(\omega))^2} H_0^{(2)}\left(\myfrac{ \omega r}{c_s(\omega)}\right) - \myfrac{1}{\omega r c_s(\omega)} H_1^{(2)}\left(\myfrac{\omega r}{c_s(\omega)} \right)  + \myfrac{1}{\omega r c_p(\omega)} H_1^{(2)}\left(\myfrac{\omega r}{c_p(\omega)} \right) \right], \label{B5} 
\end{align}
where $c_p(\omega)=\sqrt{\myfrac{(c_{11} +  c_{33}) M_1(\omega) + c_{33} M_2(\omega)}{\rho}}$, and $c_s(\omega)=\sqrt{\myfrac{c_{33}M_2(\omega)}{\rho}}$ where $M_1$ and $M_2$ are recovered from (\ref{eq12}), which for isotropic case are as follows,
\[
M_{\nu \in \{1,2\}} = \myfrac{\tau_\sigma^{(\nu)}}{\tau_\epsilon^{(\nu)}}\left(\myfrac{1 + \text{i}\omega \tau_{\epsilon}^{\nu}}{1 + \text{i}\omega \tau_{\sigma}^{\nu}} \right).
\] 
Now, taking the Fourier transform of (\ref{B1}) and using (\ref{B4}) and (\ref{B5}), we get
\begin{align} \label{b6}
\begin{aligned}
u_1(r,\omega,c_p, c_s) &=\left(\myfrac{F_0}{2\pi\rho}\right)\myfrac{xz}{r^2}\left[\widehat{G}_1(r,\omega, c_p, c_s) + \widehat{G}_3(r,\omega, c_p, c_s)\right]\\
u_3(r,\omega, c_p, c_s)&=\left(\myfrac{F_0}{2\pi\rho}\right)\myfrac{1}{r^2}\left[z^2\widehat{G}_1(r,\omega, c_p, c_s) - x^2 \widehat{G}_3(r,\omega, c_p, c_s)\right]
\end{aligned}
\end{align}
To ensure that solution is real in time, we express (\ref{b6}) as
\begin{align}\label{b7}
u_{1,3}(\omega) =\begin{cases}
u_{1,3}(r,\omega, c_p, c_s), \qquad &\omega \geq 0,\\
u_{1,3}^*(r,-\omega, c_p, c_s), \qquad &\omega < 0,
\end{cases}
\end{align}
where asterisk $(^*)$ denotes the complex conjugate.
Multiplying (\ref{b7}) with frequency domain representation of a source time function and then taking the inverse Fourier transform will yield time-domain analytical displacement solution of 2D viscoelastic wave equation  (\ref{eq24}). The $\widehat{G}_1$ and $\widehat{G}_3$ are considered as zeros due to Hankel's functions being singular.

\end{document}